\documentclass[a4paper]{amsart}

\pdfoutput=1

\usepackage{amsthm, amssymb, amsmath, amsfonts, mathrsfs}

\usepackage{microtype}

\usepackage[pagebackref,colorlinks=true,pdfpagemode=none,urlcolor=blue,
linkcolor=blue,citecolor=blue]{hyperref}

\usepackage{relsize}

\usepackage{amsmath,amsfonts,amssymb,amsthm}
\usepackage{amsthm, amssymb, amsmath, amsfonts, mathrsfs}

\usepackage{mathrsfs}
\usepackage{MnSymbol}
\usepackage{scalerel} 

\usepackage{color}
\usepackage{accents}

\usepackage[normalem]{ulem}
\usepackage{bbm}

\usepackage{cleveref}

\newtheorem{theorem}{Theorem}[section]
\newtheorem{lemma}[theorem]{Lemma}
\newtheorem{assumption}[theorem]{Assumption}

\newtheorem{corollary}[theorem]{Corollary}
\newtheorem{proposition}[theorem]{Proposition}

\theoremstyle{remark}
\newtheorem{remark}[theorem]{Remark}

\renewenvironment{proof}[1][Proof]{ {\itshape \noindent {#1.}} }{$\Box$
\medskip}

\numberwithin{equation}{section}
\newcommand{\R}{\mathbb{R}}

\newcommand{\Z}{\mathbb{Z}}

\newcommand{\Pb}{\mathbb{P}}
\newcommand{\PP}{\mathbf{P}}
\newcommand{\QQ}{\mathbf{Q}}
\newcommand{\E}{\mathbb{E}}

\newcommand{\F}{\mathcal{F}}

\newcommand{\M}{\mathbf{M}}

\newcommand{\Om}{\Omega}

\newcommand{\cR}{\mathcal{R}}

\newcommand{\U}{\mathcal{U}}

\newcommand{\eps}{\varepsilon}

\newcommand{\la}{\langle}
\newcommand{\ra}{\rangle}

\newcommand{\EE}{\mathbf{E}}
\newcommand{\1}{\mathbbm{1}}

\newcommand{\bbP}{\mathbb P}

\newcommand{\si}{\sigma}

\newcommand{\cP}{\mathcal{P}}

\newcommand{\bT}{\mathbb{T}}
\newcommand{\bbR}{\mathbb{R}}
\newcommand{\bbZ}{\mathbb{Z}}
\newcommand{\cal}{\mathcal}
\newcommand{\cZ}{\mathcal{Z}}
\newcommand{\om}{\omega}

\newcommand{\cD}{\mathcal{D}}
\newcommand{\cL}{\mathcal{L}}

\newcommand{\gu}{\textcolor{black}}

\newcommand{\tk}[1]{{\color{black}#1}}

\begin{document}

\title{KPZ on torus: Gaussian fluctuations}


\author{Yu Gu, Tomasz Komorowski}

\address[Yu Gu]{Department of Mathematics, University of Maryland, College Park, MD 20742, USA}

\address[Tomasz Komorowski]{Institute of Mathematics, Polish Academy of Sciences, ul.
Śniadeckich 8, 00-636 Warsaw, Poland}

\maketitle

\begin{abstract}

We study the KPZ equation on a torus and derive Gaussian fluctuations in large time.



\medskip


\end{abstract}

\section{Introduction}

\subsection{Main result}
\label{sec1.1}
We consider the stochastic heat equation (SHE) on the torus $\bT^d$:
\begin{equation}\label{e.she}
\partial_t {\cal U}=\tfrac12\Delta {\cal U}+{\cal U}\xi, 
\end{equation}
where $\xi$ is a generalized Gaussian random field over $\bbR\times\bT^d$ with the covariance function 
\[
\EE[\xi(t,x)\xi(s,y)]=\delta(t-s)R(x-y).
\]
Here $\bT^d$ is the unit torus.
We assume $\xi$ is built on the probability space $(\Omega,\F,\PP)$,
and consider two cases in the paper: (i) $R:\bT^d\to\R_{\geq0}$ is a
smooth bounded function and $d\geq 1$; (ii) $R(\cdot)=\delta(\cdot)$
is the Dirac function and $d=1$, i.e.,  $\xi$ is a $1+1$ spacetime
white noise. The product between ${\cal U}$ and $\xi$ is interpreted in the
It\^o's sense. Without loss of generality, we assume 
\begin{equation}
\label{hr0}
\int_{\bT^d} R(x)dx=1
\end{equation}
 in case (i).  We assume that the initial condition is an arbitrary probability measure $\nu$ on $\bT^d$.


%
Here is the main result of the paper:
\begin{theorem}\label{t.mainth}
There exist $\gamma,\sigma>0$, which only depend on $R,d$ and are defined in \eqref{e.defgamma} and \eqref{e.defsigma} below, such that for any $x\in\bT^d$, 
\[
\frac{\log \U(t,x)+\gamma t}{\sqrt{t}}\Rightarrow N(0,\sigma^2), \quad\quad \mbox{ as } t\to\infty.
\]
\end{theorem}
Results on process-level convergence can be found in Remark~\ref{r.process} below.

\subsection{Context}

Define $h(t,x)=\log \U(t,x)$, then $h$ formally satisfies the nonlinear SPDE
\begin{equation}\label{e.kpz}
\partial_t h=\tfrac12\Delta h+\tfrac12|\nabla h|^2+\xi,
\end{equation}
which is the Kardar-Parisi-Zhang (KPZ) equation and a default model for random surface growth \cite{kardar1986dynamic}. The study of the equation in $d=1$ with a spacetime white noise has witnessed important progress in recent years, including making sense of the equation \cite{gubinelli2015paracontrolled,gubinelli2017kpz,hairer2013solving,hairer2014theory,kupiainen2016renormalization}, studying universal behaviors  \cite{amir2011probability,balazs2011fluctuation,borodin,borodin1,duncan,kpzfix,quastelkpz,spohn,spohn1,virag}, etc. For a detailed review of singular SPDEs and the KPZ universality class, we refer to \cite{corwin2012kardar,haoshen,quastelintroduction,quastel2015one} and the references therein. For the equation \eqref{e.kpz} and other related  models in statistical physics including directed polymers and exclusion processes, the study of the large scale  behaviors is mostly focused on an infinite line and the limiting object is the so-called KPZ fixed point. It is expected that the   result on the torus is different and Gaussian fluctuations prevail in large time. Indeed, Theorem~\ref{t.mainth} shows that as $t\to\infty$,
\begin{equation}\label{e.conkpz}
\frac{h(t,x)+\gamma t}{\sqrt{t}}\Rightarrow N(0,\sigma^2).
\end{equation}
The result also holds in $d\geq 2$ with a smooth noise.

The recent work on periodic  totally asymmetric simple exclusion process \cite{zhipeng1,zhipeng2,zhipeng} considered the case where the size of the torus is scaled with time, and studied the fluctuations of the integrated particle current for several special initial data. Very precise limiting fluctuations were derived, and a transition from the Tracy-Widom distribution in the small time limit to the Gaussian distribution in the large time limit was derived in \cite[Theorem 1.5]{zhipeng}. Using $L$ to denote the size of the torus, in our problem, we have $L=1$, which corresponds to the super-relaxation time scale $t\gg L^{3/2}$. \gu{The relaxation time scale is $t\sim L^{3/2}$, the critical regime where $h(t,x)$ are correlated for all $x$ on the torus. In the sub-relaxation time scale $t\ll L^{3/2}$, the Tracy-Widom type of fluctuations are expected and proved when the equation is posed on the whole line. The Gaussian fluctuation in \eqref{e.conkpz} is due to the fast mixing in the time variable, which we will discuss in more detail later.}


There are also recent results on deriving the Edwards-Wilkinson limit of the KPZ equation in $d\geq2$ in a  weak disorder regime, with the equation posed on the whole space and driven by a smooth noise, see \cite{CRS18,chatterjee2018constructing,cometskpz,cosco,kpz1,kpz2,nikos,magnen2017diffusive}. Although the limit is also of Gaussian distribution, the mixing mechanism is rather different from our setting, and it is ultimately the result of the weak disorder that leads to the Gaussian fluctuations, which is related to the diffusive behaviors of directed polymers in high dimensions with high temperatures.

\tk{The problem we consider in this paper can be viewed as the KPZ equation on a fixed interval with a periodic boundary condition. One can ask the same question about the one with Neumann boundary conditions, which is related to the exclusion process on an interval with sources and sinks at the boundary points \cite{haoshen1}. We expect that our approach applies to this setting as well, see Remark~\ref{r.openkpz} below for more detail. Similar questions can be asked for the KPZ equation in half-space, where there is a phase transition in the parameter associated with the boundary condition. There has been a lot of progress in the study of such problems, including the construction of explicit invariant measures \cite{baco,barraquand2021steady,BKWW21,CK21}, the proof of Gaussian or non-Gaussian fluctuations in large time \cite{sasamoto,shalin}, etc. }

\tk{An interesting problem is to study the case when the size of the torus scales with the time variable, similar to the one considered for the exclusion process in \cite{zhipeng1,zhipeng2,zhipeng}. A key quantity to study is the variance $\sigma^2$ appeared in Theorem~\ref{t.mainth} and its asymptotic behavior as the size of the torus goes to infinity. Some partial progress can be found in \cite{dunlap}, where a different proof of Theorem~\ref{t.mainth} was also sketched, based on an application of the  Clark-Ocone formula, which leads to a rather explicit expression of $\sigma^2$ in terms of the invariant measure of the KPZ equation.}

Our proof of \eqref{e.conkpz} goes through the Hopf-Cole transformation so that we directly study $\U=\exp(h)$   and avoid dealing with the small scale singularity appearing in \eqref{e.kpz}. The recent work \cite{perkowski} studied the generators of the stochastic Burgers equation, the singular SPDE satisfied by $\nabla h(t,x)=\tfrac{\nabla \U(t,x)}{\U(t,x)}$, and their result implies the geometric ergodicity of $\nabla h$. Here we choose to work on a different quantity, the endpoint distribution of the directed polymer in random environment:
\[
u(t,x)=\frac{\U(t,x)}{\int_{\bT^d} \U(t,x')dx'}.
\] 
One can also view $u$ as the projective process, which appears frequently in the study of random dynamical systems. One of the main ingredients in our proof is the geometric ergodicity of the Markov process $\{u(t,\cdot)\}_{t\geq0}$, see Theorem~\ref{t.geomeEr}. The proof is inspired by the classical work of Sinai \cite{sinai}, which was actually on the Burgers equation.  Another key quantity is the intersection time of two independent polymer paths, similar to the so-called \emph{replica overlap}, which is used to distinguished different temperature regimes and  plays a crucial role in the study of thermodynamic limits of directed polymers, and in our case  corresponds to \[
\int_0^t \int_{\bT^{2d}} R(x-y)u(s,x)u(s,y)dxdyds.
\]
It is an additive functional of $\{u(t,\cdot)\}_{t\geq0}$, and the idea is to utilize the aforementioned geometric ergodicity to construct a corrector and extract a martingale from the above term. The geometric ergodicity is then used to drive a martingale central limit theorem which eventually leads to \eqref{e.conkpz}. A  closely related work is \cite{rosati}, where similar results were proved by a different approach, and a one force - one solution principle (modulo constants) for the KPZ equation was established.

For a more detailed discussion on directed polymers in random environment, we refer to the book \cite{comets2017directed} and the references therein.  For the endpoint distribution we are interested in, our approach is more related to that of  \cite{bakhtin1,bates2016endpoint,mukherjee}. For the polymer lying on the torus, we refer to the physics literature \cite{brunet,brunet1,brunet2} for some relevant discussions.

\subsection{Organization of the paper}

The   paper is organized as follows. In Section~\ref{s.2}, we sketch
the main ingredients used in the proof, including  the geometric
ergodicity of $\{u(t,\cdot)\}_{t\geq0}$ and  a martingale
decomposition for the free energy of the directed polymer. In
Section~\ref{s.3}, we derive a nonlocal SPDE satisfied by
$u$. Section~\ref{s.polymer} is devoted to proving the geometric
ergodicity of $\{u(t,\cdot)\}_{t\geq0}$ by following the approach of
\cite{sinai}. In Section~\ref{s.5}, we prove the main theorem through
constructing the corrector and applying the martingale central limit
theorem. In the appendix we recall some facts concerning the
Fortet-Mourier metric on the space of probability measures on polish
metric spaces, see Section \ref{s.FM} and recall some fairly standard facts
concerning  the stochastic heat equation,
see Section \ref{s.SHE}.

\subsection*{Acknowledgements}
We would like to thank two anonymous referees for many helpful suggestions which improved the presentation of the paper. Y.G. was partially supported by the NSF through DMS-1907928/DMS-2203014.
 T.K. acknowledges the support of NCN grant 2020/37/B/ST1/00426.

\section{Ideas of the proof}

\label{s.2}

\subsection{Endpoint distribution of directed polymer}

To study the SHE \eqref{e.she}, a closely related object is the directed polymer in random environment. We introduce some notations. For any $y\in \bT^d$, Let $\U(t,x;y)$ be  the solution to
\begin{equation}\label{e.defGreen}
\begin{aligned}
&\partial_t \U(t,x;y)=\tfrac12\Delta_x \U(t,x;y)+\U(t,x;y)\xi(t,x),\quad\quad t>0,\\
&\U(0,x;y)=\delta(x-y).
\end{aligned}
\end{equation}
In other words, $\U$ is the Green's function of the SHE.
For any probability measure $\nu$, denote the solution to  \eqref{e.she} with the initial data $\nu(dx)$ by $\U(t,x;\nu)$, then with the Green's function, we can write
\begin{equation}\label{e.UnuUy}
\U(t,x;\nu)=\int_{\bT^d} \U(t,x;y)\nu(dy).
\end{equation}
If $\nu(dx)=f(x)dx$ for some $0\leq f\in L^1(\bT^d)$, we   write
instead $\U(t,x;f)=\U(t,x;\nu)$. For a directed polymer of length $t$, with the starting point distributed according to $\nu$, its endpoint distribution is given by  
\begin{equation}\label{e.defPo}
u(t,x;\nu)=\frac{\U(t,x;\nu)}{\int_{\bT^d} \U(t,x';\nu)dx'}.
\end{equation}
Throughout the paper, we will omit the dependence of either $\U$, or
$u$ on the initial data $\nu$ when there is no danger of
confusion. Let  $\mathcal{M}_1(\bT^d)$ be the space of probability
measures on $\bT^d$, then for $t\geq 0$ and each realization of the
noise $\xi$, we can view $u(t)=u(t,\cdot)$ as the density of some
probability measure on $\bT^d$ and sometimes, with some abuse of notations, we write $u(t)\in \mathcal{M}_1(\bT^d)$.

\begin{remark}\label{r.fkpolymer}
In the case  when $R$ is  smooth, by the Feynman-Kac formula, \tk{see e.g.
\cite[Section 2.2]{bertini1995stochastic},} we have 
\begin{equation}
  \label{FK}
\U(t,x;y)=\E_x[e^{\int_0^t\xi(t-s,B_s)ds-\frac12R(0)t}\delta(B_t-y)],
\end{equation}
where $\E_x$ is the expectation with respect to the Wiener  measure $\bbP_x$  on $ C([0,+\infty);\bbR^d)$, with the paths starting at $x$,
i.e. $\bbP_x[B_0=x]=1$. 
Through a change of variable and the time reversal
property of a Brownian bridge, we have 
\[
\begin{aligned}
\U(t,x;y)=&\E_x[e^{\int_0^t\xi(s,B_{t-s})ds-\frac12R(0)t}\delta(B_t-y)]\\
=&\E_y[e^{\int_0^t\xi(s,B_{s})ds-\frac12R(0)t}\delta(B_t-x)],
\end{aligned}
\]
which implies
\tk{\begin{align}
      \label{031011-22}
&\U(t,x;\nu)=\int_{\bT^d} \nu(dy) \E_y[e^{\int_0^t
                 \xi(s,B_s)ds-\frac12R(0)t}\delta(B_t-x)]\notag\\
  &
    =\int_{\bT^d} p_t(x-y)\E_{y,x}^t[e^{\int_0^t
                 \xi(s,B_s)ds-\frac12R(0)t}]\nu(dy),
\end{align}
where $\E_{y,x}^t$ is the expectation with respect to the Brownian
bridge measure on $C([0,t];\bT^d)$ corresponding to $B_0=y$ and
$B_t=x$ and
\begin{equation}\label{e.defptx}
p_t(x)=\frac{1}{(2\pi t)^{d/2}}\sum_{n\in\bbZ^d}\exp\left\{-\frac{|x+n|^2}{2t}\right\}
\end{equation}
is the heat kernel on $\bT^d$. }
Therefore, one can write 
\[
u(t,x;\nu)=\frac{\int_{\bT^d} \nu(dy) \E_y[e^{\int_0^t \xi(s,B_s)ds}\delta(B_t-x)]}{\int_{\bT^d} \nu(dy) \E_y[e^{\int_0^t \xi(s,B_s)ds } ]},
\]
which is  the endpoint distribution of the directed polymer.
\end{remark}

Sometimes to
indicate the dependence of $\xi$ on the random event we shall write
\[
\xi=\big\{\xi(t,x;\om)\big\}_{(t,x;\om)\in \bbR\times\bT^d\times\Omega}.
\]  
In
addition, since the noise is space-time homogeneous, we
may and shall assume that there is an additive group
$\{\theta_{t,x}\}_{(t,x)\in  \bbR\times\bT^d}$ of ${\cal F}/{\cal F}$
measurable maps $ \theta_{t,x}:\Om\to\Om$ such that $\PP\circ
\theta_{t,x}^{-1}=\PP$ for all $(t,x)\in \bbR\times\bT^d$ and 
\begin{equation}
\label{012903-21}
\xi(t,x; \theta_{s,y}(\om))=\xi(t+s,x+y; \om), \quad (t,x),(s,y)\in \bbR\times\bT^d,\,\om\in\Omega.
\end{equation}
  As a consequence of the Markov property of $\{{\cal
    U}(t)\}_{t\geq0}$ and the white-in-time nature of the noise, we conclue the following.
\begin{lemma}\label{l.markov}
$\{u(t)\}_{t\geq0}$ is   a Markov process taking values in $\mathcal{M}_1(\bT^d)$.
\end{lemma}
The proof of the lemma can be found in Section \ref{MPd} of the Appendix.


One of the main technical parts of the paper is to prove the geometric ergodicity of the Markov process $\{u(t;\nu)\}_{t\geq 0}$. To state the result, we introduce some notations.

For any $p\ge1$ we let
$$
\mathcal{M}_1^{(p)}(\bT^d):=\Big\{\nu\in
\mathcal{M}_1(\bT^d):
\nu(dx)=f(x)
dx,\,f\in L^p(\bT^d)\Big\},
$$
equipped with the relative topology from $L^p(\bT^d)$. We
shall identify $\mathcal{M}_1^{(p)}(\bT^d)$ with the appropriate subset
of densities, i.e. 
\begin{equation}
\label{012403-21}
D^p(\bT^d)=D(\bT^d)\cap L^p(\bT^d)\quad \mbox{and }D^\infty(\bT^d)=D(\bT^d)\cap C(\bT^d).
\end{equation}
Here $ D(\bT^d)$ is the set of densities w.r.t. the
 Lebesgue measure on $\bT^d$. We will consider the Fortet-Mourier
 metric on $\mathcal{M}_1(\bT^d)$, see Section~\ref{s.FM} for a brief
 recollection of the facts concerning this metric. For any $F: \mathcal{M}_1(\bT^d)\to \bbR$, we define
\begin{equation}
\label{FLip}
\|F\|_{{\rm Lip}}:=\|F\|_\infty+ \sup\limits_{\mu\not=\nu,\mu,\nu\in \mathcal{M}_1(\bT^d)}\frac{|F(\mu)-F(\nu)|}{{\rm d}_{\mathrm{TV}}(\mu,\nu)}.
\end{equation}
Denote by ${\rm Lip}(\mathcal{M}_1(\bT^d))$ the space of all
functions $F$ for which $\|F\|_{{\rm Lip}}<+\infty$. Let $\mathcal{B}(\mathcal{M}_1(\bT^d))$ and $\mathcal{B}(D^p(\bT^d))$ be the Borel $\sigma-$algebra of $\mathcal{M}_1(\bT^d)$ and $D^p(\bT^d)$ respectively.

Define the
transition
probability densities
\begin{equation}
\label{PtA}
{\cal P}_t(\nu,A):=\EE 1_A(u(t;\nu)),\quad A\in {\cal B}({\cal
  M}_1 (\bT^d)),\,\nu\in {\cal
  M}_1 (\bT^d).
\end{equation}
\tk{For any $t>0$, $u(t,\cdot;\nu)$ is a
continuous function almost surely. This can be seen from formula
\eqref{e.defPo} and, in the case when $R$ is smooth, from
\eqref{031011-22}.
In the case
$R=\delta$ this follows from e.g. \cite[Proposition 5.1,
p. 515]{davar} for  the initial data is a positive and bounded
measurable function and \cite[Theorem 3.1]{chendalang}, when the initial data is a
measure (the latter paper considers the SHE on $\R$, but the proof
carries out verbatim  to the case of the equation on $\bT$).}
Thus, for any $t>0$, $\mathcal{P}_t(\nu,\cdot)$ is actually supported on $D^p(\bT^d)$ for any $p\in[1,+\infty]$, and we can consider transition probabilities ${\cal P}_t(\nu,\cdot)$
on ${\cal B}(D^p(\bT^d))$.

Define the transition probability operator by
\begin{equation}
\label{PtF}
{\cal P}_tF(\nu):=\int_{{\cal M}_1(\bT^d)}F(v){\cal P}_t(\nu,dv),\quad F\in
B_b({\cal M}_1 (\bT^d)),\,\nu\in  {\cal M}_1 (\bT^d).
\end{equation}
We shall also consider a semigroup $\{{\cal P}_t\}_{t\geq0}$ on
$B_b(D_p(\bT^d))$ for any $p\in[1,+\infty]$,  \tk{where $B_b(\cdot)$ denote the
space of all bounded, Borel measurable functions on a respective
metric space.}



\begin{theorem}
\label{t.geomeEr}
There exist a Borel probability measure $\pi_\infty$ on ${\cal M}_1(\bT^d)$ and
constants $C,\lambda>0$, which only depend on the
covariance function $R(\cdot)$ and dimension $d$  such that 
\begin{itemize}
\item[(i)] 
\begin{equation}
\label{032303-21}
\pi_\infty\big( D^p(\bT^d)\big)=1,
\quad p\in[1,+\infty],
\end{equation}
\item[(ii)] for any  $ F\in
{\rm Lip}({\cal M}_1(\bT^d)) $ we have
\begin{equation}
\label{011203-21}
\left\|{\cal
    P}_tF-\int_{{\cal M}_1(\bT^d)}F(u)\pi_\infty(du)\right\|_\infty \leq
Ce^{-\lambda t}\|F\|_{{\rm Lip}},\quad  t\geq0   .
\end{equation}
The above result also holds for $D^p(\bT^d)$ in place of ${\cal M}_1(\bT^d)$ for any $p\in[1,+\infty]$.
\end{itemize}
\end{theorem}
The proof of the theorem is given in Section~\ref{s.polymer}.

\subsection{Martingale decomposition}

We will first prove Theorem~\ref{t.mainth} for the partition function of the
directed polymer, defined as  
\begin{equation}\label{e.defZ}
Z_t=\int_{\bT^d} \U(t,x)dx.
\end{equation}
By the mild formulation of SHE, one can write 
\begin{equation}\label{e.mild}
\U(t,x)=\int_{\bT^d} p_t(x-y)\nu(dy)+\int_0^t\int_{\bT^d} p_{t-s}(x-y)\U(s,y)\xi(s,y)dyds,
\end{equation}
where $p_t(x)$ - the heat kernel on $\bT^d$ - is given by \eqref{e.defptx}. Thus, recalling that $\nu$ is a probability measure, we have
\begin{equation}\label{e.4211}
Z_t= 1+\int_0^t\int_{\bT^d}\U(s,y)\xi(s,y)dyds,
\end{equation}
which is a positive  local martingale, whose quadratic variation
equals  
\begin{equation}\label{e.qvZ}
\la Z\ra_t=\int_0^t \cR\big(\U(s)\big) ds.
\end{equation}
Here for $R\in L^\infty(\bT^d)$ we let
\begin{equation}
\label{fR}
\begin{split}
&\cR(u,v):=\int_{\bT^{2d}}R(x-y)u(x)v(y)dxdy,\\
&
\cR(u):=\cR(u,u)=\int_{\bT^{2d}}R(x-y)u(x)u(y)dxdy,\quad  u,v\in L^1(\bT^d).
\end{split}
\end{equation}
If $R(\cdot)=\delta(\cdot)$, the above definition holds for  $u,v\in L^2(\bT^d)$.

By Lemma~\ref{l.bdQVZ}, we know that $\EE\la Z\ra_t^p \leq C(t,p)$, for some
constant $C(t,p)<+\infty$ depending only on $t$ and $p\in[1,+\infty)$.
Hence  $\big(Z_t\big)_{t\ge0}$ is   a martingale. 
Define a local martingale $(M_t)_{t\geq0}$ through the SDE
\begin{equation}
dM_t=Z_t^{-1}dZ_t, \quad\quad M_0=0.
\end{equation}
By the It\^o formula
\begin{equation}\label{e.relationZM}
Z_t=\exp\big\{M_t-\tfrac12\la M\ra_t\big\}.
\end{equation}
We can write  $M_t$ and its quadratic variation in terms of the endpoint distribution of the directed polymer:
\begin{equation}\label{e.ExpM}
\begin{aligned}
M_t=\int_0^t\int_{\bT^d} Z_s^{-1}\U(s,y)\xi(s,y)dyds=\int_0^t\int_{\bT^d} u(s,y)\xi(s,y)dyds,
\end{aligned}
\end{equation}
and 
\begin{equation}\label{e.qvM}
\la M\ra_t 
=\int_0^t \cR(u(s))ds.
\end{equation}
The following lemma shows that $M$ is a square-integrable martingale.
\begin{lemma}
For any $t>0$ and $p\in[1,+\infty)$ we have $\EE[\la M\ra_t^p]\leq C(t,p,d,R)$.
\end{lemma}
\begin{proof}
For the case of $R$ being bounded, we have $\la M\ra_t\leq
t\|R\|_{L^\infty(\bT^d)}$. For the case of $R(\cdot)=\delta(\cdot)$,
we have $\la M\ra_t= \int_0^t  \|u(s)\|_{L^2(\bT^d)}^2 ds$.  Applying
the triangle inequality we further derive 
\[
\EE[\la M\ra_t^p]^{1/p} \leq C\int_0^t \EE[\|u(s)\|_{L^2(\bT^d)}^{2p}]^{1/p}ds \leq \gu{C\int_0^t \EE[\|u(s)\|_{L^{2p}(\bT^d)}^{2p}]^{1/p}ds},
\]then the lemma is a direct consequence of Lemma~\ref{l.mmbdu}.
\end{proof}

Then we can write 
\begin{equation}\label{e.delogZt}
\begin{aligned}
\log Z_t=&M_t-\tfrac12\la M\ra_t\\
=&\int_0^t\int_{\bT^d} u(s,y)\xi(s,y)dyds-\tfrac12\int_0^t \cR(u(s))ds.
\end{aligned}
\end{equation}
The term $\int_0^t \cR(u(s))ds$ is an additive functional of the
Markov process $(u(s))_{s\geq0}$, which will be the key object to
study in the paper. 
Define 
\begin{equation}\label{e.defgamma}
\gamma:=\tfrac12\int_{\mathcal{M}_1(\bT^d)} \cR(u)\pi_\infty(du),
\end{equation}
where $\pi_\infty$ is the  unique invariant measure for the process (Theorem  \ref{t.geomeEr}),
and \begin{equation}\label{e.cenR}
\tilde{\cR}=\cR-2\gamma.
\end{equation}
It will become clear later (see Lemma~\ref{l.nondege1} below) that $\gamma\in(0,+\infty)$.

Now we write
\begin{equation}\label{e.delogZt1}
\log Z_t+\gamma t=\int_0^t\int_{\bT^d} u(s,y)\xi(s,y)dyds-\tfrac12\int_0^t \tilde{\cR}(u(s))ds.
\end{equation}
The idea is to utilize Theorem~\ref{t.geomeEr} to solve the Poisson equation associated with $\tilde{\cR}$ and extract the martingale part from $\int_0^t \tilde{\cR}(u(s))ds$, which combines with the martingale $M_t$ to converge to the limiting Gaussian random variable, after a rescaling:
\begin{theorem}\label{t.partitionfunction}
As $t\to\infty$, we have 
\[
\frac{\log Z_t+\gamma t}{\sqrt{t}}\Rightarrow N(0,\sigma^2).
\]
In addition, we have
\begin{equation}
\label{012204-21}
\si^2\ge\int_{\bT^d}R(x)dx.
\end{equation}
\end{theorem}

\begin{remark}[Process-level convergence]
\label{r.process}
Combining Theorems~\ref{t.geomeEr} and \ref{t.partitionfunction}, we actually have the following picture: as $t\to\infty$,  \begin{equation}\label{e.4271}
\log u(t,x)=\log \tfrac{\U(t,x)}{Z_t}=h(t,x)-\log Z_t\Rightarrow \log \rho(x)
\end{equation}
in distribution in $C(\bT^d)$, where the limit $\rho$ is sampled from the invariant measure $\pi_\infty$. Therefore, in large time, the random height function $h$ is approximately $\log \rho$ shifted by the random constant $\log Z_t$, and this random constant is roughly  Gaussian distributed after a centering and rescaling. 
 We also have 
\[
h(t,x)-h(t,0)\Rightarrow \log \rho(x)-\log \rho(0)
\]
in distribution in $C(\bT^d)$. 

According to the results of  \cite{funaki}, see Theorems 1.1 and
  2.1,    in the case of $1+1$ spacetime white noise, for each $t\geq0$ and in stationarity
  the law of
  $\left(h(t,x)-h(t,0)\right)_{x\in\bT}$ over $C(\bT)$ coincides with  that of
  $(B(x))_{x\in\bT}$ - a
 Brownian bridge  satisfying $B(0)=B(1)=0$.
The invariant measure $\pi_\infty$ in this case is given by the law of
$$
\rho(x)=\frac{e^{B(x)}}{\int_{\bT}e^{B(x')}dx'}, \quad x\in\bT
$$
over $D^\infty(\bT)$,
  see also \cite{bertini,gubinelli2017kpz}. 
\end{remark}

\section{SPDE for endpoint distribution of directed polymer}

\label{s.3}

Define the Fourier coefficients of the spatial covariance of the noise
\begin{equation}\label{e.defrk}
 \hat r_k= \int_{\bT^d}R(x)e^{-i 2\pi k\cdot x}dx,\quad k\in\bbZ^d.
\end{equation}
They are non-negative. \gu{The following is a standard Fourier representation of the cylindrical Wiener process, see \cite[Section 4.1]{daza}.}

We can write formally that 
\begin{equation}
\label{xtx}
\xi(t,x)=\sum_{k\in \Z^d} e_k(x)\dot{w}_k(t), \quad\quad e_k(x)=\sqrt{\hat{r}_k}e^{i2\pi k\cdot x},
\end{equation}
where $(w_k)_{k\in\Z^d}$ are   one dimensional complex-valued
zero mean Wiener processes satisfying
\begin{equation}
\label{022204-21}
\gu{w_k^*}(t)=w_{-k}(t),\quad \EE[w_k(t)
\gu{w_\ell^*}(s)]=\delta_{k,\ell}\min(t,s),\quad k,\ell\in\Z^d,\,t,s\ge0.
\end{equation}
Here $\delta_{k,\ell}$ is the Kronecker symbol.  We also define 
\begin{equation}\label{e.defW}
W(t,x)=\sum_{k\in\Z^d} e_k(x)w_k(t)
\end{equation}
and write 
\[
W(t)=W(t,\cdot),\quad\quad dW(t,x)=\xi(t,x)dt.
\]
We formulate the following assumption:
\begin{assumption}\label{a.smooth}
 There exists $N>0$ such that $\hat{r}_k=0$ for $  |k|>N$.
\end{assumption}
Under the above assumption, the summations in \eqref{xtx} and
\eqref{e.defW} are finite. 

The main result of this section is to derive an SPDE satisfied by the
endpoint distribution of the directed polymer, \gu{under the Assumption~\ref{a.smooth}. When the noise is smooth in the spatial variable, one could directly apply It\^o formula to obtain the equation satisfied by $u$. It may be possible to study the singular case, e.g., when $d=1$ and the noise is white in space, and to make sense of the SPDE derived in \eqref{e.spde} below, in a mild formulation. For our purpose, it is enough to consider the case of a smooth noise here, and the derived SPDE will be used later in proving the semi-martingale decomposition in Proposition~\ref{p.madefinal} below. Then through an approximation argument, the Assumption~\ref{a.smooth} will be removed to cover the white noise case}. Recall that
$u(t,x;\nu)$ has been defined in \eqref{e.defPo} for any $\nu\in\mathcal{M}_1(\bT^d)$.
In this section  the initial data is assumed to be of the form
$\nu(dx)=v(x)dx$, where $v\in D(\bT^d)\cap C^\infty(\bT^d)$. To simplify the notation, we
omit the initial data and write   $u(t)=u(t,\cdot)$ and $\U(t)=\U(t,\cdot)$.  

\gu{Throughout the paper, we will use $\star$ to represent the spatial convolution, i.e., 
\[
f\star g(x)=\int_{\bT^d} f(x-y)g(y)dy.
\]}
\tk{Denote by $H^k(\bT^d)$ the usual Sobolev space, the completion of $C^\infty(\bT^d)$ in the
  norm
  $$\|v\|_{H^k(\bT^d)}:=\left\{\sum_{0\le m_1+\ldots+m_d\le k}\|\partial_{x_1}^{m_1}\ldots \partial_{x_d}^{m_d}
    v\|_{L^2(\bT^d)}^2\right\}^{1/2}.
$$}
\begin{proposition}\label{p.spde}
Suppose that
Assumption~\ref{a.smooth} is in force, and the initial data $v\in
D(\bT^d)\cap C^\infty(\bT^d)$. \gu{Then} $u$ belongs to
$C([0,\infty),H^k(\bT^d))$, for any $k\geq1$. In addition  it is a strong solution to 
\begin{equation}\label{e.spde}
\begin{aligned}
du(t)=&\Big\{\tfrac12\Delta u(t) +u(t) {\cal R}\big(u(t)\big)-u(t)R\star u(t)\Big\}dt\\
&+u(t)dW(t)-u(t)\int_{\bT^d}u(t,x)dW(t,x),\quad u(0)=v.
\end{aligned}
\end{equation}
\end{proposition}

\begin{proof}
Under Assumption~\ref{a.smooth} and by the fact that $v\in C^\infty(\bT^d)$, we have
\begin{equation}\label{e.4131a}
\U(t)=v+\int_0^t \tfrac12\Delta \U(s)ds+\int_0^t \U(s)dW(s).
\end{equation}
\tk{Thanks to \cite[Theorem 6.7, p. 157]{daza} (that can be applied
  thanks to Assumption~\ref{a.smooth})  we conclude that there
  exists the unique mild solution to the SHE \eqref{e.4131a} and
  such that ${\cal U}(t)\in H^k(\bT^d)$, $t\ge0$ for any
  $k\ge0$. Furthermore, by \cite[Proposition 6.29, p. 176 and
  Theorem 6.30, p. 177-178]{daza}, the
  solution is strong.}

By \eqref{e.4211}, we have (which can also be obtained by integrating the above equation in $x$)
\begin{equation}\label{e.4131b}
Z_t= \int_{\bT^d}\U(t,x)dx=1+\int_0^t \int_{\bT^d}\U(s,x) dW(s,x).
\end{equation}
Using the above two equations, we have 
\begin{equation}\label{e.4131}
\begin{aligned}
&\la Z\ra_t=\int_0^t \cR\big(\U(s)\big)ds,\\
&\la \U,Z\ra_t=\int_0^t\U(s) \gu{R}\star \U(s)ds.
\end{aligned}
\end{equation}
By \eqref{e.negammZ}, we know that  $(Z_t)_{t\geq0}$ is continuous and positive almost surely. Then by It\^o's formula, we have 
\[
\begin{aligned}
du(t)=&d[\U(t)Z_t^{-1}]\\
=&Z_t^{-1}d\U(t)-\U(t)Z_t^{-2}dZ_t+\U(t)Z_t^{-3}d\la Z\ra_t-Z_t^{-2}d\la \U,Z\ra_t.
\end{aligned}
\]
Invoking \eqref{e.4131a} - \eqref{e.4131}, the proof is complete.
\end{proof}

We introduce the following operators:
\begin{equation}\label{e.defAB}
\begin{aligned}
&\mathscr{A}(f):=\tfrac12\Delta f+f \la f,R\star f\ra_{L^2(\bT^d)}-fR\star f, \quad\quad f\in H^2(\bT^d),\\
&\mathscr{B}_k(f): =fe_k-f\la f,e_{-k}\ra_{L^2(\bT^d)},\quad\quad f\in L^2(\bT^d).
\end{aligned}
\end{equation}
Note that $e_k^*=e_{-k}$, so $\la f,e_{-k}\ra_{L^2(\bT^d)}=\int_{\bT^d} f(x)e_k(x)dx$, and for real-valued $f$, we have 
\[
\mathscr{B}_k(f)^*=\mathscr{B}_{-k}(f).
\] With the above notations, the equation \eqref{e.spde} can be written in a more compact form:
\[
du(t)=\mathscr{A}(u(t))dt+\sum_{k\in\Z^d } \mathscr{B}_k(u(t))dw_k(t),
\]
and we have the following lemma:
\begin{lemma}\label{l.bdspde}
Under Assumption~\ref{a.smooth} and further assume $v\in D(\bT^d)\cap C^\infty(\bT^d)$, then for any $T>0$ and $p\in[1,+\infty)$, there exists $C=C\tk{(T,p,R,d,v)}$ such that 
\begin{equation}
  \label{X010711-22}
\begin{aligned}
&\sup_{t\in[0,T]}\EE\big[\|\mathscr{A}(u(t))\|_{L^2(\bT^d)}^p\big] \leq C,\\
&\sum_{k\in\Z^d } \sup_{t\in [0,T]}\EE\big[\|\mathscr{B}_k(u(t))\|_{L^2(\bT^d)}^p\big] \leq C.
\end{aligned}
\end{equation}
\end{lemma}
\begin{proof}
  \tk{Using \eqref{FK} we can write
\begin{equation}
  \label{FK1}
\U(t,x;v)=\E_0[e^{\int_0^tdW(s,x+B_{t-s})-\frac12R(0)t}v(x+B_t)],
\end{equation}
where we recall that $\E_0$ is the expectation on the Brownian motion starting at $0$.
Thanks to Assumption~\ref{a.smooth}, we conclude from \eqref{FK1} 
  that for each $p\in[1,+\infty)$ and $T>0$ there exists a constant
  $C(T,p,d,v)$ such that
   \begin{equation}
   \label{X010711-22a}
 \sup_{t\in [0,T]}\EE\big[\|{\cal U}(t;v)\|_{H^2(\bT^d)}^p\big]
 \leq C.
 \end{equation}
  Since $u(t)=\U(t)Z_t^{-1}$, the result   follows from
    \eqref{X010711-22a} and 
\eqref{e.negammZ}.}
  \end{proof}

\section{Geometric ergodicity of directed polymer on torus}
\label{s.polymer}


In this section, we study the endpoint distribution of the directed polymer, defined in \eqref{e.defPo}, and prove Theorem~\ref{t.geomeEr}. \gu{The approach we adopt in this section is motivated by the classical work of Sinai \cite{sinai} in which he proved geometric ergodicity of the stochastic Burgers equation on a torus. The difference is that, we consider the polymer endpoint process $u(t,\cdot)=\U(t,\cdot)/\int \U(t,\cdot)$ instead of $\nabla \U(t,\cdot)/\U(t,\cdot)$. Suppose there is no random environment, then $u(t,x)$ is the standard heat kernel on the torus, thus, on a heuristic level, one can view our approach as a complicated way of proving the geometric ergodicity of a standard Brownian motion on a torus: we discretize the time variable and study the related random walk on a torus, and since the transition density of the random walk is bounded below by a positive constant, which is the standard Doeblin condition, we conclude the geometric ergodicity of the random walk process. With the presence of the random environment, the polymer endpoint is performing a similar random walk, with the transition kernel depending on the random environment. We will show that with a positive probability, the transition kernel is bounded below by a positive constant (the precise version is \eqref{021903-21b} below) and the same coupling argument as in the case of the Doeblin condition enables us to prove the ``renewal'' property which eventually leads to the geometric ergodicity.}

\subsection{Factorization }

Recall that we consider the solution to SHE with the initial data $\nu\in \mathcal{M}_1(\bT^d)$:
\begin{equation}\label{e.SHEf}
\partial_t \U=\tfrac12\Delta\U+\U\xi,\quad\quad \U(0,dx)=\nu(dx).
\end{equation}

For any $t>s$ and $x,y\in\bT^d$, define $\mathcal{Z}_{t,s}^\omega(x,y)$ as the solution to 
\begin{equation}\label{e.defZts}
\begin{aligned}
&\partial_t \mathcal{Z}_{t,s}^\omega(x,y)=\tfrac12\Delta_x \mathcal{Z}_{t,s}^\omega(x,y)+\mathcal{Z}_{t,s}^\omega(x,y)\xi(t,x), \quad\quad t>s,\\
&\mathcal{Z}_{s,s}^\omega(x,y)=\delta(x-y).
\end{aligned}
\end{equation}
In other words, $\mathcal{Z}_{t,s}^\omega(x,y)$ is the propagator of the SHE from $(s,y)$ to $(t,x)$. We will keep its dependence on $\omega$ throughout this section. In light of \eqref{012903-21}, we have
\begin{equation}\label{022903-21}
\mathcal{Z}_{t,s}^{\theta_{r,z}(\om)}
(x,y)=\mathcal{Z}_{t+r,s+r}^{\om} (x+z,y+z).
\end{equation}
\tk{For any interval $I\subset \R$, denote by $\F_{I}$  the $\si$-algebra
generated by the random variables
\begin{equation}
\label{cF}
X:=\Phi\Big(\int_I\int_{\bT^d}\varphi_1(s,x)\xi(s,x)dsdx,\ldots \int_I\int_{\bT^d}\varphi_{N}(s,x)\xi(s,x)dsdx\Big),
\end{equation} 
where $N\ge1$,
$\Phi:\bbR^{N}\to\bbR$ is a  bounded and Borel measurable  function, and 
$\varphi_\ell$, $1\le \ell\le N$ are  functions   that are compactly
supported and  infinitely many
times differentiable  in $I\times \bT^d$.} It is clear that if two intervals $I_1$ and $I_2$ do not intersect, i.e., $I_1\cap I_2=\emptyset$, then $\F_{I_1}$ is independent of $\F_{I_2}$.

We shall also write 
\begin{equation}
\label{cFt}
{\cal F}_t:={\cal F}_{[0,t]},\quad t\ge0.
\end{equation}
Note that the random variable   $\cZ_{t,s}^\om$
is  ${\cal F}_{[s,t]}$-measurable.

Fix $t>1$. 
We have 
\begin{equation}\label{e.3111}
\U(t,x;\nu)=\int_{\bT^d} \mathcal{Z}_{t,t-1}^\om(x,y_1) \U(t-1,y_1;\nu) dy_1.
\end{equation}

Let 
\begin{equation}\label{e.defN}
N(t):=[t]+1.
\end{equation} 
Sometimes, when it leads to no confusion, for abbreviation we shall write $N=N(t)$. By iterating \eqref{e.3111}, we further derive that
\begin{equation}\label{e.3112}
\begin{aligned}
\U(t,x;\nu)=\int_{\bT^{Nd}}& \mathcal{Z}_{t,t-1}^\om (x,y_1) \left(\prod_{j=1}^{N-2} \mathcal{Z}^\om_{t-j,t-j-1}(y_j ,y_{j+1} ) \right)\\
&\times \mathcal{Z}^\om_{t-N+1,0}(y_{N-1} ,y_{N} ) 
dy_1\ldots dy_{N-1}\nu(dy_N ).
\end{aligned}
\end{equation}
With the  convention of $
y_0=x$, 
we can write 
\begin{equation}\label{e.31121}
\begin{aligned}
\U(t,x;\nu)=  \int_{\bT^{Nd}} \prod_{j=0}^{N-1} \mathcal{Z}^\om_{t-j, {(t-j-1)_+}}(y_j ,y_{j+1} )  dy_1\ldots dy_{N-1}\nu(dy_N ).
\end{aligned}
\end{equation}
 {Here $(x)_+:=\max\{x,0\}$}.

\subsection{Construction of a Markov chain }

\label{sec5.1}

In this section, we construct a  time reversed Markov chain $\{Y_n\}_{n=1}^N$, which
depends on the random environment, and allows
to express the endpoint density  $ u(t,x;\nu)$, see \eqref{e.defPo},
as an average with respect to $Y_1$, see \eqref{e.densityY1} below. It turns out  that, with some
positive probability, the random realizations of the transition probability
density kernel of the chain are uniformly positive,
see \eqref{021903-21} and \eqref{021903-21b}. This, in turn, allows us to 
use a coupling argument that leads to the proof of the geometric
ergodicity of the  density  process $\big( u(t;\nu)\big)_{t\ge0}$.

The  {random} Markov chain takes values in $\bT^d$ and is
constructed as follows. For a fixed $\om\in\Om$, we run the chain backward in time. Let
$\pi_N^\om(y_N)$ be the  {(random)}  density of $Y_N$ and
$\pi_k^\om(y_k\,|\,y_{k+1})$ be the  {(random)}  transition density
kernels:
\begin{equation}\label{e.trandensity}
\begin{aligned}
&\pi_1^\om(y_1\,|\,y_2)=\frac{\cZ_{t-1,t-2}^\om (y_1,y_2)}{\int_{\bT^d} \cZ_{t-1,t-2}^\om (z_1,y_2)dz_1}\\
&\pi_k^\om(y_k\,|\,y_{k+1})=\cZ ^\om_{t-k,(t-k-1)_+}(y_k,y_{k+1})\\
&
\times \frac{\int_{\bT^{(k-1)d}} \cZ_{t-1,t-2}^\om (z_1,z_2)\ldots \cZ ^\om_{t-k+1,t-k}(z_{k-1},y_{k}) dz_{1,k-1}}{\int_{\bT^{kd}} \cZ ^\om_{t-1,t-2}(z_1,z_2)\ldots \cZ ^\om_{t-k,(t-k-1)_+}(z_k,y_{k+1})dz_{1,k}}
\end{aligned}
\end{equation}
for $k=1,\ldots, N-1$. Here $dz_{1,k}:=dz_1\ldots dz_k$.
Finally we let
\begin{align}
\label{piN}
&
\pi_N^\om(dy_N):=A^\om_{N, \nu}  \int_{\bT^{(N-1)d}}\left(\prod_{j=1}^{N-2}
  \mathcal{Z}^\om_{t-j,t-j-1}(z_j ,z_{j+1} ) \right)\\
&
\times \mathcal{Z}^\om_{t-N+1, 0}(z_{N-1} ,y_N )\nu(dy_N )dz_{1,N-1},\notag
\end{align}
with the (random) constant $A^\om_{N, \nu}$ given by 
\[
A_{N,\nu}^\om:=\left(\int_{\bT^{Nd}}\prod_{j=1}^{N-1}
  \mathcal{Z}^\om_{t-j,  {(t-j-1)_+}}(z_j ,z_{j+1} )  \nu(dz_N
  )dz_{1,N-1}\right)^{-1}.
\]

{Sometimes, we shall write $\pi_N^\omega(dy_N;\nu)$  when we wish to highlight the dependence of 
the distribution on the initial data $\nu$. Also, when $\nu(dx)=f(x)dx$ for
$f\in D(\bT^d)$ we shall write $\pi_N^\omega(y_N;f)$ as the density of
$\pi_N^\omega(dy_N;\nu)$ w.r.t. the Lebesgue measure}.
Note that $\cZ_{t,s}^\om (x,y)\in \F_{[s,t]}$ for all $x,y\in\bT^d$,
cf \eqref{cF}.
Let $\bbP_{\pi_N}^\om$ be the path measure on $\bT^{Nd}$ 
  constructed from the above transition probabilities. It
  is  
  the  law of a time reversed Markov chain $\{Y_n\}_{n=1}^N$, where
  $$
Y_k(y_1,\ldots,y_N):=y_k, \quad\quad  (y_1,\ldots,y_N)\in
\bT^{Nd},\,k\in\{1,\ldots,N\}.
$$
Let $\E_{\pi_N}^\om$ denote the   expectation with respect to
  this measure. Suppose that $\nu(dx) =f(x) dx$. Then,
the joint density of $(Y_1,\ldots,Y_N)$ under  $\bbP_{\pi_N}^\om$ equals 
\begin{align}
\label{012003-21a}
\varphi^\om_{1,N}(y_1,\ldots,y_N)&=
\pi_1 ^\om (y_1\,|\,y_2)\ldots \pi_{N-1}^\om (y_{N-1}\,|\,y_N)\pi_N ^\om (y_N;f)\\
&
=A_{N, \nu}^\om\left(\prod_{j=1}^{N-1} \mathcal{Z}^\om_{t-j, (t-j-1)_+}(y_j ,y_{j+1} ) \right)f(y_N ).\notag
\end{align} 
Here $A_{N, \nu}^\om$ is the (random) normalizing constant.

For a given $2\le k\le N-1$ denote by $\bbP_{k}^\om$ the path
measure on $\bT^{Nd}$, corresponding to a time reversed Markov chain
obtained from the transition probability densities
\begin{equation}
\label{042003-21}
p_j^\om(y_j\,|\,y_{j+1})=
\left\{
\begin{array}{ll}
\pi_j^\om(y_j\,|\,y_{j+1}),& 1\le j\le k-1,\\
&\\
1,& k\le j\le N-1.
\end{array}
\right.
\end{equation}
and $p_N^\om(y_N)\equiv 1$. The respective expectation shall be
denoted by $\E_{k}^\om$.
The joint density of $(Y_1,\ldots,Y_k)$ under  $\bbP_{k}^\om$ equals 
\begin{align}
\label{012003-211}
\varphi^\om_{1,k}(y_1,\ldots,y_k)&=
p_1 ^\om (y_1\,|\,y_2)\ldots p_{k-1}^\om (y_{k-1}\,|\,y_{k})\\
&
=A_{k,\1}^\om\left(\prod_{j=1}^{k-1} \mathcal{Z}^\om_{t-j, t-j-1}(y_j ,y_{j+1} ) \right), \notag
\end{align} 
with the (random) constant 
\begin{equation}\label{012003-212}
A_{k,\1}^\om:=\left(\int_{\bT^{kd}}\prod_{j=1}^{k-1}
  \mathcal{Z}^\om_{t-j,t-j-1 }(z_j ,z_{j+1} ) dz_{1,k}\right)^{-1}.
\end{equation}

In other words, comparing the Markov chains sampled from the path measures $\Pb_k^\omega$ and $\Pb_{\pi_N}^\omega$, the only difference lies in the distribution of $(Y_k,\ldots,Y_N)$. For $\Pb_k^\omega$, we sample $Y_k,\ldots,Y_N$ independently from the uniform distribution on $\bT^d$, while we follow \eqref{e.trandensity} and \eqref{piN} for $\Pb_{\pi_N}^\omega$.


By  \eqref{e.3112} and \eqref{012003-21a}, one can write
\begin{equation}\label{e.UY1}
\begin{aligned}
\U(t,x;\nu)=&
\big(A^\om_{N,\nu}\big)^{-1}\int_{\bT^{Nd}}\cZ_{t,t-1}^\om(x,y_1)\pi_1^\om(y_1\,|\,y_2)\ldots \pi_{N-1}^\om(y_{N-1}\,|\,y_N)\pi_N^\om(dy_N)dy_{1,N-1}\\
=&\big(A^\om_{N,\nu}\big)^{-1}\E_{\pi_N}^\om[\cZ^\om_{t,t-1}(x,Y_1)].
\end{aligned}
\end{equation}
This, in turn implies 
\begin{equation}\label{e.densityY1}
  u(t,x;\nu)=\frac{\U(t,x;\nu)}{\int_{\bT^d}\U(t,x';\nu)dx'}=\frac{\E_{\pi_N}^\om[\cZ_{t,t-1}^\om(x,Y_1)]}{\int_{\bT^d}\E_{\pi_N}^\om[\cZ_{t,t-1}^\om (x',Y_1)]dx'}.
\end{equation}
We note that in the above formula, the dependence on $\nu$ is through the distribution $\pi_N^\omega$. It is also clear that the dependence of $u(t,x;\nu)$ on the random environment $\{\xi(s,\cdot):s\leq t-1\}$ is only through the random variable $Y_1$, hence it suffices to study the mixing property of the Markov chain.

\subsection{Properties of the Markov chain}
\begin{lemma}\label{l.ulbd}
For any $p\in[1,+\infty)$, there exists $C>0$ only depending on $p,d,R$ such that 
\begin{equation}\label{e.ulbd}
\EE\Big[\big(\inf_{x,y\in\bT^d}
\cZ_{t,t-1}^\omega(x,y)\big)^{-p}\Big]+\EE\Big[\big(\sup_{x,y\in\bT^d}\cZ_{t,t-1}^\omega(x,y)\big)^p\Big]\leq
C,\quad \quad  t>1.
\end{equation}
\end{lemma}

It is a classical result for the stochastic heat equation, and for the convenience of readers, we present its proof in Section~\ref{s.sheOther}. \gu{
The above lemma is key to our argument and   the essential reason why the same result does NOT hold in the whole space. It basically states that with high probability $\cZ_{t,t-1}^\omega(x,y)$ is uniformly bounded away from zero and infinity. Heuristically, one can write  
\[
\cZ_{t,t-1}^\omega(x,y)=p_1(x-y) \frac{\cZ_{t,t-1}^\omega(x,y)}{p_1(x-y)},
\] and convince oneself that $\cZ_{t,t-1}^\omega$ is roughly of the same order as $p_1$. Recall that $p_t(x)$ is the standard heat kernel, so the key is, on the torus, we have $\inf_{x\in \bT^d}p_1(x)>0$, which  is not the case for the whole space.
}

\gu{\begin{remark}\label{r.openkpz}
If we consider the KPZ equation with a Neumann boundary condition, then the geometric ergodicity of the process $\{u(t,\cdot)\}_{t\geq0}$, as stated in Theorem~\ref{t.geomeEr}, follows verbatim the same argument in this section, provided that one can prove the above Lemma~\ref{l.ulbd} for the propagator of SHE with the corresponding Robin boundary condition, and this would lead to the geometric ergodicity of the process $\{h(t,\cdot)-h(t,0)\}_{t\geq 0}$. Indeed, this was done in the recent preprint \cite{shalin1} (see also \cite{alisa} on  the uniqueness of the invariant measure of the open KPZ equation in certain cases, along the line of proving the strong Feller property) .
\end{remark}}

For any $t>1$ and $\delta>0$, define the event 
\begin{equation}\label{e.defB}
{B_t(\delta)}:=\left\{\om\in\Omega:\,\frac{\inf_{x,y,z\in\bT^d}\cZ_{t+1,t}^\om(x,y)\cZ_{t,t-1}^\om (y,z)}{\sup_{x,z\in \bT^d}\int_{\bT^d} \cZ_{t+1,t}^\om (x,y')\cZ_{t,t-1}^\om (y',z)dy'}>\delta\right\},
\end{equation}
which belongs to   $\F_{[t-1,t+1]}$, cf \eqref{cF}, and, {thanks to
\eqref{012903-21},
\begin{equation}\label{e.defBa}
\theta^{-1}_{s,y}\Big(B_t(\delta)\Big)=B_{t+s}(\delta),\quad (s,y)\in\bbR_+\times\bT^d.
\end{equation}}
 Lemma~\ref{l.ulbd} implies 
\begin{lemma}\label{l.lowerbd}
There exists a constant $C>0$ only depending on $d,R$ such that  
\begin{equation}
\mathbf{P}[B_t(\delta)]\geq 1-C\delta, \quad \mbox{for all $t>1$ and $\delta>0$}.
\end{equation}
\end{lemma}

\begin{proof}
First it is clear that by the time homogeneity of the random
environment, the probability $\mathbf{P}[B_{t }(\delta)]$ does not depend on $t$. To simplify the notation, denote 
\[
W_1=\inf_{x,y,z\in\bT^d}\cZ_{t+1,t}^\omega(x,y)\cZ_{t,t-1}^\omega(y,z), \quad\quad W_2=\sup_{x,z\in \bT^d}\int_{\bT^d} \cZ_{t+1,t}^\omega(x,y)\cZ_{t,t-1}^\omega(y,z)dy,
\] 
then for any $\delta>0$, we have 
\[
\begin{aligned}
\mathbf{P}[B_{t }(\delta)]=&\mathbf{P}[W_1W_2^{-1}>\delta]=1-\mathbf{P}[W_1^{-1}W_2\geq \delta^{-1}]\\
\geq &1-\delta\EE[W_1^{-1}W_2]\geq 1-\delta \sqrt{\EE[W_1^{-2}]\EE[W_2^2]}.
\end{aligned}
\]
\gu{For the factor containing  $W_1$, we have
\[
W_1^{-2}\leq  (\inf_{x,y\in\bT^d}\cZ_{t+1,t}^\omega(x,y))^{-2}(\inf_{y,z\in\bT^d}\cZ_{t,t-1}^\omega(y,z))^{-2},
\]
and by independence and Lemma~\ref{l.ulbd}, this further implies
\[
\EE [W_1^{-2}]\leq C.
\]
A similar argument applies to the term containing $W_2$.} Thus, there exists a constant $C>0$ depending on $d,R$ such that 
\[
\sqrt{\EE[W_1^{-2}]\EE[W_2^2]}\leq C.
\]
The proof is complete.
\end{proof}

The following lemma holds.
\begin{lemma}\label{l.lbdtr}
We have
\begin{equation}
\label{011903-21}
\inf_{y_k,y_{k+1}\in\bT^d}\pi_k^\om(y_k\,|\, y_{k+1})> \delta,\quad\quad
\om\in B_{ t-k}(\delta) 
\end{equation}
for all $k=2,\ldots,N-2$.
\end{lemma}
\begin{proof}
Recall that 
\[
\begin{aligned}
\pi_k^\om(y_k\,|\,y_{k+1})=&\frac{\int_{\bT^{(k-1)d}} \cZ_{t-1,t-2}^\omega(y_1,y_2)\ldots \cZ_{t-k,t-k-1}^\om(y_k,y_{k+1})dy_1\ldots dy_{k-1}}{\int_{\bT^{kd}} \cZ_{t-1,t-2}^\om(y_1,y_2)\ldots \cZ_{t-k,t-k-1}^\om(y_k,y_{k+1})dy_1\ldots dy_{k-1}dy_k}\\
= &\frac{ \int_{\bT^d} g(y_{k-1})\cZ_{t-k+1,t-k}^\om(y_{k-1},y_k)\cZ_{t-k,t-k-1}^\om(y_k,y_{k+1})dy_{k-1}}{\int_{\bT^{2d}}g(y_{k-1})\cZ_{t-k+1,t-k}^\om(y_{k-1},y_k)\cZ_{t-k,t-k-1}^\om(y_k,y_{k+1})dy_{k-1} dy_k},
\end{aligned}
\]
for some non-trivial $g\geq0$. Then it is enough to use the fact that
\[
\begin{aligned}
&\inf_{y_k,y_{k+1}\in\bT^d}\pi_k^\om(y_k\,|\,y_{k+1})\\
&\geq  \frac{\int_{\bT^d} g(y)dy\cdot \inf_{y_{k-1},y_k,y_{k+1}\in\bT^d}\cZ_{t-k+1,t-k}^\om(y_{k-1},y_k)\cZ_{t-k,t-k-1}^\om(y_k,y_{k+1})}{\int_{\bT^d} g(y)dy\cdot \sup_{y_{k-1},y_{k+1}}\int_{\bT^d}\cZ_{t-k+1,t-k}^\om(y_{k-1},y_k)\cZ_{t-k,t-k-1}^\om(y_k,y_{k+1})dy_k }
\end{aligned}
\] 
and the definition of $B_t(\delta)$ in \eqref{e.defB} to complete the proof.
\end{proof}

Combining Lemmas~\ref{l.lowerbd} and \ref{l.lbdtr}, we conclude that
there exists $\delta>0$, depending  only on $d,R$, and  events
\begin{equation}
\label{021903-21a}
B_{t-k}(\delta)\in \F_{[t-k-1,t-k+1]}, \quad\quad  k=2,\ldots,N-2,
\end{equation}
such that 
\begin{equation}
\label{021903-21}
\mathbf{P}[B_{t-k }(\delta)]> \delta,
\end{equation}
and the transition probability of the chain $\{Y_n\}_{n=1}^N$ satisfies
\begin{equation}
\label{021903-21b}
\inf_{y_k,y_{k+1}\in\bT^d}\pi_k^\om(y_k\,|\,y_{k+1})>\delta,\quad\quad\quad \mbox{ if }  \om
\in B_{t-k}(\delta),\quad k=2,\ldots,N-2.
\end{equation}

\subsection{Coupling construction}

\gu{The constructed Markov chain in the previous section satisfies a type of ``Doeblin condition'', expressed in the form of \eqref{021903-21b}. It is not the standard Doeblin condition which gives a lower bound on the transition density. Instead, the transition density is random itself and whether the Doeblin condition holds depends on if the random realization $\omega$ belongs to certain sets or not. Nevertheless, the coupling argument in the case of the Doeblin condition still applies, and} 
in this section, we use this
standard  argument to study the dependence of $u(t,x; \nu)$ on the ``remote'' environment $\{\xi(r,\cdot): r\leq s\}$ with $t-s\gg1$.

Let $t>1$. Define
\[
I_k:=[t-k-1,t-k+1],\quad\quad k=0,\ldots,[t]-1.
\]
We decompose the interval $[0,t]$ into
subintervals 
\[
[0,t]=[0,t-2m_N-1]\cup \bigcup_{m\in \Lambda_N} I_{m} \cup [t-1,t].
\]
where, as we recall $N=[t]+1$ and 
\begin{equation}\label{e.defLambdaN}
\Lambda_N:=\{2,4,\ldots, 2m_N\},
\end{equation} and $ 2m_N$ is the largest   number so that $t-2m_N-1\geq 0$:
\[
2m_N=\left\{\begin{array}{cc}
N-2 & \text{ if } N \text{ is even,}\\
N-3 & \text{ if } N \text{ is odd.}
\end{array}
\right.
\]
From \eqref{021903-21b} we know that there exists $\delta>0$ that only depends on $d,R$ such that  $B_{t-k }(\delta)\in\F_{I_k}$, and 
\[
\inf_{y_k,y_{k+1}\in\bT^d}\pi_k^\om(y_k\,|\,y_{k+1})>\delta,\quad\quad \mbox{ if }
\om\in B_{t-k}(\delta),\quad k\in\Lambda_N.
\] 
We fix $\delta$ from now on.
Let
\begin{equation}
\label{qk}
q_k^\om(y_k\,|\,y_{k+1}):=  \frac{\pi_k^\om(y_k\,|\,y_{k+1})-\delta}{1-\delta},\quad \quad \mbox{ if }\om\in B_{t-k}(\delta).
\end{equation}
Let us introduce a sequence of i.i.d. Bernoulli random variables
$\{\tau_k\}_{k\in\bbZ}$ that is defined on a probability space $\Big(\Sigma,{\cal
  A},\QQ\Big)$. We shall denote by $\E_{\QQ}$ the corresponding 
expectation.  {We assume
  that 
\[
\mathbf{Q}[\tau_k=1]=\delta, \quad\quad \mathbf{Q}[\tau_k=0]=1-\delta.
\]
}
Consider
the following random transition  kernels: for $k=2,3,\ldots, {N-2} $,  
 we have
\begin{align}
&\tilde \pi_k^{\om,0}\Big(y_k\,|\,y_{k+1}\Big)=q_k^\om(y_k|y_{k+1}),
  \quad \mbox{ if } \om \in B_{t-k}(\delta),
\label{ctp1}\\
&
\tilde \pi_k^{\om,1}\Big(y_k \,|\, 
y_{k+1}\Big)=1,\qquad \mbox{ if }\om \in
B_{t-k}(\delta), \label{ctp2}\\
&
\tilde \pi_k^{\om,\tau}\Big(  y_k \,|\, 
y_{k+1}\Big)=\pi_{k}^\om(y_k\,|\,y_{k+1}),\qquad  \mbox{ if }\om \in B_{t-k}^c(\delta),\tau\in\{0,1\}. \label{ctp3}
\end{align}
We also let

\begin{align}
&{\tilde{\pi}_{1}^{\om,\tau}\Big(  y_1\,|\, 
y_{2}\Big)=\pi_{1}^\om(y_{1}\,|\,y_{2}), \quad \tau\in\{0,1\}.}\\
&{\tilde{\pi}_{N-1}^{\om,\tau}\Big(  y_{N-1}\,|\, 
y_{N}\Big)=\pi_{N-1}^\om(  y_{N-1}\,|\, 
y_{N}), \quad \tau\in\{0,1\}.}\\
&\tilde \pi_{N}^{\om,\tau}\Big( d y_{N}\Big)=\pi_{N}^\om(dy_{N}),\quad \tau\in\{0,1\}. \label{ctpN}
\end{align}
 We shall write $\tilde \pi_{N}^{\om,\tau}(y; \nu)$  when we
  wish to 
  highlight the dependence of $\tilde \pi_{N}^{\om,\tau}$ on the
  initial data $\nu$, with the convention of writing $\tilde
  \pi_{N}^{\om,\tau}(y; f)$, when $\nu(dx)=f(x)dx$.
For a given realization of $\{\tau_k(\si)\}_{k=1}^N$, $\si\in \Sigma$, we 
consider the path measure $ \bbP^{\om,\si}_{\tilde{\pi}_N}$ on $\bT^{Nd}$ corresponding to the time
reversed 
Markov dynamics  obtained by using $\tilde \pi_N^{\om,\tau_N(\si)}$ as
the final density and $\tilde \pi_k^{\om,\tau_k(\si)}( y_k \,|\,
y_{k+1})$ as transition probabilities between times $k+1$ and $k$,
for $k=1,\ldots,N-1$. Let $\E^{\om,\si}_{\tilde{\pi}_N}$ be the corresponding expectation, for each fixed realization of $\omega,\sigma$.
Note that
\begin{equation}
\label{032204-21}
\int_{\Sigma}\bbP^{\om,\si}_{\tilde{\pi}_N}(A)\QQ(d\si)=\Pb^\omega_{\pi_N}(A),
\quad A\in {\cal B}(\bT^{Nd}). 
\end{equation}
Therefore, from \eqref{e.densityY1} and \eqref{032204-21},
\begin{equation}\label{e.densityY1c}
  u(t,x;\nu)=\frac{\int_{\Sigma}\QQ(d\si)\E_{\tilde \pi_N}^{\om,\si}[\cZ_{t,t-1}^\om(x,Y_1)]}{\int_{\bT^d}dx'\int_{\Sigma}\QQ(d\si)\E_{\tilde{\pi}_N}^{\om,\si}[\cZ_{t,t-1}^\om (x',Y_1)]}.
\end{equation}

Define the events
\begin{equation}
\label{022203-21z}
{A}_k^t=B_{t-k}(\delta)\times \{\tau_k=1\}\in {\cal
  F}_{[t-k-1,t-k+1]}\otimes {\cal A},\quad k=2,3,\ldots,N-2.
\end{equation}
If some $A_k^t$ occurs, the transition probability is given by the uniform distribution, as can be seen from \eqref{ctp2}, so the chain ``renews''. We emphasize here that $N$ depends on $t$ as well. If there is no danger of confusion we shall omit the superscript $t$
from the notation of the events.
Note that  
\begin{equation}
\label{012003-21}
\mathbf{P}\otimes\QQ[ {A}_k^t]=\mathbf{P}[B_{t-k}(\delta)]\mathbf{Q}[\tau_k=1]>\delta^2,\quad
 {k=2,3,\ldots,N-2}.
\end{equation}
It is also clear that the events $A_m^t$,  $m\in \Lambda_{N}$ are independent under $\mathbf{P}\otimes\QQ$ (recall that $\Lambda_N$ was defined in \eqref{e.defLambdaN}).
In addition, thanks to \eqref{e.defBa},   for any
$1<s<t$ and non-negative $\ell\in\bbZ$ we have
\begin{equation}
\label{032903-21}
\theta_{t-s,y}^{-1}\Big(A_k^s\Big)=A_{k}^t,\quad k=2,3,\ldots,N(s)-2,\,y\in\bT^d.
\end{equation}
Here we have denoted $\theta_{r,y}(\om,\si):=\Big(\theta_{r,y}(\om),\si\Big)$.

Note that for any $F\in B_b(\bT^d)$ we have
\begin{equation}
\label{022003-21}
\begin{aligned}
\E_{\QQ}\otimes \E_{\tilde \pi_N}^{\om,\cdot}
\Big[F(Y_1)\Big|\tau_k=1\Big]=& {\tfrac{1}{\delta}\int_{\{\tau_k=1\}}\QQ(d\si)\E_{\tilde \pi_N}^{\om,\si}[F(Y_1)]}\\
=&\E_{k}^{\om}
\Big[F(Y_1)\Big],\quad \quad 2\le k\le N-2,\,\om \in B_{t-k}(\delta),
\end{aligned}
\end{equation}
where, as we recall, $\E_k^\omega$ is the expectation with respect to the Markov chain constructed in \eqref{042003-21} and \eqref{012003-211}.
Now it is clear 
that for any $t>1$ and $2\leq k\leq N-2$,  the random variable
\begin{equation}
\label{022003-21a}
 1_{ B_{t-k}(\delta)}(\omega) \int_{\{\tau_k=1\}}\QQ(d\si)\E_{\tilde \pi_N}^{\om,\si}[\cZ_{t,t-1}^\om(x,Y_1)],\quad\quad  \om\in\Om
\end{equation}
is ${\cal F}_{ {[t-k-1,t]}}$-measurable.

%

\subsection{Approximation}
For any $\nu\in {\cal M}_1(\bT^d)$, define
\begin{equation}\label{e.defXom}
X^\omega(t,x; \nu):=\E_{ \pi_N(\nu)}^{\om}[\cZ_{t,t-1}^\om(x,Y_1)]=\E_{\QQ}\E_{\tilde \pi_N(\nu)}^{\om,\cdot}[\cZ_{t,t-1}^\om(x,Y_1)].
\end{equation}
When $\nu(dx) =f(x)dx$ for $f\in D(\bT^d)$ we shall write $X^\omega(t,x; f):=X^\omega(t,x; \nu)$. 
We omit writing the initial data $\nu$ from our notation when they are obvious from the context.

Fix $t>10$ and $2\leq k\leq  {N(t)-2}$. We shall approximate $X^\omega(t,x)$
by some random variable that is $\F_{[t-k,t]}-$measurable. 
Define 
\begin{equation}\label{e.defLambdaNk}
\Lambda_{N(t),k}=\{m: m \leq k-1\}\cap \Lambda_{N(t)},
\end{equation} 
\begin{equation}\label{e.deftildeAk}
\tilde{A}_k^t=\bigcup_{m\in \Lambda_{N(t),k}}A_m^t\in {\cal
  F}_{[t-k,t-1]}\otimes {\cal A},
\end{equation}
and
\begin{equation}\label{e.defXkom}
X_k^\omega(t,x):= \int_{\Sigma}\QQ(d\si)1_{\tilde{A}_k^t}(\omega,\sigma)\E_{\tilde \pi_{N(t)}}^{\om,\si}[\cZ_{t,t-1}^\om(x,Y_1)],
\end{equation}
and, thanks to
 \eqref{032903-21},
\begin{equation}\label{e.deftildeAk1}
\theta_{t-s,0}^{-1}\Big(\tilde{A}_k^s\Big)=\tilde{A}_{k}^t.
\end{equation}
In particular, as can be seen from \eqref{022003-21}, the random
  field $\{X_k(t,x)\} $ does not depend on the
  initial data.
Moreover, thanks to \eqref{022903-21} and \eqref{e.deftildeAk1}, we have
\begin{equation}\label{XtXs}
X_k^{\theta_{t-s,0}(\omega)}(s,x)=X_k^\omega(t,x),\quad k\le N(s)-2.
\end{equation}
As a direct consequence of the above and   formula \eqref{022003-21}, 
we conclude:
\begin{lemma}
\label{prop012303-21}
Suppose that $\nu_1,\nu_2\in D(\bT^d)$.
For any $t\ge s>k$ and $2\leq k\leq N(s)-2$,  
the laws of the fields
$X_k(t,\cdot\,;\nu_1)$ and $X_k(s,\cdot\,;\nu_2)$ over the space
$L^1(\bT^d)$ are identical.
\end{lemma}


\begin{lemma}\label{l.appX}
For any   $2\leq k\leq N(t)-2$ and $\nu\in {\cal M}_1(\bT^d)$,
the field $X_k (t,\cdot;\nu)$ is
$\F_{[t-k,t]}-$measurable. Moreover, for any $p\in[1,+\infty)$, there exists
$C,\lambda>0$ only depending on $p,d$ and $R$ such that 
\begin{equation}\label{e.3211}
\EE\Big[\sup_{\nu\in{\cal M}_1(\bT^d)}\|X^\omega (t;\nu)-X_k^\omega (t;\nu)\|^p_{L^\infty(\bT^d)}\Big] \leq
Ce^{-\lambda k}\quad \mbox{for   $2\leq k\leq {N(t)-2}$}.
\end{equation}
\end{lemma}
\begin{proof}
First we obviously have
\begin{equation}
\E_{\tilde \pi_N}^{\om,\si}[\cZ_{t,t-1}^\om(x,Y_1)]\leq \sup_{x,y\in\bT^d} \cZ_{t,t-1}^\omega(x,y),
\end{equation} 
which implies
\[
\sup_{x\in\bT^d}|X^\omega(t,x;\nu)-X_k^\omega(t,x;\nu)| \leq \sup_{x,y\in\bT^d} \cZ_{t,t-1}^\omega(x,y)\int_{\Sigma}\QQ(d\sigma) 1_{(\tilde{A}_k^t)^c}(\omega,\sigma).
\]
Thus, by Jensen's inequality
\begin{align*}
&\EE\Big[\sup_{\nu\in{\cal
  M}_1(\bT^d)}\|X^\omega(t;\nu)-X_k^\omega(t;\nu)\|^p_{L^\infty(\bT^d)}\Big]\\
&
 \leq \sqrt{\EE \Big(\sup_{x,y\in\bT^d} \cZ_{t,t-1}^\omega(x,y)\Big)^{2p} }\sqrt{\mathbf{P}\otimes\QQ[(\tilde{A}_k^t)^c]}.
\end{align*}
Since $(\tilde{A}_k^t)^c=\bigcap_{m\in\Lambda_{N,k}} (A_m^t)^c$ and the events
$A_m^t$, $m\in\Lambda_{N,k}$ are independent, from \eqref{012003-21} we have 
\begin{equation}
\label{022203-21}
\mathbf{P}\otimes\QQ[(\tilde{A}_k^t)^c]=\prod_{m\in \Lambda_{N,k}} \mathbf{P}\otimes\QQ[(A_m^t)^c]<(1-\delta^2)^{\#\Lambda_{N,k}},
\end{equation}
where $\delta$ depends only on $R$ and $d$.
Combining \eqref{022203-21} with Lemma~\ref{l.ulbd}, we derive that 
\[
\EE\Big[\sup_{\nu\in{\cal M}_1(\bT^d)}\|X^\omega(t;\nu)-X_k^\omega(t;\nu)\|^p_{L^\infty(\bT^d)}\Big] \leq  C(1-\delta^2)^{\#\Lambda_{N,k}/2},
\]
where the constant $C >0$ does not depend on $t$.
This proves \eqref{e.3211}. 

We prove now that $X_k^\omega(t,x;\nu)$ is $\F_{[t-k,t]}-$measurable.
Consider now the event $\tilde{A}_k$ given by
\eqref{e.deftildeAk}. {From this point on we shall drop $t$ from the
notation of the events and integers.} Since $N,k$ are both fixed here, to simplify the
notations, let 
$$
M:=\max\{m: m\in \Lambda_{N,k}\}
$$ and for any $m\in \Lambda_{N,k}$,  define 
\[
D_m=\left\{\begin{array}{ll}
A_m\cap\left(  \mathlarger{\mathlarger{\bigcup}}\limits_{m'>m,\,m'\in \Lambda_{N,k}} A_{m'}\right)^c,  & \text{ if  } m<M,\\
 A_M, & \text{ if } m=M.
 \end{array}
 \right. 
 \]
Then,
\begin{equation}
\label{042203-21}
1_{\tilde{A}_k}(\omega,\sigma)=\sum_{m\in \Lambda_{N,k}}1_{D_m}(\omega,\sigma),
\end{equation}
 and it suffices to show that for any $m\in \Lambda_{N,k}$, 
 \[
  \int_{\Sigma}\QQ(d\si)1_{D_m}(\omega,\sigma)\E_{\tilde \pi_N}^{\om,\si}[\cZ_{t,t-1}^\om(x,Y_1)] 
    \]
    is $\F_{[t-k,t]}-$measurable. We can write
        \[
    \begin{aligned}
  &\int_{\Sigma}\QQ(d\si)1_{D_m}(\omega,\sigma)\E_{\tilde \pi_N}^{\om,\si}[\cZ_{t,t-1}^\om(x,Y_1)] \\
  &= \int_{\Sigma}\QQ(d\si)1_{A_m}(\omega,\sigma)
\left(\prod_{m'>m,m'\in \Lambda_{N,k}}
  1_{A_{m'}^{c}}(\omega,\sigma)\right)\E_{\tilde
  \pi_N}^{\om,\si}[\cZ_{t,t-1}^\om(x,Y_1)] .
    \end{aligned}
    \]  
  By convention we understand a product over an empty set to be equal to $1$.  Invoking \eqref{022003-21}, the right hand side of the above equation can be written as 
    \begin{align}
\label{052203-21}
&
\left(\int_{\Sigma}\QQ(d\si)1_{A_m}(\omega,\sigma)
\E_{ m}^{\om}[\cZ_{t,t-1}^\om(x,Y_1)]\right)
 \int_{\Sigma}\QQ(d\si)\left(\prod_{m'>m,m'\in \Lambda_{N,k}}
  1_{A_{m'}^{c}}(\omega,\sigma)\right)\notag\\
&
   = 
\E_{ m}^{\om}[\cZ_{t,t-1}^\om(x,Y_1)] \delta 1_{B_{t-m}(\delta)}(\omega)\prod_{m'>m,m'\in
  \Lambda_{N,k}}\Big(1-\delta 1_{B_{t-m'}(\delta)}(\omega)\Big),
\end{align}
which is $\F_{[t-k,t]}-$measurable. This completes the proof.
\end{proof}

From the above proof, we can actually conclude that
$X_k^{\omega}(t,x)$ only depends on $t,x,k$ and the random environment
$\{\xi(s,\cdot):s\in [t-k,t]\}$. Recall that $u(t,\cdot;\nu)$ denotes
the endpoint density of the polymer path, assuming that the starting distribution is an arbitrary probability measure $\nu\in \mathcal{M}_1(\bT^d)$. We have the following key proposition which underlies Theorem~\ref{t.geomeEr}.
\begin{proposition}\label{p.couplingCon}
For any  $2\leq k\leq N(t)-2$, there exists $  u_k(t,x)$ that is \linebreak $\F_{[t-k,t]}-$measurable and does not depend on the initial  distribution $\nu$, such that for any $p\in[1,+\infty)$,
\begin{equation}
\label{072203-21}
\sup_{\nu\in{\cal M}_1(\bT^d)}\EE\Big[\|  u(t;\nu)- u_k(t)\|^p_{L^\infty(\bT^d)}\Big] \leq Ce^{-\lambda k},
\end{equation}
where the constants $C,\lambda>0$ depend only  on $p,d,R$.
\end{proposition}
\begin{proof}
By \eqref{e.densityY1c} and \eqref{e.defXom}, we have 
\[
  u(t,x;\nu)=\frac{X^\omega(t,x;\nu)}{\int_{\bT^d}X^\omega(t,x';\nu)dx'}.
\]
Let 
$$
C_k(\delta):=\bigcup_{m\in\Lambda_{N,k}}B_{t-m}(\delta),
$$
where we recall that $\Lambda_{N,k}$ was defined in \eqref{e.defLambdaNk}. 
 {Assuming $\omega\in B_{t-m}(\delta)$ for some $m\in \Lambda_{N,k}$, then $1_{\tilde{A}_k}(\omega,\sigma)\geq 1_{\{\tau_m=1\}}(\sigma)$ ($\tilde{A}_k$ was defined in \eqref{e.deftildeAk}), which implies 
\begin{equation}
\label{062203-21}
\begin{aligned}
\inf_{x\in\bT^d} X_k^\omega(t,x)\geq &\left(\int_{\Sigma}\QQ(d\si)1_{\{\tau_m=1\}}(\sigma)\right)\inf_{x,y\in\bT^d} \cZ_{t,t-1}^\omega(x,y)\\
&=\delta\inf_{x,y\in\bT^d} \cZ_{t,t-1}^\omega(x,y),\quad \quad \mbox{ if } \om\in C_k(\delta),
\end{aligned}
\end{equation}}
with $X_k^\omega$ defined in \eqref{e.defXkom}.
Define 
\begin{equation}\label{e.defuk}
  u_k(t,x):=
  \left\{
  \begin{array}{ll}
  \frac{X^\omega_k(t,x)}{\int_{\bT^d}X^\omega_k(t,x')dx'}, &\mbox{ if }\om\in C_k(\delta)\\
  0, &\mbox{ if } \omega\in C_k^c(\delta).
  \end{array}
  \right.
\end{equation}
From
the definition \eqref{e.defXom},  we conclude
\begin{equation}
  \label{X011011-22}
\inf_{x,y\in\bT^d}\cZ^\omega_{t,t-1}(x,y)\leq\inf_{x\in\bT^d} X^\omega(t,x;\nu)\leq \sup_{x\in\bT^d}X^\omega(t,x;\nu) \leq \sup_{x,y\in\bT^d} \cZ^\omega_{t,t-1}(x,y).
\end{equation}
Concerning  $X^\omega_k(t,x)$, we conclude that
\begin{align}
   \label{X021011-22}
\delta 1_{C_k(\delta)}(\om)\inf_{x,y\in\bT^d}\cZ^\omega_{t,t-1}(x,y)&\le  \inf_{x\in\bT^d}
  X^\omega_k(t,x)\notag\\
&
\leq \sup_{x\in\bT^d}X^\omega_k(t,x) \leq \sup_{x,y\in\bT^d} \cZ^\omega_{t,t-1}(x,y) .
\end{align}
\tk{The first inequality follows from \eqref{062203-21},  while the last one is a direct consequence of \eqref{e.defXkom}.}

\tk{Thus, for $\om\in C_k(\delta)$, we have (cf \eqref{e.defuk})
\begin{align*}
   \sup_{x\in\bT^d}|  u(t,x;\nu)- &u_k(t,x)|\le \frac{|X^\omega(t,x;\nu)-
  X^\omega_k(t,x)|}{\int_{\bT^d}X^\omega(t,x';\nu)dx'}\\
  &+X^\omega_k(t,x)\frac{\int_{\bT^d}|X^\omega_k(t,x')-X^\omega(t,x';\nu)|dx'}{\Big(\int_{\bT^d}X^\omega_k(t,x')dx'\Big) \Big(\int_{\bT^d}X^\omega(t,x';\nu)dx'\Big)}.
\end{align*}
Using \eqref{X011011-22} and \eqref{X021011-22}, we conclude  }
\begin{equation}
\label{082203-21}
\begin{split}
\sup_{x\in\bT^d}|  u(t,x;\nu)-  u_k(t,x)|\leq& \sup_{x\in\bT^d}
|X^\omega(t,x;\nu)-X_k^\omega(t,x)|\\
&
\times \left(\inf_{x,y\in\bT^d}\cZ^\omega_{t,t-1}(x,y) \right)^{-1}\left(1+\frac{1}{\delta}\cdot
    \frac{\sup_{x,y\in\bT^d}
      \cZ^\omega_{t,t-1}(x,y)}{\inf_{x,y\in\bT^d}\cZ^\omega_{t,t-1}(x,y)}\right)
\end{split}
\end{equation}
for $\om\in C_k(\delta)$.
We can therefore write
\begin{equation}
\label{072203-21a}
\EE[\sup_{x\in\bT^d} |  u(t,x;\nu)-  u_k(t,x)|^p] =I_1+I_2,
\end{equation}
where
\begin{equation}
\label{072203-21b}
\begin{split}
&
I_1:=\EE\Big[\sup_{x\in\bT^d} |  u(t,x;\nu)-  u_k(t,x)|^p1_{C_k(\delta)}\Big],\\ 
&
I_2:=\EE\Big[\sup_{x\in\bT^d} |  u(t,x;\nu)|^p1_{C_k^c(\delta)}\Big].
\end{split}
\end{equation}
To estimate $I_1$, we use the bound \eqref{082203-21} together with 
  Lemmas~\ref{l.ulbd} and \ref{l.appX}, and, since $\delta$ only depends on $d$ and $R$, we obtain that 
 there exist $C,\lambda>0$ depend only  on $p,d,R$ such that $I_1\leq Ce^{-\lambda k}$. 
To deal with  $I_2$, note that by the H\"older estimate for any $r>p$
we can write
\begin{equation}
\label{072203-21d}
\begin{aligned}
I_2\leq &\left\{\EE\Big[\left(\sup_{x\in\bT^d}
    u(t,x;\nu)\right)^r\Big]\right\}^{p/r}\Big(\PP[C_k^c(\delta)]\Big)^{1-p/r}\\
  &= \left\{\EE\Big[\left(\sup_{x\in\bT^d}
    u(t,x;\nu)\right)^r\Big]\right\}^{p/r}\Big(\PP[\bigcap_{m\in\Lambda_{N,k}}B^c_{t-m}(\delta)]\Big)^{1-p/r}.
  \end{aligned}
\end{equation}
Since $B^c_{t-m}(\delta)$ are independent and satisfy
\eqref{021903-21}, we conclude that the right hand side can be
estimated by 
$$
\left\{\EE\Big[\left(\sup_{x\in\bT^d}
    u(t,x;\nu)\right)^r\Big]\right\}^{p/r}\Big(1-\delta\Big)^{(1-p/r)\#\Lambda_{N,k}}.
$$
Using Lemma \ref{l.ulbd} we conclude that
 there exist $C,\lambda>0$ depend only  on $p,d,R$ such that $I_2\leq Ce^{-\lambda k}$.
The proof is complete.
\end{proof}


\subsection{Proof of Theorem~\ref{t.geomeEr}}

Recall that given an initial distribution $\nu\in {\cal M}_1(\bT^d)$
the polymer endpoint density $u(t;\nu)$ is given   by \eqref{e.defPo}. 
Fix $ F\in
{\rm Lip}({\cal M}_1(\bT^d)) $. 
 By virtue
of Lemma~\ref{prop012303-21} and \eqref{e.defuk}, for $t>s>k>s/2$, and
any bounded and measurable function $F$ we have
\begin{equation}
  \label{X031011-22}
  \EE[F(  u_k(t))]=\EE[F(  u_k(s))].
\end{equation}
Therefore, for any $\nu_1,\nu_2\in \mathcal{M}_1(\bT^d)$, we have
\begin{align*}
&
\Big|{\cal P}_tF(\nu_1)-{\cal P}_sF(\nu_2)\Big|=\Big|
  \EE[F(u(t; \nu_1))]-\EE[F(u(s; \nu_2))]\Big|\\
&
=\Big|\EE[F(u(t; \nu_1))]-\EE[F(  u_k(t))]-\EE[F(u(s; \nu_2))]+\EE[F(
  u_k(s))]\Big|\\
&
\le
\Big|\EE[F(  u(t; \nu_1))]-\EE[F(  u_k(t))]\Big|+\Big|\EE[F(  u(s;
  \nu_2))]-\EE[F(  u_k(s))]\Big|\\
&
\le
  \|F\|_{{\rm Lip}}\Big(\sup_{\mu\in {\cal M}_1(\bT^d)}\EE\|  u(t;
  \mu)-  u_k(t)\|_{L^1(\bT^d)}+\sup_{\mu\in {\cal
  M}_1(\bT^d)}\EE\|  u(s; \mu)-  u_k(s)\|_{L^1(\bT^d)}\Big).
\end{align*}
Using Proposition \ref{p.couplingCon}
 we can further write   that 
\begin{equation}
\label{012304-21}
\Big|{\cal P}_tF(\nu_1)-{\cal P}_sF(\nu_2)\Big|\leq C\|F\|_{{\rm
    Lip}}e^{-\lambda k}\le C\|F\|_{{\rm
    Lip}}e^{-\lambda s/2},
\end{equation}
which in turn implies that there exist $C,\lambda>0$ such that
\[
\sup_{\nu_1,\nu_2\in {\cal M}_1(\bT^d)}\Big\|\delta_{\nu_1}{\cal P}_t-\delta_{\nu_2}{\cal P}_s\Big\|_{\mathrm{FM}}\leq
C e^{-\lambda s/2},\quad t>s>0.
\]
\gu{Here we used the  Fortet-Mourier metric, see \eqref{FM} below.}
For any $\Pi\in {\rm M}_1( {\cal M}_1(\bT^d))$ and $F\in {\rm Lip}({\cal
  M}_1(\bT^d))$, we can write
$$
\int_{{\cal M}_1(\bT^d)} F(\nu)\Pi{\cal P}_t(d\nu)=\int_{{\cal M}_1(\bT^d)}\Pi(d\nu) \left(\int_{\mathcal{M}_1(\bT^d)}\delta_{\nu}{\cal P}_t(d\mu)F(\mu)\right),
$$
so we also conclude that
{\begin{equation}
\label{042303-21}
\sup_{\Pi_1,\Pi_2\in {\rm M}_1( {\cal M}_1(\bT^d))}\Big\|\Pi_1{\cal P}_t-\Pi_2{\cal P}_s\Big\|_{\mathrm{FM}}\leq
C e^{-\lambda s/2},\quad t>s>0,
\end{equation}}
which implies that there
exists a unique $\pi_\infty\in {\rm M}_1( {\cal M}_1(\bT^d))$ such that
\begin{equation}\label{e.452}
\lim_{t\to\infty}\Pi{\cal P}_t=\pi_\infty
\end{equation}
for any $\Pi\in  {\rm M}_1( {\cal M}_1(\bT^d))$, in the
Fortet-Mourier metric on ${\rm M}_1( {\cal
  M}_1(\bT^d))$. \gu{We also have
$
\pi_\infty =\pi_\infty{\cal P}_t $  for all $t>0$,
see a rather standard proof in Appendix \ref{appC} below.}
Suppose that $p\in[1,+\infty]$. Since $\delta_{\nu}{\cal
  P}_t(D^p(\bT^d))=1$ for all $\nu\in {\cal M}_1(\bT^d)$ and any $t>0$, 
 we have
$$
\pi_\infty\Big(D^p(\bT^d)\Big)=\pi_\infty{\cal P}_t\big(D^p(\bT^d)\Big)=\int_{{\cal M}_1(\bT^d)}\pi_\infty(d\nu) \left(\delta_{\nu}{\cal
  P}_t\Big(D^p(\bT^d)\Big)\right)=1,
$$
which finishes the proof of part (i).

%
Part (ii) is a direct consequence of \eqref{042303-21} and \eqref{e.452}.

Suppose now that   $p\in(1,\infty]$ and  $F\in {\rm
    Lip}(D^p(\bT^d))$. Proposition \ref{p.couplingCon} also applies in
  this case with the $L^1(\bT^d)$ norm, used above, replaced
  by the $L^p(\bT^d)$ norm. In consequence we conclude
  \eqref{012304-21} with the Lipschitz norm of $F$ in the respective
  $L^p(\bT^d)$ metric. Hence, we  infer the estimate
  \eqref{042303-21} with the respective Fortet-Mourier norm on the l.h.s. The remaining part of the argument stays the same and, as
  a result, we conclude the statement of the theorem in the case of the
  space $D^p(\bT^d)$.
\qed

\subsection{Asymptotic stability of  the transition probability operator}

A straightforward consequence of Theorem~\ref{t.geomeEr} is 
\begin{corollary}
\label{cor020104-21}
Suppose that $F\in C_b ({\cal M}_1(\bT^d))$ and $\Pi\in {\rm
  M}_1({\cal M}_1(\bT^d))$  then
\begin{equation}
\label{040104-21}
\lim_{t\to+\infty}\int_{{\cal M}_1(\bT^d)}{\cal
    P}_tFd\Pi= \int_{{\cal M}_1(\bT^d)}F(u)\pi_\infty(du).
\end{equation}
\end{corollary}

\begin{proposition}
\label{prop021203-21}
Suppose that $F\in {\rm Lip} ({\cal M}_1(\bT^d))$. Then,
\begin{equation}
\label{021203-21}
\lim_{t\to\infty}\EE_{\Pi}\left(\frac{1}{t} \int_0^t
  F(u(s))ds-\int_{{\cal M}_1(\bT^d)}F(u)\pi_\infty(du) \right)^2=0. 
\end{equation}
where $\EE_{\Pi}$ denotes the expectation w.r.t. the noise $\xi$ realization and the initial
data $u(0)$ distributed according to a Borel probability 
measure $\Pi$ on ${\cal M}_1(\bT^d)$. {The same result also holds for $D^p(\bT^d)$ in place of ${\cal M}_1(\bT^d)$ for any $p\in[1,+\infty]$.}
\end{proposition}
\proof
\gu{To prove \eqref{021203-21}, } it suffices to show that as $t\to\infty$,
\begin{equation}
\label{021203-21a}
\EE_{\Pi}\left(\frac{1}{t} \int_0^t
  F(u(s))ds\right)\to \int_{{\cal M}_1(\bT^d)}F(u)\pi_\infty(du),
\end{equation}
and
\begin{equation}
\label{021203-21b}
\EE_{\Pi}\left(\frac{1}{t} \int_0^t
  F(u(s))ds\right)^2\to\left(\int_{{\cal M}_1(\bT^d)}F(u)\pi_\infty(du)\right)^2. 
\end{equation}
Equality
\eqref{021203-21a} follows directly from \eqref{040104-21}. Concerning
\eqref{021203-21b}, 
we can write
\begin{equation}
\label{021203-21c}
\begin{split}
&\EE_{\Pi}\left(\frac{1}{t} \int_0^t
  F(u(s))ds\right)^2=\frac{2}{t^2} \EE_{\Pi}\left(\int_0^t
  F(u(s))ds\int_0^{s}F(u(s'))ds' \right)\\
&
 {=\frac{2}{t^2} \int_0^t
  ds \int_0^{s}ds' \left[\int_{{\cal M}_1(\bT^d)} {\cal
      P}_{s'}\Big(F{\cal P}_{s-s'}F
      \Big)d\Pi \right]}\\
      &=\frac{2}{t^2} \int_0^t
  ds \int_0^{s}ds' \left[\int_{{\cal M}_1(\bT^d)} {\cal
      P}_{s-s'}\Big(F{\cal P}_{s'}F
      \Big)d\Pi \right].
\end{split}
\end{equation}
Using Theorem~\ref{t.geomeEr} part (ii), we conclude that
\[
\Big\|{\cal
      P}_{s-s'}\Big(F{\cal P}_{s'}F
      \Big)-{\cal
      P}_{s-s'}\Big(F\int_{{\cal M}_1(\bT^d)} F(u)\pi_\infty(du)\Big)
     \Big\|_\infty\le Ce^{-\lambda s'}.
      \]
      Therefore, 
  \begin{equation}
\label{021203-21cc}
\begin{split}
& \left|\frac{2}{t^2}\int_0^t
  ds \int_0^{s}ds' \left[\int_{{\cal M}_1(\bT^d)} {\cal
      P}_{s-s'}\Big(F{\cal P}_{s'}F
      \Big)d\Pi \right]\right.\\
      &
      \left.-\frac{2}{t^2}\left(\int_{{\cal M}_1(\bT^d)} Fd\pi_\infty\right)
 \int_0^t
  ds \int_0^{s}ds' \int_{{\cal M}_1(\bT^d)} {\cal
      P}_{s-s'}Fd\Pi\right|\\
      &
      \le\frac{C}{t^2} \int_0^t ds\int_0^s e^{-\lambda s'}ds'\to0, \quad\quad\quad \mbox{ as } t\to\infty.
      \end{split}
\end{equation}    
On the other hand, upon another
application of  Theorem~\ref{t.geomeEr} part (ii), we conclude that
\begin{equation}
\label{021203-21d}
\begin{aligned}
&  \frac{2}{t^2} \int_0^t
  ds \int_0^{s}ds' \int_{{\cal M}_1(\bT^d)} {\cal
      P}_{s-s'}Fd\Pi\\
      &=2 \int_0^1
  ds \int_0^{s}ds' \int_{{\cal M}_1(\bT^d)} {\cal
      P}_{t(s-s')}Fd\Pi\\
      &\to2\int_{{\cal M}_1(\bT^d)} F(u)\pi_\infty(du)  \int_0^1
  ds \int_0^{s}ds'=\int_{{\cal M}_1(\bT^d)} F(u)\pi_\infty(du).
  \end{aligned}
\end{equation}
Combining \eqref{021203-21cc} with \eqref{021203-21d}, we conclude  \eqref{021203-21b}, thus finishing the proof of the proposition. The proof for the case of $F\in {\rm Lip} (D^p(\bT^d))$ is the same.
\qed

\subsection{More on Theorem~\ref{t.geomeEr}}

In the statement of Theorem~\ref{t.geomeEr}, we have only considered
the Lipschitz functional on $\mathcal{M}_1(\bT^d)$, or $D^p(\bT^d)$.
In what follows we shall need the stability of the semigroup ${\cal
  P}_t$ on functionals that are only local Lipschitz.  For that
purpose we shall need  the following lemma:
\begin{lemma}\label{l.bdut1}
For any $p\in[1,+\infty)$, there exists $C=C(d,R,p)$ such that \[
\sup_{t>1,\nu\in\mathcal{M}_1(\bT^d)}\EE\Big[\big(\sup_{x\in\bT^d}u(t,x;\nu)\big)^p\Big]\leq C.
\]
\end{lemma}
\begin{proof}
By \eqref{e.densityY1}, we have when $t>1$ that 
\[
\sup_{x\in\bT^d}u(t,x;\nu) \leq \frac{\sup_{x,y\in\bT^d} \cZ_{t,t-1}^\omega(x,y)}{\inf_{x,y\in\bT^d} \cZ_{t,t-1}^\omega(x,y)}, 
\]
then the result is a direct consequence of Lemma~\ref{l.ulbd}.
\end{proof}

Now consider the functional $\cR: D^2(\bT^d)\to \R$ defined in \eqref{fR}:
\begin{equation}
\label{022304-21}
\cR(v)=\int_{\bT^{2d}} R(x-y) v(x)v(y)dxdy\quad\mbox{ and  } \quad \tilde{\cR}=\cR-2\gamma,
\end{equation}
with $\gamma$ given by \eqref{e.defgamma}. \tk{The functional is globally Lipschitz
  in the case when $R(\cdot)$ is bounded, as then
\begin{align*}
  |\cR(v_2)-\cR(v_1)|&\le
  \|R\|_\infty\Big(\|v_2\|_{L^1(\bT^d)}+\|v_1\|_{L^1(\bT^d)}\Big)
  \|v_2-v_1\|_{L^1(\bT^d)}\\
  &
    \leq  2\|R\|_\infty
  \|v_2-v_1\|_{L^2(\bT^d)},\quad v_1,v_2\in D^2(\bT^d).
  \end{align*}
It is
   only locally Lipschitz
when $R=\delta$, as then $\cR(v)=\|v\|_{L^2(\bT^d)}^2$.} 

\begin{proposition}\label{p.tailPtR}
There exist $C,\lambda$ only depending on $R,d$ such that 
\[
\sup_{t>1,\nu\in \mathcal{M}_1(\bT^d)}|\mathcal{P}_t\tilde{\cR}(\nu)|\leq
Ce^{-\lambda t}.
\]
\end{proposition}
\begin{proof}
By Lemma~\ref{l.bdut1} and the Markov property, it suffices to prove the result for any $v\in D^2(\bT^d)$ and $t>1$. First we note that $\mathcal{R}(\cdot)$ is continuous from $D^2(\bT^d)$ to $\R$, therefore by Theorem~\ref{t.geomeEr}, we have as $t\to\infty$ that 
\[
\cR(u(t;v))\Rightarrow \cR(\tilde{u}) \quad\quad \mbox{ in distribution}, 
\]
with $\tilde{u}$ sampled from $\pi_\infty$. Lemma~\ref{l.bdut1} ensures the uniform integrability so we have
\begin{equation}\label{e.461}
\cP_t\cR(v)=\EE[\cR(u(t;v))]\to \EE[\cR(\tilde{u})]=2\gamma.
\end{equation}
Next, \tk{using \eqref{X031011-22},} for any $t_1>t_2>k>t_2/2$, we have
\[
\begin{aligned}
|\cP_{t_1}\cR(v)-\cP_{t_2}\cR(v)|
\leq \sum_{j=1}^2|\EE[\cR(u(t_j;v))-\cR(u_k(t_j))]|,
\end{aligned}
\]
with $u_k$ defined in \eqref{e.defuk}. By Proposition~\ref{p.couplingCon} and the fact that 
\[
\begin{aligned}
&|\cR(u(t_j;v))-\cR(u_k(t_j))|\\
&\leq
\|u(t_j;v)-u_k(t_j)\|_{L^\infty(\bT^d)}\Big(\|u(t_j;v)\|_{L^1(\bT^d)}+\|u_k(t_j)\|_{L^1(\bT^d)}\Big)\\
&\leq  2\|u(t_j;v)-u_k(t_j)\|_{L^\infty(\bT^d)},
\end{aligned}
\]
\tk{(here we have used the fact that $\int_{\bT^d}R(x)dx=1$)} we derive $|\cP_{t_1}\cR(v)-\cP_{t_2}\cR(v)|\leq  Ce^{-\lambda t_2}$. Sending $t_1\to\infty$ and applying \eqref{e.461}, we complete the proof.
\end{proof}

\tk{\begin{proposition}
  \label{prop012211-22}
  For any $p\in[1,+\infty)$ we have
  \[
\int_{{\cal M}_1(\bT^d)} \|u\|_{L^\infty(\bT^d)}^p\pi_\infty(du)<\infty.
\]
  \end{proposition}
  \proof Fix any $\nu\in {\cal M}_1(\bT^d)$. Thanks to \eqref{e.452},
  we have $\lim_{t\to\infty}\delta_{\nu}{\cal P}_t=\pi_\infty$, in the
Fortet-Mourier metric on ${\rm M}_1( {\cal
  M}_1(\bT^d))$. By Lemma \ref{l.bdut1}, we have
 \[
\limsup_{t\to\infty}\int_{{\cal M}_1(\bT^d)} \|u\|_{L^\infty(\bT^d)}^p \delta_{\nu}{\cal P}_t (du)<+\infty.
\]
The conclusion of the proposition then follows from an application of
\cite[Theorem 3.5, p. 31]{billingsley}.
  \qed}

\section{Proofs of Theorem~\ref{t.partitionfunction} and Theorem~\ref{t.mainth}}

\label{s.5}

Recall from \eqref{e.delogZt1} that 
\[
\log Z_t+\gamma t=\int_0^t\int_{\bT^d} u(s,y)\xi(s,y)dyds-\tfrac12\int_0^t \tilde{\cR}(u(s))ds
\]
where $\tilde{\cR}=\cR-2\gamma$, 
\begin{equation}\label{e.uU}
u(s)=u(s;\nu)=Z_s^{-1}\U(s;\nu),
\end{equation}
where $\nu$ is an arbitrary probability measure on $\bT^d$. By the definition of $\gamma$ in \eqref{e.defgamma}, we have 
\[
\int_{\mathcal{M}_1(\bT^d)} \tilde{\cR}(u) \pi_\infty(du)=0.
\]
The goal in this section is to prove Theorem~\ref{t.partitionfunction}: as $t\to\infty$, 
\begin{equation}\label{e.4201}
\frac{\log Z_t+\gamma t}{\sqrt{t}}\Rightarrow N(0,\sigma^2).
\end{equation}
With the above convergence, the proof of Theorem~\ref{t.mainth} goes as follows: for any $\nu\in\mathcal{M}_1(\bT^d)$ and $x\in\bT^d$, we write
\begin{equation}
\frac{\log \U(t,x;\nu)+\gamma t}{\sqrt{t}}=\frac{\log Z_t+\gamma t}{\sqrt{t}}+\frac{\log u(t,x;\nu)}{\sqrt{t}}.
\end{equation}
For $t>1$,  by \eqref{e.densityY1} we have
\begin{equation}\label{e.utUL}
\frac{\inf_{x,y\in\bT^d} \cZ_{t,t-1}^\omega(x,y)}{\sup_{x,y\in\bT^d} \cZ_{t,t-1}^\omega(x,y)}\leq u(t,x;\nu) \leq \frac{\sup_{x,y\in\bT^d} \cZ_{t,t-1}^\omega(x,y)}{\inf_{x,y\in\bT^d} \cZ_{t,t-1}^\omega(x,y)}.
\end{equation}
\tk{Using the elementary inequality  $\max\{\log r,0\}<r$, valid
for all $r>0$, 
we conclude that  
$$
|\log u(t,x;\nu)|\le u(t,x;\nu)+u^{-1}(t,x;\nu)
$$
Applying \eqref{e.utUL} and Lemma~\ref{l.ulbd},} we conclude therefore that there exists a constant $C>0$ such that  
$$
\EE|\log u(t,x;\nu)|\leq C\quad \mbox{ for all }
t>1,x\in\bT^d.
$$
 As a result we obtain
\[
\frac{\log u(t,x;\nu)}{\sqrt{t}}\to0,\quad\mbox{as }t\to+\infty,
\]
in probability and the proof of Theorem~\ref{t.mainth} is complete.  

\subsection{Construction of the corrector}

The  purpose of the present section is to construct the solution of
the equation
$-{\cal L}\chi=\tilde\cR$, where $\tilde\cR$ is defined in
\eqref{022304-21} and ${\cal L}$ is a properly understood
generator of the transition probability semigroup $(\cP_t)_{t\geq0}$, see
Remark \ref{rmk012304-21} below. The field $\chi:D^\infty(\bT^d)\to\bbR$ 
is called the  corrector and is a crucial object in the 
proof of the CLT for $\log Z_t$.
%
%
%
We start with the following lemma:
\begin{lemma}\label{l.bdRv}
For each  $p\in[1,+\infty)$ there exists a constant $C>0$ such that 
\begin{equation}\label{e.bduinfinity}
\sup_{t\geq0} \EE[\|u(t;v)\|_{L^\infty(\bT^d)}^p]\leq C(1+\|v\|_{L^\infty(\bT^d)}^p) ,\quad \tk{\mbox{for any }v\in D^\infty(\bT^d)}
\end{equation}
and 
\begin{equation}
\label{042304-21}
|{\cal P}_t \tilde{\cR}(v)|\leq  C(e^{-\lambda
  t}1_{t>1}+(1+\|v\|_{L^\infty(\bT^d)})1_{t\in[0,1]}),\quad \mbox{for any }v\in D^\infty(\bT^d),\,t\ge0
.
\end{equation}
\end{lemma}

\begin{proof}
The case of $t>1$ is implied by Lemma~\ref{l.bdut1} and Proposition~\ref{p.tailPtR}. For $t\in[0,1]$, we have
\[
\begin{aligned}
\cP_t\cR(v)&=\EE \int_{\bT^{2d}} R(x-y)u(t,x;v)u(t,y;v)dxdy\\
& \leq\EE \|u(t;v)\|_{L^\infty(\bT^d)} \int_{\bT^{2d}} R(x-y)u(t,x;v)dxdy =\EE \|u(t;v)\|_{L^\infty(\bT^d)},
\end{aligned}
\]
and it suffices to apply \eqref{e.she14} to conclude the proof.
\end{proof}


\begin{lemma}\label{l.conPRv}
For any $v\in D^\infty(\bT^d)$, the function $\mathcal{P}_t\tilde{\cR}(v)$ is continuous in $t\geq0$.
\end{lemma}

\begin{proof}
\tk{  Since
  $$
  \cP_t\cR(v)=\EE \int_{\bT^{2d}} R(x-y)u(t,x;v)u(t,y;v)dxdy,
  $$
  we can easily conclude the assertion of the lemma using  Lemma~\ref{l.uconti}, when $R(\cdot)$ is bounded.
 In the case of 
$R(\cdot)=\delta(\cdot)$ and $d=1$, then
$\cP_t\cR(v)=\EE\|u(t;v)\|_{L^2(\bT^d)}^2$ and the conclusion follows
from    Lemma~\ref{l.uconti} and an application of \eqref{e.bduinfinity}.}
\end{proof}
 

The time dependent corrector field is defined as  
\begin{equation}
\label{chit}
\chi(t,v):=\int_0^t\cP_s\tilde{\cR}(v)ds,\quad \quad t\ge0,\,v\in D^\infty(\bT^d).
\end{equation}

 

We have the following result.
\begin{proposition}
\label{prop032403-21}
The function $\chi:[0,+\infty)\times D^\infty(\bT^d)\to\R$ satisfies
\begin{equation}
\label{e.pt}
\frac{d\chi(t,v)}{dt}=\cP_t\tilde{\cR}(v),\quad t\ge0.
\end{equation}
In addition, for any $t\ge0$ and $v\in D^\infty(\bT^d)$, we have
\begin{equation}
\label{cellt}
\lim_{\delta\to0} \delta^{-1}[\cP_\delta \chi(t,v)-\chi(t,v)]  ={\cal P}_t \tilde{\cR}(v)-\tilde{\cR}(v).
\end{equation}
\end{proposition}

\begin{proof}
 First, \eqref{e.pt} is a consequence of Lemma~\ref{l.conPRv}. To show \eqref{cellt}, for any $\delta>0, t\geq0$ and $v\in D^\infty(\bT^d)$, we have
\[
\begin{aligned}
\delta^{-1}[\cP_\delta\chi(t,v)-\chi(t,v)]=&\delta^{-1} [\int_0^t{\cal P}_{s+\delta} \tilde{\cR}(v)ds-\int_0^t{\cal P}_s \tilde{\cR}(v)ds]\\
=&\delta^{-1}[\int_t^{t+\delta} {\cal P}_{s} \tilde{\cR}(v)ds-\int_0^\delta {\cal P}_s \tilde{\cR}(v)ds].
\end{aligned}
\]
Sending $\delta\to0$ and applying Lemma~\ref{l.conPRv} again, the proof is complete.
\end{proof}

By virtue of Lemma~\ref{l.bdRv}, we can define $\chi: D^\infty(\bT^d)\to \R$
\begin{equation}\label{e.defcorrector}
\chi(v):=\int_0^\infty \cP_t\tilde{\cR}(v) dt.
\end{equation}
 \begin{remark}
\label{rmk012304-21}
Analogously to \eqref{cellt}, we can show that 
\begin{equation}
\label{032304-21}
\lim_{\delta\to0} \delta^{-1}[\cP_\delta \chi(v)-\chi(v)]   =-\tilde{\cR}(v),\quad\quad  v
\in D^\infty(\bT^d).
\end{equation}
The left hand side of \eqref{032304-21} can be treated as a pointwise
definition of the generator ${\cal L}$  of the semigroup   $(\cP_t)_{t\geq0}$. 
\end{remark}
As a
 consequence of Lemma~\ref{l.bdRv}, the following proposition holds: 
\begin{proposition}
\label{prop022503-21}
We have for any $v\in D^\infty(\bT^d)$,
\begin{equation}\label{e.bdchi}
|\chi(v)|\leq C(1+\|v\|_{L^\infty(\bT^d)}),
\end{equation}
\begin{equation}
\label{cellt1}
\sup_{v\in D^\infty(\bT^d)} |\chi(v)-\chi(t,v)|\leq  \frac{C}{\lambda}e^{-\lambda
  t}, \quad\quad t\ge 1.
\end{equation}
\end{proposition}
\begin{proof}
It is a direct consequence of \eqref{042304-21}.
\end{proof}

\subsection{Frechet gradient of corrector}

Using the corrector $\chi$ defined in \eqref{e.defcorrector} we
shall be able to write $\log Z_t+\gamma t$ as a sum of a continuous square
integrable martingale and a term of order  $1$, see \eqref{ztNt} -
\eqref{e.4183} below. This fact  is  crucial in our proof of the CLT.
The quadratic variation of the martingale can be expressed in terms of
an appropriately defined   Frechet gradient of the corrector. The
present section is devoted to provide a definition of such an object. Note that
the latter  is not completely obvious, due to the fact that $\chi$ is
defined on a set $D^\infty(\bT^d)$ that has no interior points.

We extend the definition of $\cP_t\cR$ to an open subset of
$L^2(\bT^d)$ containing $D^\infty(\bT^d)$. 
For any $t\geq0$ and $v\in L^2(\bT^d)$ such that $\int_{\bT^{d}}\U(t,x;v)dx\neq 0$,  define 
\begin{equation}\label{e.defPtomega}
\mathscr{P}_t(v):=\int_{\bT^{2d}} R(x-y) \frac{\U(t,x;v)\U(t,y;v)}{\left(\int_{\bT^d}\U(t,x';v)dx'\right)^2}dxdy.
\end{equation}
Therefore,  $\cP_t\cR(\cdot)$ coincides with $\EE[\mathscr{P}_t(\cdot)]$ on $D^\infty(\bT^d)$, but the latter is defined on a  bigger set. Define
\begin{equation}\label{e.L2plus}
L_+^2(\bT^d):= \{f\in L^2(\bT^d): f\geq 0,\, \,\|f\|_{L^2(\bT^d)}\neq 0\}.
\end{equation}
For any $v\in L_+^2(\bT^d)$, we let
\begin{equation}\label{e.deftildev}
\tilde{v}(x)=\|v\|_{L^1(\bT^d)}^{-1}v(x). 
\end{equation}

For any $t\geq0$, $z\in\bT^d$, $v\in D^\infty(\bT^d)$ and $h\in
L^2(\bT^d)$, we let
\begin{equation}\label{e.defUz}
U(t,x;z,v)=\frac{\U(t,x;z)}{\int_{\bT^d} \U(t,x';v)dx'}
\end{equation}
and
\begin{equation}\label{e.defUh}
U(t,x;h,v)=\int_{\bT^d} U(t,x;z,v)h(z)dz
=\frac{\U(t,x;h)}{\int_{\bT^d} \U(t,x';v)dx'},\quad x\in\bT^d.
\end{equation}
We have the following lemma for the Frechet derivative of $\mathscr{P}_t(v)$.
\begin{lemma}\label{l.fredeomega}
For any $t\geq0, v\in D^\infty(\bT^d)$ and almost every realization of
$\xi$, $\mathscr{P}_t$ is twice Frechet differentiable at $v$ and the
respective first and second order Frechet derivatives equal
\begin{equation}\label{e.fredePomega}
\begin{aligned}
\cD\mathscr{P}_t(v)(z)=&\,2\int_{\bT^{2d}}R(x-y)U(t,x;z,v )u(t,y;v )dxdy\\
&-2\|U(t;z,v)\|_{L^1(\bT^d)}\int_{\bT^{2d}}R(x-y)u(t,x;v )u(t,y;v ) dxdy,
\end{aligned}
\end{equation}
\begin{equation}\label{e.frede2}
\begin{aligned}
&\cD^2\mathscr{P}_t(v)(z_1,z_2)\\
&=2\int_{\bT^{2d}} R(x-y)U(t,x;z_1,v)U(t,y;z_2,v)dxdy \\
&-4\|U(t;z_2,v)\|_{L^1(\bT^d)}\int_{\bT^{2d}} R(x-y)U(t,x;z_1,v)u(t,y;v)dxdy\\
&-4\|U(t;z_1,v)\|_{L^1(\bT^d)}\int_{\bT^{2d}} R(x-y)U(t,x;z_2,v)u(t,y;v)dxdy \\ 
&+6\|U(t;z_1,v)\|_{L^1(\bT^d)}\|U(t;z_2,v)\|_{L^1(\bT^d)}\int_{\bT^{2d}} R(x-y)u(t,x;v)u(t,y;v)dxdy.  
\end{aligned}
\end{equation}

\end{lemma}

\begin{proof}
We only consider $t>0$ here (the case of $t=0$ is easy to analyze).  One can write
\[
\|\U(t;v )\|_{L^1(\bT^d)}=\int_{\bT^{2d}} \U(t,x;y )v(y)dxdy,
\] where  $\U(t,x;y )$ is the Green's function of SHE, which, for almost all $\xi$,  is a continuous and positive function in $(x,y)\in\bT^{2d}$. Fix such a realization of $\xi$ and $v\in D^\infty(\bT^d)$, we have $\|\U(t;v )\|_{L^1(\bT^d)}>0$, thus, there exists $\delta>0$ (which could depend on $\xi$) such that $
\|\U(t;\tilde{v})\|_{L^1(\bT^d)}>0$ for any $\tilde{v}$ with $\|\tilde{v}-v\|_{L^2(\bT^d)}<\delta$. Then $\mathscr{P}_t(\cdot)$ is well-defined for all such $\tilde{v}$, and it is a straightforward calculation to check that for any $h\in L^2(\bT^d)$, as $\delta\to0$
\[
\delta^{-1}[\mathscr{P}_t(v+\delta h )-\mathscr{P}_t(v)]\to \la\cD\mathscr{P}_t(v),h\ra_{L^2(\bT^d)}.
\]
The proof of \eqref{e.frede2} is similar.
\end{proof}

Using the functional $\cR$, see \eqref{fR}, we can rewrite
\eqref{e.fredePomega} in the form
\begin{equation}\label{e.4121}
\begin{aligned}
\la \cD\mathscr{P}_t(v),h\ra_{L^2(\bT^d)}
=&2 \cR(U(t;h,v),u(t;v))\\
&-2\cR(u(t;v))\int_{\bT^d} U(t,x;h,v)dx,
\end{aligned}
\end{equation}
for any $v\in D^\infty(\bT^d), h \in L^2(\bT^d)$. 
By further restricting to $h\in L^2_+(\bT^d)$, with $\tilde{h}=h/\|h\|_{L^1(\bT^d)}$, we can write
\begin{equation}\label{e.4122}
\begin{aligned}
&\la \cD\mathscr{P}_t(v),h\ra_{L^2(\bT^d)}\\
&=2\int_{\bT^{3d}} R(x-y)[u(t,x;\tilde{h})-u(t,x;v)]u(t,y;v)U(t,x';h,v)dxdydx'\\
&=2\,\cR(u(t;\tilde{h})-u(t;v),u(t;v))\|U(t;h,v)\|_{L^1(\bT^d)}.
\end{aligned}
\end{equation}


Next, we prove several technical lemmas on the Frechet derivatives,
where we will only use the following estimate:  
\begin{equation}\label{e.conditionR}
|\cR(f,g)| \leq \hat R_*\|f\|_{L^2(\bT^d)}\|g\|_{L^2(\bT^d)},\quad f,g\in L^2(\bT^d),
\end{equation}
where 
\begin{equation}
\label{Rstar}
\hat R_*:=\sup_{k\in\bbZ^d}\hat{r}_k 
\end{equation} and $(\hat r_k)_{k\in\Z^d}$ are the Fourier coefficients of $R(\cdot)$, cf. \eqref{e.defrk}. By our assumption of $R\geq0$ and $\int R=1$, we actually have $\hat R_*=\hat{r}_0=1$.

\begin{lemma}\label{l.bdDe1}
For any $p\in[1,+\infty)$ and $T>0$, there exists $C>0$, depending only on
  $ p$ and $T$, such that  
\[
\begin{aligned}
&\sup_{t\in[0,T]}\big(\EE\la |\cD\mathscr{P}_t(v)|,h\ra_{L^2(\bT^d)}^p\big)^{1/p} \leq C\|h\|_{L^2(\bT^d)}\|v\|_{L^2(\bT^d)}(1+\|v\|_{L^2(\bT^d)})\\
&\leq  C\|h\|_{L^2(\bT^d)}\|v\|_{L^\infty(\bT^d)},\quad   \, v\in D^\infty(\bT^d), \, 0\leq h\in L^2(\bT^d).
\end{aligned}
\]
\end{lemma}

\begin{lemma}\label{l.bdDe2}
For any $p\in[1,+\infty)$ and $T>0$, there exists $C>0$, depending only on
  $  p$ and $T$, such that
\[
\begin{aligned}
\sup_{t\in[0,T]}\big(\EE  \la h_1,
|\cD^2\mathscr{P}_t(v)|h_2\ra_{L^2(\bT^d)}^p\big)^{1/p}\leq 
C\|h_1\|_{L^2(\bT^d)}\|h_2\|_{L^2(\bT^d)}(1+\|v\|_{L^\infty(\bT^d)} )
\end{aligned}
\]
for all $v\in D^\infty(\bT^d)$, $0\leq h_1,h_2\in L^2(\bT^d)$.
\end{lemma}

\begin{proof}[Proof of Lemma~\ref{l.bdDe1}]
From \tk{\eqref{e.4121} and \eqref{e.conditionR}}
we have
\[
\begin{aligned}
\la |\cD\mathscr{P}_t(v)|,h\ra_{L^2(\bT^d)} \leq &C\|U(t;h,v)\|_{L^2(\bT^d)}\|u(t;v)\|_{L^2(\bT^d)}(1+\|u(t;v)\|_{L^2(\bT^d)})\\
\leq & C\|U(t;h,v)\|_{L^2(\bT^d)}\|u(t;v)\|_{L^\infty(\bT^d)}.
\end{aligned}
\]
Here we have used the fact that $\|f\|_{L^2(\bT^d)}^2 \leq
\|f\|_{L^\infty(\bT^d)}$ for any $f\in D(\bT^d)$. The proof is
complete by  applying Lemma~\ref{l.lshe1} and   the H\"older inequality.
\end{proof}

\begin{proof}[Proof of Lemma~\ref{l.bdDe2}]
From \eqref{e.frede2} \tk{and \eqref{e.conditionR}}, we have
\begin{equation}\label{e.4141}
\begin{aligned}
 &\la h_1, |\cD^2\mathscr{P}_t(v)|h_2\ra_{L^2(\bT^d)}\\
 &\leq C\|U(t;h_1,v)\|_{L^2(\bT^d)}\|U(t;h_2,v)\|_{L^2(\bT^d)}(1+\|u(t;v)\|_{L^\infty(\bT^d)} ).
 \end{aligned}
 \end{equation}
 Similarly, the proof is complete by applying Lemma~\ref{l.lshe1} and
 the H\"older inequality.
\end{proof}

Now we define the gradient of the time dependent corrector: for any
$v\in D^\infty(\bT^d)$, $t,\,T\geq0$, we let
\begin{equation}\label{e.defDe11}
\cD\cP_t\cR(v):=\EE[\cD\mathscr{P}_t(v)],\quad\quad \cD^2\cP_t\cR(v):=\EE[\cD^2\mathscr{P}_t(v)],
\end{equation}
where $\cD\mathscr{P}_t$ and $\cD^2\mathscr{P}_t$ are given by
\eqref{e.fredePomega} and \eqref{e.frede2} respectively. Let
\begin{equation}\label{e.defgraCor}
\cD\chi(T,v):=\int_0^T\cD\cP_t\cR(v)dt.
\end{equation}
The following two lemmas concern the continuity of the Frechet
derivatives, and their proofs also only rely on the property of $\cR$
through \eqref{e.conditionR}.
\begin{lemma}\label{l.contDe1}
There exists $C=C(T)>0$ depending on $T>0$ such that
\begin{equation}
\label{032304-21z}
\begin{aligned}
&\sup_{t\in[0,T]}|\la \cD\cP_t\cR(v_1)-\cD\cP_t\cR(v_2),h\ra_{L^2(\bT^d)}|  \\
&\leq  C\|h\|_{L^2(\bT^d)}\|v_1-v_2\|_{L^2(\bT^d)}P(\|v_1\|_{L^2(\bT^d)},\|v_2\|_{L^2(\bT^d)}),
\end{aligned}
\end{equation}
for all  $v_1,v_2\in D^\infty(\bT^d)$ and $h\in L^2(\bT^d)$. 
Here $P(\cdot,\cdot)$ is some polynomial function.

Moreover, we have also
\begin{equation}
\label{032304-21a}
\begin{aligned}
&\sup_{t\in[0,T]}|\la h_1, \big(\cD^2\cP_t\cR(v_1)-\cD^2\cP_t\cR(v_2)\big)h_2\ra_{L^2(\bT^d)}| \\
&\leq  C \|h_1\|_{L^2(\bT^d)}\|h_2\|_{L^2(\bT^d)}\|v_1-v_2\|_{L^2(\bT^d)}P(\|v_1\|_{L^2(\bT^d)},\|v_2\|_{L^2(\bT^d)}),
\end{aligned}
\end{equation}
for all $v_1,v_2\in D^\infty(\bT^d)$ and $h_1,h_2\in
L^2(\bT^d)$. Here, again,  
 $P(\cdot,\cdot)$ is some polynomial.
\end{lemma}

\begin{proof}[Proof of Lemma~\ref{l.contDe1}]
We only prove \eqref{032304-21z}, as the argument for
\eqref{032304-21a} follows the same lines. From Lemma~\ref{l.lshe1} and the fact that $\U$ solves a linear equation, we know that for any $t\in[0,T]$, $h\in L^2(\bT^d), v_1,v_2\in D^\infty(\bT^d)$, we have

(i) $\big(\EE[\|u(t;v_1)-u(t;v_2)\|_{L^2(\bT^d)}^p]\big)^{1/p}\leq C\|v_1-v_2\|_{L^2(\bT^d)}(1+\|v_1\|_{L^2(\bT^d)})$,

(ii) $\big(\EE[ \|U(t;h,v_1)-U(t;h,v_2)\|_{L^2(\bT^d)}^p] \big)^{1/p}\leq C \|h\|_{L^2(\bT^d)}\|v_1-v_2\|_{L^2(\bT^d)}$.

%
To see   (i) note that
\begin{align*}
\|u(t;v_1)-u(t;v_2)\|_{L^2(\bT^d)}\le \frac{
  \|\U(t;v_1-v_2)\|_{L^2(\bT^d)}}{\|\U(t;v_2)\|_{L^1(\bT^d)}}+\frac{\|u(t;v_1)\|_{L^2(\bT^d)}
  \|\U(t;v_1-v_2)\|_{L^1(\bT^d)}}{\|\U(t;v_2)\|_{L^1(\bT^d)}},
\end{align*}
\tk{Here $\U(t;v_1-v_2):=\U(t;v_1)-\U(t;v_2)$.} The estimate follows then directly from Lemma~\ref{l.lshe1} and
estimate \eqref{e.negammZ} (recall that
$Z_t=\|\U(t;v)\|_{L^1(\bT^d)}$).
 Formula (ii) follows by similar considerations.

 \tk{Using \eqref{e.fredePomega} we conclude that
\[
\begin{aligned}
&\la \EE\cD\mathscr{P}_t(v_1),h\ra_{L^2(\bT^d)}-\la
\EE\cD\mathscr{P}_t(v_2),h\ra_{L^2(\bT^d)}=I_1+I_2,\quad\mbox{where}\\
&
I_1:=2\int_{\bT^{3d}}R(x-y)h(z)\EE\Big[U(t,x;z,v_1 )u(t,y;v_1
)-U(t,x;z,v_2 )u(t,y;v_2 ) \Big]dxdydz,\\
&
I_2:=2\int_{\bT^{3d}}R(x-y)h(z)\EE\Big[\|U(t;z,v_2)\|_{L^1(\bT^d)}u(t,x;v_2
)u(t,y;v_2 )\\
&
-\|U(t;z,v_1)\|_{L^1(\bT^d)}u(t,x;v_1 )u(t,y;v_1 )\Big] dxdydz.
\end{aligned}
\]
Note that
 \[
\begin{aligned}
&
|I_1|\le 2\Big|\int_{\bT^{3d}}R(x-y)h(z)\EE\Big[U(t,x;z,v_1 ) -U(t,x;z,v_2 ) \Big]u(t,y;v_1 )dxdydz\Big|\\
&
+2\Big|\int_{\bT^{3d}}R(x-y)h(z)\EE\Big[U(t,x;z,v_2 )\Big[u(t,y;v_1
)- u(t,y;v_2 ) \Big]dxdydz\Big|.
\end{aligned}
\]
Using \eqref{e.conditionR}, estimate \eqref{e.she13} together with the
already proved point   (ii) we conclude, upon an application of  the H\"older inequality, that
 \[
|I_1|\leq
C\|h\|_{L^2(\bT^d)}\|v_1-v_2\|_{L^2(\bT^d)} \big(\|v_1\|_{L^2(\bT^d)}+1\big)\|v_1\|_{L^2(\bT^d)}.
\]
A similar argument can be used to estimate $|I_2|$, therefore we infer that}
\[
\begin{aligned}
&|\la \EE\cD\mathscr{P}_t(v_1),h\ra_{L^2(\bT^d)}-\la \EE\cD\mathscr{P}_t(v_2),h\ra_{L^2(\bT^d)}|\\
&\leq C \|h\|_{L^2(\bT^d)}\|v_1-v_2\|_{L^2(\bT^d)}P(\|v_1\|_{L^2(\bT^d)},\|v_2\|_{L^2(\bT^d)})
\end{aligned}
\]
for some polynomial function $P(\cdot)$. The proof is complete.
\end{proof}



 \subsubsection{Estimate on $ \cD\cP_t\cR(v)$}
For the gradient of the corrector $\cD\chi(T,v)$ to be well-defined with $T$ approaching to infinity, one needs to refine the estimates on $\cD\cP_t\cR(v)$ for large $t$. Here is the main result  of this section.
  \begin{proposition}\label{p.gracor}
There exists $C,\lambda>0$ depending only on $R,d$ such that 
\begin{equation}
\label{011703-21z}
\begin{aligned}
\Big|\langle {\cal D}\cP_t\cR( v),h\rangle_{L^2(\bT^d)}\Big|\leq C\|h\|_{L^1(\bT^d)}\big(1_{[1,+\infty)}(t)e^{-\lambda t} +1_{[0,1]}(t)\|v\|_{L^\infty(\bT^d)}     \big)
\end{aligned}
\end{equation}
for all $t\geq 0$, $v\in D^\infty(\bT^d) $ and $ h\in L^2(\bT^d)$.
\end{proposition}

Directly from Proposition \ref{p.gracor} we conclude the following.
\begin{corollary}\label{cor013103-21}
For $C,\lambda>0$ as in Proposition \ref{p.gracor}  we have 
\begin{equation}
\label{011703-21}
\begin{aligned}
&\|{\cal D}\cP_t\cR( v)\|_{L^\infty(\bT^d)} \leq C\big(1_{t>1}e^{-\lambda t} +1_{t\in[0,1]}\|v\|_{L^\infty(\bT^d)}     \big),\\
&
\|\cD\chi(t,v)\|_{L^\infty(\bT^d)} \leq \tk{C\Big[\frac{1}{\lambda}
+\|v\|_{L^\infty(\bT^d)}\Big]}
\end{aligned}
\end{equation}
for all $t\geq 0$, $v\in D^\infty(\bT^d)$.
\end{corollary}

%

The proof of Proposition~\ref{p.gracor} relies on the following lemmas.
\begin{lemma}\label{l.471}
For any $p\in[1,+\infty)$ there exist $C,\lambda>0$ depending only on $R,d,p$ such that
\[
\EE\|u(t;\nu_1)-u(t;\nu_2)\|_{L^\infty(\bT^d)}^p \leq Ce^{-\lambda
  t},\quad\quad  t> 1,\,\nu_1,\nu_2\in \mathcal{M}_1(\bT^d).
\]
\end{lemma}
\begin{proof}
By
Proposition~\ref{p.couplingCon}, for any $\tfrac{t}{2}<k<N(t)-2$, we have
\[
\begin{aligned}
&\EE\|  u(t;\nu_1)-  u(t;\nu_2)\|_{L^\infty(\bT^d)}^p\\
&\leq C\big( \EE\|  u(t;\nu_1)-  u_k(t)\|_{L^\infty(\bT^d)}^p+\EE\|  u(t;\nu_2)-
u_k(t)\|_{L^\infty(\bT^d)}^p\big)\leq Ce^{-\lambda k}
\end{aligned}
\]
for some $C,\lambda>0$ only depending on $R,d,p$. The proof is complete.
\end{proof}

We have the following estimate on $U$:
\begin{lemma}\label{l.bdU}
If $p\in[1,+\infty)$, then there exists $C$ depending on $R,d,p$ such that 
\[
\EE\left[\|U(t;h,v)\|_{L^1(\bT^d)}^p\right]\leq
C\|h\|_{L^1(\bT^d)}^p,\quad \quad v\in D(\bT^d),\, 0\leq h\in L^1(\bT^d),\,t\ge0.
\]
\end{lemma}
\begin{proof}
By definition, we have
\[
\int_{\bT^d} U(t,x;h,v)dx=\frac{\int_{\bT^d} \U(t,x;h)dx}{\int_{\bT^d}\U(t,x;v)dx}. 
\]
Through a time reversal of the noise $\xi(\cdot,\cdot)\mapsto \xi(t-\cdot,\cdot)$, we have, for each $t\geq 0$,
\[
\frac{\int_{\bT^d} \U(t,x;h)dx}{\int_{\bT^d}\U(t,x;v)dx}\stackrel{\rm law}{=}\frac{\int_{\bT^d}\U(t,x;\1)h(x)dx}{\int_{\bT^d} \U(t,x;\1)v(x)dx}=\frac{\int_{\bT^d}u(t,x;\1)h(x)dx}{\int_{\bT^d} u(t,x;\1)v(x)dx},
\]
where  $\1(x)\equiv1$. Thus, it suffices to estimate the moments of the r.h.s. We have
\[
\begin{aligned}
\frac{\int_{\bT^d}u(t,x;\1)h(x)dx}{\int_{\bT^d} u(t,x;\1)v(x)dx} \leq \|h\|_{L^1(\bT^d)} \sup_{x\in\bT^d} u(t,x;\1) (\inf_{y\in\bT^d} u(t,x;\1))^{-1}.
\end{aligned}
\]
For $t>1$, we use \eqref{e.utUL} and apply Lemma~\ref{l.ulbd} to complete the proof. For $t\leq 1$, we use the estimate 
\begin{align*}
&
\sup_{x\in\bT^d} u(t,x;\1) (\inf_{x\in\bT^d} u(t,x;\1))^{-1} =
       \sup_{x\in\bT^d}\U(t,x;\1) \sup_{x\in\bT^d} \U(t,x;\1)^{-1}\\
&
\le \tfrac{1}{2} \sup_{x\in\bT^d}\U(t,x;\1)^2+ \tfrac{1}{2} \sup_{x\in\bT^d} \U(t,x;\1)^{-2}
 \end{align*}
 and Lemma~\ref{l.U1} to complete the proof.
\end{proof}
 
\begin{proof}[Proof of Proposition~\ref{p.gracor}]
Assuming first that $h\in L_+^2(\bT^d)$, recall  from  \eqref{e.4122}, we have 
\[
\begin{aligned}
\langle {\cal D}\cP_t\cR(v),h\rangle_{L^2(\bT^d)}=2\EE[ \cR(u(t;\tilde{h})-u(t;v),u(t;v))\|U(t;h,v)\|_{L^1(\bT^d)}].
\end{aligned}
\]
Here, as we recall, $\tilde{h}=\|h\|_{L^1(\bT^d)}^{-1}h$.
For $t>1 $, we use the fact that \[
|\langle {\cal D}\cP_t\cR(v),h\rangle_{L^2(\bT^d)}| \leq 2\EE[\|u(t;\tilde{h})-u(t;v)\|_{L^\infty(\bT^d)}\|u(t;v)\|_{L^1(\bT^d)}\|U(t;h,v)\|_{L^1(\bT^d)}]
\]
and apply Lemmas~\ref{l.471}, \ref{l.bdU}
and the H\"older inequality to derive that 
\[
|\langle {\cal D}\cP_t\cR(v),h\rangle_{L^2(\bT^d)}| \leq Ce^{-\lambda t}\|h\|_{L^1(\bT^d)}.
\]
For $t\leq 1$, 
we use a different expression: in \eqref{e.fredePomega}, we  use the fact that 
\[
\int_{\bT^{2d}}R(x-y)U(t,x;h,v)u(t,y;v)dxdy\leq\|R\|_{L^1(\bT^d)} \|u(t;v)\|_{L^\infty(\bT^d)}\|U(t;h,v)\|_{L^1(\bT^d)}
\] to derive that 
\[
|\la \cD\mathscr{P}_t(v), h\ra_{L^2(\bT^d)}| \leq 4 \|R\|_{L^1(\bT^d)}\|u(t;v)\|_{L^\infty(\bT^d)}\|U(t;h,v)\|_{L^1(\bT^d)},
\]
Taking the expectation, applying Lemma~\ref{l.bdU}, estimate
\eqref{e.she14}   and
the H\"older inequality, we obtain 
\[
\begin{aligned}
|\langle {\cal D}\cP_t\cR(v),h\rangle_{L^2(\bT^d)}| 
&\leq C \|v\|_{L^\infty(\bT^d)}\|h\|_{L^1(\bT^d)}.
\end{aligned}
\]
 
To generalize the proof to $h\in
L^2(\bT^d)$, it suffices to write $h=h_+-h_-$ and apply the above discussion to $h_+,h_-$ separately. The proof is complete.
\end{proof}

Define \begin{equation}\label{e.defDchi}
\cD\chi(v)=\int_0^\infty\cD\cP_t\cR(v)dt, \quad\quad v\in D^\infty(\bT^d).
\end{equation}
As a direct consequence of \eqref{e.4122} and Corollary~\ref{cor013103-21}, we have 
\begin{proposition}\label{p.bddeCor}
For any $v\in D^\infty(\bT^d)$, we have 
\begin{equation}\label{e.vDv}
\la v,\cD\chi(v)\ra_{L^2(\bT^d)}=0.
\end{equation}
In addition there exist $C,\lambda>0$ such that
\begin{equation}\label{e.4171}
\|\cD \chi(v)\|_{L^\infty(\bT^d)} \leq
C(1+\|v\|_{L^\infty(\bT^d)}),\quad v\in D^\infty(\bT^d),
\end{equation}
and
\begin{equation}
\label{010104-21}
\sup_{v\in D^\infty(\bT^d)}\| {\cal D}\chi(t,v)-{\cal D}\chi
(v)\|_{L^\infty(\bT^d)}\leq \frac{C}{\lambda}e^{-\lambda t},\quad t>1.
\end{equation}
\end{proposition}

\subsection{A semimartingale decomposition} The goal in this section is to obtain a semimartingale decomposition of the process $\{\cP_s\cR(u(t;v))\}_{t\geq0}$, 
for any fixed $s\geq0$ and $v\in D^\infty(\bT^d)$. Recall that, for any (fixed) $v\in D^\infty(\bT^d)$, we have
\[
\cP_s\cR(v)=\EE[\cR(u(s;v))]=\EE[\mathscr{P}_s(v)].
\]
From now on, assume that, in \eqref{e.defPtomega}, we use an
independent copy of noise, denoted by $\tilde{\xi}$, to generate the
random variable
$\mathscr{P}_s(v)$. To emphasize the dependence, we will write $\mathscr{P}_s(v)=\mathscr{P}_s(v;\tilde{\xi})$, so for any $v\in D^\infty(\bT^d)$, 
\[
\cP_s\cR(v)=\EE_{\tilde{\xi}}[\mathscr{P}_s(v;\tilde{\xi})],
\]
and 
\[
\cP_s\cR(u(t;v,\xi))=\EE_{\tilde{\xi}}[\mathscr{P}_s(u(t;v,\xi);\tilde{\xi})],
\]
where $\EE_{\tilde{\xi}}$ denotes the expectation only with respect to
$\tilde{\xi}$, and the $u(t;v,\xi)$ in the above expression is based on the original noise $\xi$.

By Lemma~\ref{l.fredeomega}, we know that for almost every realization
of $\tilde{\xi}$, $\mathscr{P}_s(\cdot;\tilde{\xi})$ is twice Frechet
differentiable at any $v\in D^\infty(\bT^d)$. The same proof actually
shows that $\mathscr{P}_s(\cdot;\tilde{\xi})$ is infinitely Frechet
differentiable. 
The idea is to fix the realization of  $\tilde{\xi}$, apply
the It\^o formula to $\mathscr{P}_s(u(t;v);\tilde{\xi})$, then take the expectation with respect to $\tilde{\xi}$ to obtain a semimartingale decomposition of $\cP_s\cR(u(t;v))$.

\subsubsection{The case of smooth noise and initial data}

In this section, we assume that $v\in D(\bT^d)\cap C^\infty(\bT^d)$ (so it is automatically in $D^\infty(\bT^d)$) and the noise $\xi$  satisfies Assumption~\ref{a.smooth}, i.e., it contains only finitely many Fourier modes. Since the initial data is fixed, we omit its dependence to write $u(t)=u(t;v)$. By Proposition~\ref{p.spde}, we know that $u(t)$ is a strong solution to the following SPDE
\begin{equation}\label{e.spde1} 
\begin{aligned}
du(t)=\mathscr{A}(u(t))dt+\sum_{k\in\Z^d} \mathscr{B}_k(u(t))dw_k(t),
\end{aligned}
\end{equation}
with the operators $\mathscr{A},\mathscr{B}_k$ defined in \eqref{e.defAB}. 

Fix $s\geq0$ and a realization of $\tilde{\xi}$, applying It\^o
formula, see e.g. \cite[Theorem 4.32, p. 106]{daza}, 
we have 
\begin{equation}\label{e.itoPs}
\begin{aligned}
\mathscr{P}_s(u(t);\tilde{\xi})-\mathscr{P}_s(v;\tilde{\xi})&=\int_0^t \la \cD\mathscr{P}_s(u(r);\tilde{\xi}),u(r)dW(r)\ra_{L^2(\bT^d)} \\
&+\int_0^t \la \cD\mathscr{P}_s(u(r);\tilde{\xi}), \mathscr{A}(u(r))\ra_{L^2(\bT^d)} dr\\
&+\tfrac12\sum_{k\in\Z^d}\int_0^t \la \mathscr{B}_k(u(r)), \cD^2\mathscr{P}_s(u(r);\tilde{\xi})\mathscr{B}_{k}(u(r))\ra_{L^2(\bT^d)}dr.
\end{aligned}
\end{equation}
Here we have used the fact that
\[
\la \cD\mathscr{P}_s(v;\tilde{\xi}),v\ra_{L^2(\bT^d)}=0, \quad\quad \mbox{ if } v\in D^\infty(\bT^d),
\]
which comes from \eqref{e.4122}.

Freeze $(u(r;\xi))_{r\in[0,t]}$ and take the expectation with respect to $\tilde{\xi}$ in \eqref{e.itoPs}. By  Lemmas~\ref{l.bdDe1} and \ref{l.bdDe2}, we interchange the expectation and the   integral to obtain 
\begin{equation}\label{e.itoPs1}
\begin{aligned}
\cP_s\cR(u(t) )-\cP_s\cR(v)=&\int_0^t \la
\cD\cP_s\cR(u(r)),u(r)dW(r)\ra_{L^2(\bT^d)} \\
&+\int_0^t
\cL\cP_s\cR(u(r))dr,
\end{aligned}
\end{equation}
with 
\begin{equation}\label{e.defL1}
\begin{aligned}
\cL\cP_s\cR(v):=&\la \cD\cP_s\cR(v),\mathscr{A}(v)\ra_{L^2(\bT^d)}\\
&+\frac12\sum_{k\in\Z^d} \la \mathscr{B}_k(v), \cD^2\cP_s\cR(v )\mathscr{B}_{k}(v)\ra_{L^2(\bT^d)}.
\end{aligned}
\end{equation}
for any $s\geq0$ and $v\in  D(\bT^d)\cap C^\infty(\bT^d)$. In addition
$\cD\cP_s\cR$ and $\cD^2\cP_s\cR$ are defined in \eqref{e.defDe11}.

The following lemma holds:
\begin{lemma}\label{l.generator1}
For any $T>0$, $v\in D(\bT^d)\cap C^\infty(\bT^d)$ and $p\in[1,+\infty)$, we have 
\begin{equation}\label{e.bddrift}
\sup_{s,r\in[0,T]}\EE|\cL\cP_s\cR( u(r))|^p \leq C(T,p,v).
\end{equation}
For any $s\geq0$ and $v$ as above,  we have 
\begin{equation}\label{e.generatorcon}
\begin{aligned}
\cL\cP_s\cR(v )=\lim_{\delta\to0}\delta^{-1}[\cP_{s+\delta}\cR(v)-\cP_s\cR(v)].
\end{aligned}
\end{equation}
  \end{lemma}
  \begin{proof}
  First, by Lemmas~\ref{l.bdDe1} and \ref{l.bdDe2}, we have 
  \[
  \begin{aligned}
  |\cL\cP_s\cR( u(r))| &\leq C  \|\mathscr{A}(u(r))\|_{L^2(\bT^d)}\|u(r)\|_{L^\infty(\bT^d)}\\
  &+C \sum_{k\in\Z^d}\|\mathscr{B}_k(u(r))\|_{L^2(\bT^d)}^2(1+\|u(r)\|_{L^\infty(\bT^d)}).
  \end{aligned}
  \]
  Then we apply Lemmas~\ref{l.bdspde} and \ref{l.bdRv} to conclude
  \eqref{e.bddrift}. To show \eqref{e.generatorcon}, fix
  $s\geq0$. Define a continuous local martingale
  \[
  \M_s(t):=\int_0^t \la \cD\cP_s\cR(u(r)),u(r)dW(r)\ra_{L^2(\bT^d)},
  \] 
   By the equality \eqref{e.defDe11} and  Lemma~\ref{l.bdDe1}, its quadratic variation satisfies 
  \[
  \la \M_s\ra_t=\int_0^t\cR\Big(\cD\cP_s\cR(u(r))u(r)\Big)dr\leq C\int_0^t \|u(r)\|_{L^\infty(\bT^d)}^4dr.
  \]
 Next we apply Lemma~\ref{l.bdRv} to conclude that $\EE[\la
 \M_s\ra_t]\leq Ct$, so $\M_s$ is a bona fide, square integrable,
 continuous trajectory martingale. Therefore from \eqref{e.itoPs1} we have 
  \[
  \begin{aligned}
  \delta^{-1}[\cP_{s+\delta}\cR(v)-\cP_s\cR(v)]=&\delta^{-1}\big(\EE[\cP_s\cR(u(\delta))]-\cP_s\cR(v)\big)\\
  =&\delta^{-1}\int_0^\delta\EE[\cL\cP_s\cR(u(r))]dr.
  \end{aligned}
  \]
By \eqref{e.bddrift}, Lemma~\ref{l.contDe1}, and the fact that the sample path of $u(\cdot)$ lies in $C([0,\infty),H^k(\bT^d))$ for any $k\geq1$, we conclude the proof.
    \end{proof}

 \begin{remark}
 Formula \eqref{e.itoPs1} 
yields  the semimartingale decomposition of the process
$(\cP_s\cR(u(t)))_{t\ge0}$ for any $s\ge0$.
\end{remark}

The drift term in \eqref{e.itoPs1} can be simplified after a time integration:
\begin{lemma}\label{l.interchangeIn}
For any $T,t\geq0$, 
\[
\int_0^T \big(\int_0^t  \cL\cP_s\cR(u(r))  dr\big) ds=\int_0^t \cP_T\cR(u(r))dr-\int_0^t\cR(u(r))dr.
\]
\end{lemma}

\begin{proof}
First, by \eqref{e.bddrift}, we have
\[
\int_0^T \big(\int_0^t  \cL\cP_s\cR(u(r))  dr\big) ds=\int_0^t \big(\int_0^T  \cL\cP_s\cR(u(r))  ds\big)dr.
\]
It remains to use \eqref{e.generatorcon} to derive that for any $v\in D(\bT^d)\cap  C^\infty(\bT^d)$,
\[
{\cal P}_T \cR(v)-\cR(v)=\int_0^T \cL \cP_s\cR(v)ds.
\]
The proof is complete.
\end{proof}

Using the semimartingale decomposition in \eqref{e.itoPs1} and Lemma~\ref{l.interchangeIn}, we have
\begin{proposition}\label{p.madefinal}
Suppose that Assumption~\ref{a.smooth} is in force and further assume
that the initial condition $v\in D(\bT^d)\cap C^\infty(\bT^d)$. Then,
for any $T,t\geq0$ we have
\begin{equation}\label{e.matdecom}
\begin{aligned}
\int_0^t \cR(u(r;v))dr=&\int_0^t \la \cD\chi(T,u(r;v)),u(r;v)dW(r)\ra_{L^2(\bT^d)}\\
&-\chi(T,u(t;v))+\chi(T,v)+\int_0^t \cP_T \cR(u(r;v))dr.
\end{aligned}
\end{equation}
\end{proposition}

\begin{proof}
In \eqref{e.itoPs1}, we integrate $s$ from $0$ to $T$ and apply Lemma~\ref{l.interchangeIn} to derive
\[
\begin{aligned}
\chi(T,u(t))-\chi(T,v)&=\int_0^T\big(\int_0^t \la \cD\cP_s\cR(u(r)),u(r)dW(r)\ra_{L^2(\bT^d)}\big) ds\\
&+\int_0^t \cP_T \cR(u(r))dr-\int_0^t \cR(u(r))dr.
\end{aligned}
\] 
For the first term in the right hand side, we use the stochastic
Fubini theorem, see e.g. \cite[Theorem 4.33, p. 110]{daza},  and obtain
\[
\begin{aligned}
&\int_0^T\big(\int_0^t \la \cD\cP_s\cR(u(r)),u(r)dW(r)\ra_{L^2(\bT^d)}\big) ds\\
&=\int_0^t \la \int_0^T  \cD\cP_s\cR(u(r))ds, u(r)dW(r)\ra_{L^2(\bT^d)}=\int_0^t \la \cD\chi(T,u(r)),u(r)dW(r)\ra_{L^2(\bT^d)},
\end{aligned}
\]
where $\cD\chi$ was defined in \eqref{e.defgraCor}. The proof is complete.
\end{proof}


\subsubsection{The general case}
The following is the main result of this section, which removes the smoothness assumptions on the noise and the initial data in Proposition~\ref{p.madefinal}.
\begin{proposition}\label{p.41711}
For any $T,t\geq0$ and $v\in D^\infty(\bT^d)$, the decomposition \eqref{e.matdecom} holds.
\end{proposition}

A direct consequence of the above proposition is the decomposition of
the additive functional appearing in the formula for $\log Z_t$, see
\eqref{e.delogZt}, into a stochastic integral with respect to $dW(t)$
and   boundary terms.
\begin{corollary}
For any $t\geq0$ and $v\in D^\infty(\bT^d)$, we have 
\begin{equation}\label{e.matdecom2}
\begin{aligned}
\int_0^t \tilde{\cR}(u(r;v))dr=\int_0^t &\la \cD\chi(u(r;v)),u(r;v)dW(r)\ra_{L^2(\bT^d)}\\
&-\chi(u(t;v))+\chi(v).
\end{aligned}
\end{equation}
\end{corollary}
\begin{proof}
First, in both sides of \eqref{e.matdecom}, we subtract $2\gamma t$ to obtain 
\begin{equation}
\begin{aligned}
\int_0^t \tilde{\cR}(u(r;v))dr=&\int_0^t \la \cD\chi(T,u(r;v)),u(r;v)dW(r)\ra_{L^2(\bT^d)}\\
&-\chi(T,u(t;v))+\chi(T,v)+\int_0^t \cP_T\tilde{\cR}(u(r;v))dr.
\end{aligned}
\end{equation}
Sending $T\to\infty$ in the above equation, and applying
Lemma~\ref{l.bdRv} together with estimates \eqref{cellt1} and \eqref{010104-21}, we conclude the proof.
\end{proof}

The rest of the section is devoted to proving Proposition~\ref{p.41711}. Fix $v\in D^\infty(\bT^d)$ from now on. 
Let us introduce some notations. Recall that 
$(\hat r_k)_{k\in\bbZ^d}$ are the Fourier coefficients of $R(\cdot)$, see \eqref{e.defrk}.
For any $N>0$, define 
\begin{equation}
R_N(x)=\sum_{|k|\leq N} \hat{r}_k e^{i2\pi k\cdot x}.
\end{equation}
Let $\xi_N$ be a finite dimensional noise that is white in time and smooth in space with the spatial covariance function  $R_N$, and 
\[
dW_N(t,x)=\xi_N(t,x)dt=\sum_{|k|\leq N} \sqrt{\hat{r}_k} e^{i2\pi k\cdot x}d{w}_k(t).
\] 
Given $v\in D^\infty(\bT^d)$, we choose and fix a sequence of $v_N\in D(\bT^d)\cap C^\infty(\bT^d)$ such that 
\begin{equation}
\label{vN}
\|v_N-v\|_{L^\infty(\bT^d)}\to0, \quad\quad N\to\infty.
\end{equation}

Let $\U_N(t;v)$ be the solution to SHE with noise $\xi_N$ and initial
data $v$. Let furthermore 
\[
u_N(t,x;v)=\frac{\U_N(t,x;v)}{\int_{\bT^d} \U_N(t,x';v)dx'},
\]
\[
U_N(t,x;h,v)=\frac{\U_N(t,x;h)}{\int_{\bT^d} \U_N(t,x';v)dx'}.
\]
Define 
\[
\begin{aligned}
&\cR_N(v_1,v_2)=\int_{\bT^{2d}}R_N(x-y)v_1(x)v_2(y)dxdy, \quad\quad v_1,v_2\in L^2(\bT^d),\\
&\cR_N(v)=\cR_N(v,v), \quad\quad v\in L^2(\bT^d),
\end{aligned}
\]
and denote $\cP_t^N$ the semigroup associated with $u_N$, so 
\[
\cP_t^N\cR_N(v)=\EE[\cR_N(u_N(t;v))].
\]
Define
\[
\begin{aligned}
&\chi_N(T,v)=\int_0^T \cP_t^N\cR_N(v)dt-2\gamma T,\\
&\mathscr{P}_t^N(v):=\int_{\bT^{2d}} R_N(x-y) \frac{\U_N(t,x;v)\U_N(t,y;v)}{\left(\int_{\bT^d}\U_N(t,x';v)dx'\right)^2}dxdy,
\end{aligned}
\]
and
\[
\cD\chi_N(T,v)=\int_0^T \EE[\cD\mathscr{P}_t^N(v)]dt.
\]
Then by Proposition~\ref{p.madefinal}, we have 
\begin{equation}\label{e.4151}
\begin{aligned}
\int_0^t \cR_N(u_N&(r;v_N))dr=\int_0^t \la \cD\chi_N(T,u_N(r;v_N)),u_N(r;v_N)dW_N(r)\ra_{L^2(\bT^d)}\\
&-\chi_N(T,u_N(t;v_N))+\chi_N(T,v_N)+\int_0^t \cP_T^N \cR_N(u_N(r;v_N))dr.
\end{aligned}
\end{equation}

The goal is to pass to the limit  in \eqref{e.4151} as
$N\to\infty$. We first state the following elementary lemmas.

\begin{lemma}\label{l.conuN}
For any $T\geq0$, $p\in[1,+\infty)$, $h\in L^2(\bT^d)$, 
  the following convergences take place
\[
\begin{aligned}
&\sup_{t\in[0,T]}\EE\|u_N(t;v)-u(t;v)\|_{L^2(\bT^d)}^p\to0,\\
&\sup_{t\in[0,T]}\EE\|u_N(t;v_N)-u(t;v)\|_{L^2(\bT^d)}^p\to0,\\
&\sup_{t\in[0,T]}\EE\|U_N(t;h,v)-U(t;h,v)\|_{L^2(\bT^d)}^p\to0,\quad
\mbox{as $N\to\infty$}.
\end{aligned}
\]
Furthermore, for any $t>0$, as $N\to\infty$, 
\begin{equation}
\label{012604-21}
\EE \big(\int_{\bT^{2d}} [\U_N(t,x;y)-\U(t,x;y)]^2 dxdy\big)^p\to0.
\end{equation}
\end{lemma}

The proof is  standard and we present it in
Section~\ref{s.sheOther} of the appendix.


%

\begin{lemma}\label{l.conRN}
(i) For any $v_1,v_2\in L^2(\bT^d)$, we have 
\[
\begin{aligned}
&|\cR_N(v_1,v_2)-\cR(v_1,v_2)|\leq  \|v_1\|_{L^2(\bT^d)}\|v_2\|_{L^2(\bT^d)},\\
&\cR_N(v_1,v_2)\to \cR(v_1,v_2), \quad\quad N\to\infty.
\end{aligned}
\]
 (ii) For any $v_j\in L^2(\bT^d),j=1,\ldots,4$, we have 
\[
\tk{|\cR(v_1,v_2)-\cR(v_3,v_4)|} \leq  \|v_1-v_3\|_{L^2(\bT^d)}\|v_2\|_{L^2(\bT^d)}+ \|v_2-v_4\|_{L^2(\bT^d)}\|v_3\|_{L^2(\bT^d)}.
\]
\end{lemma}

\begin{proof}
It suffices to write the integrals in the Fourier domain 
\[
\cR_N(v_1,v_2)=\sum_{|k|\leq N}\hat{r}_k\hat{v}_1(k)\hat{v}_2^*(k), \quad\quad \cR(v_1,v_2)=\sum_{k\in \Z^d}\hat{r}_k\hat{v}_1(k)\hat{v}_2^*(k),
\]
and use the fact that $0\leq \hat{r}_k\leq \hat R_*=1$ to conclude the
proof of (i). Similarly, for (ii), by virtue of the Cauchy-Schwarz
inequality, we have
\[
\begin{aligned}
&|\cR_N(v_1,v_2)-\cR_N(v_3,v_4)|\leq  \sum_{k\in\Z^d}   \big|\hat{v}_1(k)\hat{v}_2^*(k)-\hat{v}_3(k)\hat{v}_4^*(k)\big|\\
&\leq \|v_1-v_3\|_{L^2(\bT^d)}\|v_2\|_{L^2(\bT^d)}+ \|v_2-v_4\|_{L^2(\bT^d)}\|v_3\|_{L^2(\bT^d)}.
\end{aligned}
\]
\tk{The proof of part (ii) is complete, upon passing to the limit as $N\to\infty$.}
\end{proof}

%

With the above two lemmas, we have
\begin{proposition}\label{p.4161}
For any $t,T\geq0$, 
  we have
\begin{equation}
\begin{aligned}
&\chi_N(T,v_N)\to \chi(T,v),\\
&\chi_N(T,u_N(t;v_N))\to \chi(T,u(t;v)), \\
&\int_0^t \cR_N(u_N(r;v_N))dr\to \int_0^t \cR(u(r;v))dr,\\
&\int_0^t \cP_T^N \cR_N(u_N(r;v_N))dr\to  \int_0^t \cP_T\cR(u(r;v))dr,\quad\mbox{as $N\to\infty$}. 
\end{aligned}
\end{equation}
The convergence holds in $L^p(\Omega)$ for any $p\in[1,+\infty)$.
\end{proposition}

\begin{proof}
The proofs for all cases are similar, so we take the last one as an example. One can write
\[
\int_0^t \cP_T^N \cR_N(u_N(r;v_N))dr=\int_0^t \EE[\cR_N(u_N(r+T;v_N))|\F_r]dr,
\]
and
\[
\int_0^t \cP_T\cR(u(r;v))dr=\int_0^t \EE[\cR(u(r+T;v))|\F_r]dr.
\]
So the difference can be decomposed into two terms:
\[
\begin{aligned}
&I_1=\int_0^t \EE[\cR_N(u_N(r+T;v_N))|\F_r]dr-\int_0^t \EE[\cR_N(u(r+T;v))|\F_r]dr,\\
&I_2=\int_0^t \EE[\cR_N(u(r+T;v))|\F_r]dr-\int_0^t \EE[\cR(u(r+T;v))|\F_r]dr.
\end{aligned}
\]
We apply Lemmas~\ref{l.conuN} and \ref{l.conRN} to conclude that as $N\to\infty$, $I_1,I_2\to0$ in $L^p(\Omega)$ for any $p\in[1,+\infty)$.
\end{proof}


To deal with the martingale term in \eqref{e.4151}, we need the following lemmas:

\begin{lemma}\label{l.appFDe}
For any $f\in D^\infty(\bT^d)$ and $t>0$, we have as $N\to\infty$, 
\[
\cD\cP_t^N\cR_N(f)\to \cD\cP_t\cR(f) \quad  \mbox{ in } \quad  L^2(\bT^d).
\]
\end{lemma}
\begin{proof}
To simplify the notation we shall   assume that $\hat R_*=1$. Fix $t>0$. Recall from \eqref{e.fredePomega} that 
\[
\begin{aligned}
 \cD\mathscr{P}_t^N(f)(z)=&2\int_{\bT^{2d}}R_N(x-y)U_N(t,x;z,f)u_N(t,y;f)dxdy\\
&-2\|U_N(t;z,f)\|_{L^1(\bT^d)}\int_{\bT^{2d}}R_N(x-y)u_N(t,x;f)u_N(t,y;f)dxdy
\end{aligned}
\]
and $\cD\cP_t^N\cR_N(f)=\EE\cD\mathscr{P}_t^N(f)$. 

To prove the result, it suffices to show as $N\to\infty$,
\begin{equation}
\EE\int_{\bT^d}|\cD\mathscr{P}_t^N(f)(z)-\cD\mathscr{P}_t(f)(z)|^2dz\to0.
\end{equation}
One can write the error as $\cD\mathscr{P}_t^N(f)(z)-\cD\mathscr{P}_t(f)(z)=2(I_1^{(N)}(z)+I_2^{(N)}(z))$ with
\[
I_1^{(N)}(z)=\int_{\bT^{2d}}R_N(x-y)U_N(t,x;z,f)u_N(t,y;f)dxdy-\int_{\bT^{2d}}R(x-y)U(t,x;z,f)u(t,y;f)dxdy
\]
and
\[
\begin{aligned}
I_2^{(N)}(z)=&\|U(t;z,f)\|_{L^1(\bT^d)}\int_{\bT^{2d}}R(x-y)u(t,x;f)u(t,y;f)dxdy\\
&-\|U_N(t;z,f)\|_{L^1(\bT^d)}\int_{\bT^{2d}}R_N(x-y)u_N(t,x;f)u_N(t,y;f)dxdy.
\end{aligned}
\]

By Lemma~\ref{l.conRN}, we have 
\[
\begin{aligned}
|I_1^{(N)}(z)|\leq&\|U_N(t;z,f)-U(t;z,f)\|_{L^2(\bT^d)}\|u_N(t;f)\|_{L^2(\bT^d)}\\
&+\|U(t;z,f)\|_{L^2(\bT^d)}\|u_N(t;f)-u(t;f)\|_{L^2(\bT^d)}\\
&+|\int_{\bT^{2d}} (R_N(x-y)-R(x-y))U(t,x;z,f)u(t,y;f)dxdy|\\
&=:J_1^{(N)}(z)+J_2^{(N)}(z)+J_3^{(N)}(z)
\end{aligned}
\]
To deal with $J_1$, we first
note that
\begin{align*}
&
|U_N(t,x;z,f)-U(t,x;z,f)|\le
  |\U_N(t,x;z)-\U(t,x;z)|\cdot\|\U_N(t;f)\|_{L^1(\bT^d)}^{-1}\\
&
+U(t,x;z,f) \|\U_N(t;f)\|_{L^1(\bT^d)}^{-1}\|\U_N(t;f)-\U(t;f)\|_{L^1(\bT^d)}.
\end{align*}
Using Lemma~\ref{l.conuN}, \eqref{e.estiGreen1} and Lemma~\ref{l.MMZ},
together with the H\"older inequality, we
conclude that for any $p\in[1,+\infty)$
\begin{align*}
\lim_{N\to+\infty}\EE\left(\int_{\bT^{2d}}|U_N(t,x;z,f)-U(t,x;z,f)|^2dxdz\right)^p=0.
\end{align*}
Invoking Lemma~\ref{l.conuN} and  the H\"older inequality we conclude that
\[
\EE\|J_1^{(N)}\|_{L^2(\bT^d)}^2=\EE\Big[\|u_N(t;f)\|_{L^2(\bT^d)}^2\int_{\bT^{2d}}|U_N(t,x;z,f)-U(t,x;z,f)|^2dxdz\Big]
\to0.
\]
Similarly, we have
\[
\EE\|J_2^{(N)}\|_{L^2(\bT^d)}^2=\EE\Big[\|u_N(t;f)-u(t;f)\|_{L^2(\bT^d)}^2\int_{\bT^{2d}} U(t,x;z,f)^2dxdz\Big]\to0.
\]
Concerning $J_3^{(N)}$, denote by $\hat{U}_k(t;z,f), \hat{u}_k(t;f)$  the
Fourier coefficients of $U(t;z,f)$ and $u(t;f)$, respectively. We have
\[
\EE\|J_3^{(N)}\|_{L^2(\bT^d)}^2 \leq \EE \left[\sum_{|k|>N} \int_{\bT^d}|\hat{U}_k(t;z,f)|^2dz \sum_{|k|>N} |\hat{u}_k(t;f)|^2\right]\to0,
\]
by virtue of the monotone convergence theorem.
This shows that $\EE\|I_1^{(N)}\|_{L^2(\bT^d)}^2 \to0$. The proof that
$\EE\|I_2^{(N)}\|_{L^2(\bT^d)}^2 \to0$ can be carried out along similar lines.
\end{proof}
\begin{lemma}\label{l.42011}
For any $T,t\geq0$, $p\in[1,+\infty)$, 
 we have 
\begin{equation}\label{e.41711}
\sup_{r\in[0,t],N>1}\EE\big( \|\cD\chi_N(T,u_N(r;v_N))\|_{L^2(\bT^d)}^p+\|\cD\chi(T,u(r;v))\|_{L^2(\bT^d)}^p\big) \leq C.
\end{equation}
In addition, for any $r\in[0,t]$, 
\begin{equation}\label{e.41712}
 \EE\|\cD\chi_N(T,u_N(r;v_N))-\cD\chi(T,u(r;v))\|_{L^2(\bT^d)}^p\to0, \mbox{ as } N\to\infty.
\end{equation}
\end{lemma}
\begin{proof}
To simplify the notations, we let $f_N=u_N(r;v_N)$ and $f=u(r;v)$. 
Concerning \eqref{e.41711}, note that by Lemma~\ref{l.bdDe1} for a
fixed $T>0$ we can find $C>0$ such that
\[
\begin{aligned}
\|\cD\chi_N(T,f_N)\|_{L^2(\bT^d)}+\|\cD\chi(T,f)\|_{L^2(\bT^d)} \leq C\Big(\|f_N\|_{L^\infty(\bT^d)} +\|f\|_{L^\infty(\bT^d)}\Big).
\end{aligned}
\]
Applying Lemma~\ref{l.bdRv}, we conclude estimate \eqref{e.41711}. 
For \eqref{e.41712}, we write 
\[
\begin{aligned}
&\cD\chi_N(T,f_N)-\cD\chi(T,f)\\
&=\int_0^T\big[ \cD\cP_t^N\cR_N(f_N)-\cD\cP_t\cR(f)\big]dt\\
&=\int_0^T \big[\cD\cP_t^N\cR_N(f_N)-\cD\cP_t^N\cR_N(f)\big]dt+\int_0^T \big[\cD\cP_t^N\cR_N(f)-\cD\cP_t\cR(f)\big]dt.
\end{aligned}
\]

For the first term, by Lemma~\ref{l.contDe1}, we have 
\[
\begin{aligned}
\|\cD\cP_t^N\cR_N(f_N)-\cD\cP_t^N\cR_N(f)\|_{L^2(\bT^d)} \leq C\|f_N-f\|_{L^2(\bT^d)}P(\|f_N\|_{L^2(\bT^d)},\|f\|_{L^2(\bT^d)})
\end{aligned}
\]
where $P$ is some given polynomial function. Then we apply Lemmas~\ref{l.bdRv} and \ref{l.conuN} to derive that 
\[
\EE\int_0^T \|\cD\cP_t^N\cR_N(f_N)-\cD\cP_t^N\cR_N(f)\|_{L^2(\bT^d)}^pdt\to0,  \mbox{ as } N\to\infty.
\]

For the second term, we first note that, by Lemma~\ref{l.bdDe1}, we have
\[
\|\cD\cP_t^N\cR_N(f)-\cD\cP_t\cR(f)\|_{L^2(\bT^d)}\leq C\|f\|_{L^\infty(\bT^d)}.
\]
Then we apply Lemmas~\ref{l.bdRv} and \ref{l.appFDe}  to derive that 
\[
\EE\int_0^T \|\cD\cP_t^N\cR_N(f)-\cD\cP_t\cR(f)\|_{L^2(\bT^d)}^pdt\to0, \mbox{ as }N\to\infty.
\]
The proof is complete.
\end{proof}

Finally, we can prove the convergence of the martingale term in \eqref{e.4151}:
\begin{proposition}\label{p.41712}
For any $t,T\geq0$, 
   we have 
\[
\begin{aligned}
\int_0^t& \la \cD\chi_N(T,u_N(r;v_N)),u_N(r;v_N)dW_N(r)\ra_{L^2(\bT^d)}\\
&\to \int_0^t \la \cD\chi(T,u(r;v)),u(r;v)dW(r)\ra_{L^2(\bT^d)},
\mbox{ as $N\to\infty$, in $L^2(\Omega)$.}
\end{aligned}
\]
\end{proposition}
\begin{proof}
To simplify the notation, we denote
\[
\begin{aligned}
&f_N(r,x)=\cD\chi_N(T,u_N(r;v_N))(x)u_N(r,x;v_N),\\
&f(r,x)=\cD\chi(T,u(r;v))(x)u(r,x;v),
\end{aligned}
\] 
then the goal is to show
$\int_0^t\int_{\bT^d}f_N(r,x)dW_N(r,x)\to\int_0^t \int_{\bT^d}
f(r,x)dW(r,x)$ in $L^2(\Omega)$. The error can be written as a sum of  two terms:
\[
\begin{aligned}
I_{1,N}=&\int_0^t\int_{\bT^d} [ f_N(r,x)-f(r,x)]dW_N(t,x),\\
I_{2,N}=&\int_0^t \int_{\bT^d} f(r,x)dW_N(t,x)-\int_0^t \int_{\bT^d}f(r,x)dW(t,x).
\end{aligned}
\]

For the first term, we have 
\[
\begin{aligned}
\EE[I_{1,N}^2]=&\EE\int_0^t\int_{\bT^d}[f_N(r,x)-f(r,x)][f_N(r,y)-f(r,y)]R_N(x-y)dxdy \\
&\leq \int_0^t\EE\|f_N(r)-f(r)\|_{L^2(\bT^d)}^2dr.
\end{aligned}
\]
We bound the $L^2$ norm by
\[
\begin{aligned}
\|f_N(r)-f(r)\|_{L^2(\bT^d)} \leq& \|\cD\chi_N(T,u_N(r;v_N))-\cD\chi(T,u(r;v))\|_{L^2(\bT^d)}\|u_N(r;v_N)\|_{L^\infty(\bT^d)}\\
+&\|\cD\chi(T,u(r;v))\|_{L^\infty(\bT^d)}\|u_N(r;v_N)-u(r;v)\|_{L^2(\bT^d)}
\end{aligned}
\]
Then it suffices to apply Lemmas~\ref{l.bdRv}, \ref{l.conuN}, \ref{l.42011} and \eqref{011703-21} to conclude that $\EE[I_{1,N}^2]\to0$.

For the second term, with $\hat{f}_k(r)=\int_{\bT^d} f(r,x)e^{-i2\pi k\cdot x}dx$, we have
\[
\EE[I_{2,N}^2]=\sum_{|k|>N}\hat{r}_k\int_0^t  \EE |\hat{f}_k(r)|^2 dr \leq \sum_{|k|>N}\int_0^t  \EE |\hat{f}_k(r)|^2 dr,
\]
 which goes to zero as $N\to\infty$, by invoking the fact that 
 \[
 \begin{aligned}
 \EE\|f(r)\|_{L^2(\bT^d)}^p \leq &\EE \Big[\|\cD\chi(T,u(r;v))\|_{L^2(\bT^d)}^p\|u(r;v)\|_{L^\infty(\bT^d)}^p\Big]\\
 \leq & C\EE\Big[ (1+\|u(r;v)\|_{L^\infty(\bT^d)})^p\|u(r;v)\|_{L^\infty(\bT^d)}^p \Big]\leq C,
 \end{aligned}
 \]
 where we have used \eqref{011703-21} in the second inequality and
 \eqref{e.bduinfinity} in the third one. Here $C$ is a
 generic constant \tk{that may depend on $T$ and $v$, but is} independent of $N$.
\end{proof}

Now we can combine Propositions~\ref{p.4161} and \ref{p.41712} to conclude Proposition~\ref{p.41711}.

\subsection{Martingale CLT}
Define
\begin{equation}\label{e.4181}
{\frak z}_t:=\log Z_t+\gamma t=\int_0^t\int_{\bT^d} u(s,y;\nu)\xi(s,y)dyds-\tfrac12\int_0^t \tilde{\cR}(u(s;\nu))ds,
\end{equation}
where $\nu\in \mathcal{M}_1(\bT^d)$. Recall that the goal is to show 
\[
\frac{{\frak z}_t}{\sqrt{t}}\Rightarrow N(0,\sigma^2), \quad\quad \mbox{ as } t\to\infty.
\]
\gu{It is clear that ${\frak z}_2$ is a random variable that is finite almost surely, so we have ${\frak z}_2/\sqrt{t}\to0$ almost surely as $t\to\infty$.}
Hence, it suffices to consider 
\[
{\frak z}_t-{\frak z}_2=\int_2^t \int_{\bT^d} u(s,y;\nu)\xi(s,y)dyds-\tfrac12\int_2^t \tilde{\cR}(u(s;\nu))ds.
\]
Since $u(2;\nu)$ is a continuous function almost surely,  by the Markov property, we can assume that in \eqref{e.4181}, we start from $\nu(dx)=v(x)dx$ with some fixed $v\in D^\infty(\bT^d)$.

By \eqref{e.matdecom2}, we write 
\begin{equation}
\label{ztNt}
{\frak z}_t={\cal N}_t +{\frak r}_t,
\end{equation}
 with the martingale term
\begin{equation}\label{e.4182}
{\cal N}_t=
\int_0^t \left\la u(r;v)(1-\tfrac12\cD\chi(u(r;v))),dW(r)\right\ra_{L^2(\bT^d)},
\end{equation}
and the remainder term
\begin{equation}\label{e.4183}
{\frak r}_t=\tfrac12\chi(u(t;v))-\tfrac12\chi(v).
\end{equation}
The quadratic variation of $\mathcal{N}_t$ is
\[
\la \mathcal{N}\ra_t=
\int_0^t \cR\big(u(r;v)(1-\tfrac12\cD\chi(u(r;v)))\big)dr,
\]
which satisfies
\[
\la \mathcal{N}\ra_t \leq \int_0^t\|u(r;v)\|_{L^\infty(\bT^d)}^2(1+\tfrac12\|\cD\chi(u(r;v))\|_{L^\infty(\bT^d)})^2 dr.
\]
Applying Lemma~\ref{l.bdRv} and \eqref{e.4171}, we conclude that for any $p\in[1,+\infty)$,
\begin{equation}\label{e.bdNQVt}
\E[\la \mathcal{N}\ra_t^p] \leq Ct^p.
\end{equation}
Define
\begin{equation}\label{e.defsigma}
\sigma^2=\int_{\mathcal{M}_1(\bT^d)}\cR\big(u(1-\tfrac12\cD\chi(u))\big)\pi_\infty(du).
\end{equation}
The following proposition  shows the convergence of the (rescaled) quadratic variation:
\begin{proposition}\label{p.conQV}
As $t\to\infty$,
\[
\EE[|t^{-1}\la \mathcal{N}\ra_t-\sigma^2|^2]\to0.
\]
\end{proposition}

\begin{proof}
 For any $T>1$, define
\[
\la \mathcal{N}^T\ra_t=\int_0^t \cR\big(u(r;v)(1-\tfrac12\cD\chi(T,u(r;v)))\big)dr,
\]
and
\[
\sigma_T^2=\int_{\mathcal{M}_1(\bT^d)}\cR\big(u(1-\tfrac12\cD\chi(T,u))\big)\pi_\infty(du).
\]
By Lemma~\ref{l.conRN}, we have 
\[
\begin{aligned}
|\la \mathcal{N}^T\ra_t-\la \mathcal{N}\ra_t|\leq\int_0^t &\|u(r;v)\|_{L^\infty(\bT^d)}^2\|\cD\chi(T,u(r;v))-\cD\chi(u(r;v))\|_{L^\infty(\bT^d)}\\
&\times (1+\|\cD\chi(T,u(r;v))\|_{L^\infty(\bT^d)} +\|\cD\chi(u(r;v))\|_{L^\infty(\bT^d)})dr.
\end{aligned}
\]
Further applying \eqref{e.4171} and \eqref{010104-21}, we derive
\[
|\la \mathcal{N}^T\ra_t-\la \mathcal{N}\ra_t|\leq Ce^{-\lambda T}\int_0^t \|u(r;v)\|_{L^\infty(\bT^d)}^2(1+\|u(r;v)\|_{L^\infty(\bT^d)})dr.
\]
Similarly, we have 
\[
|\sigma_T^2-\sigma^2|\leq  Ce^{-\lambda T}\int_{\mathcal{M}_1(\bT^d)}\|u\|_{L^\infty(\bT^d)}^2(1+\|u\|_{L^\infty(\bT^d)})\pi_\infty(du).
\]
\tk{Proposition \ref{prop012211-22} implies that the right hand side is
finite.}

Thus, applying Lemma~\ref{l.bdRv}, we derive
\[
\sup_{t>0}\EE[t^{-2}|\la \mathcal{N}^T\ra_t-\la \mathcal{N}\ra_t|^2]+|\sigma_T^2-\sigma^2| \leq Ce^{-\lambda T},
\]
so we only need to show that for each $T$,  as $t\to\infty$, 
\begin{equation}\label{e.conQVT}
\EE[|t^{-1}\la \mathcal{N}^T\ra_t-\sigma_T^2|^2]\to0.
\end{equation}
Fix $T>0$. By Lemmas~\ref{l.contDe1}, \ref{l.conRN} and Corollary~\ref{cor013103-21}, the functional $\mathsf{F}:D^\infty(\bT^d)\to\R$
\[
\begin{aligned}
\mathsf{F}(u):=&\cR\big(u(1-\tfrac12\cD\chi(T,u))\big)\\
\end{aligned}
\]
satisfies 
\[
|\mathsf{F}(u_1)-\mathsf{F}(u_2)|\leq \|u_1-u_2\|_{L^\infty(\bT^d)} P(\|u_1\|_{L^\infty(\bT^d)},\|u_2\|_{L^\infty(\bT^d)}),
\]
where $P$ is some polynomial function. In other words, $\mathsf{F} $ is locally Lipschitz. For any $M>0$, let $\mathsf{F}_M:D^\infty(\bT^d)\to\R$ be a global Lipschitz function that is a cutoff of $\mathsf{F}$ in the sense that  $\mathsf{F}_M(u)=\mathsf{F}(u)$ if $\|u\|_{L^\infty(\bT^d)}\leq M$, and $|\mathsf{F}_M|\leq |\mathsf{F}|$. Applying Proposition~\ref{prop021203-21}, we have \[
\EE[|t^{-1}\int_0^t \mathsf{F}_M(u(r;v))dr-\sigma_{T,M}^2|^2]\to0,
\]
where $\sigma_{T,M}^2=\int_{\mathcal{M}_1(\bT^d)}\mathsf{F}_M(u)\pi_\infty(du)$. For the error induced by $M$, we have
\[
\begin{aligned}
&t^{-1} \int_0^t \EE|\mathsf{F}_M(u(r;v))-\mathsf{F}(u(r;v))|^2 dr\\
&\leq 4t^{-1}\int_0^t \EE\big[|\mathsf{F}(u(r;v))|^21_{\|u(r;v)\|_{L^\infty(\bT^d)}>M}\big]dr \leq \frac{C}{\sqrt{M}},
\end{aligned}
\]
where for the last estimate we have applied the H\"older inequality and Lemma~\ref{l.bdRv}.
Similarly, we can show $\sigma_{T,M}^2\to \sigma_T^2$ as $M\to\infty$. The proof is complete.
\end{proof}

The following two lemmas  combine to complete the proof of Theorem~\ref{t.partitionfunction}.
\begin{lemma}\label{l.418remainder}
As $t\to\infty$, $\frac{1}{\sqrt{t}}{\frak r}_t\to0$ in probability.
\end{lemma}
\begin{lemma}\label{l.418martingale}
As $t\to\infty$, $\frac{1}{\sqrt{t}}\mathcal{N}_t\Rightarrow N(0,\sigma^2)$ in distribution.
\end{lemma}

\begin{proof}[Proof of Lemma~\ref{l.418remainder}]
By \eqref{e.bdchi}, we have 
\[
|{\frak r}_t|\leq C(1+\|v\|_{L^\infty(\bT^d)}+\|u(t;v)\|_{L^\infty(\bT^d)}).
\]
Taking expectation and applying Lemma~\ref{l.bdRv}, we derive that $\EE|{\frak r}_t|\leq C(1+\|v\|_{L^\infty(\bT^d)})$. The proof is complete.
\end{proof}

\begin{proof}[Proof of Lemma~\ref{l.418martingale}]
By the same proof for \eqref{e.bdNQVt}, we have for any $t\geq s$
\[
\EE (\la \mathcal{N}\ra_t-\la \mathcal{N}\ra_s)^p\leq C(t-s)^p.
\]
Together with Proposition~\ref{p.conQV}, we conclude that the process $(\eps^2\la \mathcal{N}\ra_{t/\eps^2})_{t\geq0}$ converges in $C[0,\infty)$ to $(\sigma^2t)_{t\geq0}$, as $\eps\to0$. Then by the martingale central limit theorem \gu{\cite[Theorem IX.3.21]{jsbook}}, we derive that \[
\eps \mathcal{N}_{t/\eps^2}\Rightarrow \sigma B_t
\]
in $C[0,\infty)$ in distribution, where $B$ is a standard Brownian motion. The proof is complete.
\end{proof}


The following lemma guarantees the non-degeneracy of $\gamma,\sigma$.
\begin{lemma}\label{l.nondege1}
We have $\gamma,\sigma\in(0,\infty)$.
\end{lemma}

\begin{proof}
Recall that 
\begin{equation}
  \label{022211-22}
\begin{aligned}
\gamma=\tfrac12\int_{\mathcal{M}_1(\bT^d)} \cR(u)\pi_\infty(du),\quad\quad \sigma^2=\int_{\mathcal{M}_1(\bT^d)}\cR\big(u(1-\tfrac12\cD\chi(u))\big)\pi_\infty(du).
\end{aligned}
\end{equation}
By \eqref{e.4171}, Lemma~\ref{l.bdRv}, Theorem~\ref{t.geomeEr} and the
fact that 
$$
\cR(f)\leq \|f\|_{L^2(\bT^d)}^2 \leq \|f\|_{L^\infty(\bT^d)}
\quad\mbox{ for any $f\in D^\infty(\bT^d)$},
$$
we conclude that $\gamma,\sigma^2<\infty$.


Next, for any $u\in D^\infty(\bT^d)$, we let
$$
g=u(1-\tfrac12\cD\chi(u))\in L^\infty(\bT^d).
$$ 
Then (recall that
$\hat{r}_0 =1$, see \eqref{hr0})
\[
\cR(g)=\sum_{k\in\Z^d} \hat{r}_k |\hat{g}_k|^2\geq |\hat{g}_0|^2.
\]
Since 
\[
\hat{g}_0=\int_{\bT^d} g(x)dx=1-\tfrac12\la u,\cD\chi(u)\ra_{L^2(\bT^d)}=1,
\]
where the last equality comes from \eqref{e.vDv}, we conclude that
$\sigma^2\geq1$. \tk{Likewise, since $
\cR(u)=\sum_{k\in\Z^d} \hat{r}_k |\hat{u}_k|^2\geq |\hat{u}_0|^2=1
$, for any  $u\in D^\infty(\bT^d)$, we conclude immediately from
\eqref{022211-22} that $\gamma\ge 1/2$.}
\end{proof}

\appendix

\section{Fortet-Mourier metric on the space of Borel probability measures}
\label{s.FM}

Suppose that $(\mathbb X,{\rm d})$ is a  metric space.
Given a set $A\subset \mathbb X$  and a function $F:A\to\bbR$ we
denote
$$
\|F\|_{\infty,A}:=\sup_{x\in A}|F(x)|.
$$
If $A=\mathbb X$, we omit writing the subsript $A$ in the notation of
the norm.

Let ${\rm M}(\mathbb X) $ be the space of
all Borel signed measures on $\mathbb X$ with a finite total variation.
By ${\rm
  M}_1(\mathbb X) \subset {\rm M}(\mathbb X) $ we denote the subset consisting of all probability measures.

For any $F: \mathbb X\to \bbR$
we let
\begin{equation}
\label{FLip0}
\|F\|_{{\rm Lip}}:=\|F\|_\infty+ \sup\limits_{f\not=g,f,g\in \mathbb
  X}\frac{|F(f)-F(g)|}{{\rm d}(f,g)}.
\end{equation}
Denote by ${\rm Lip}(\mathbb X)$ the space of all
functions $F$ for which $\|F\|_{{\rm Lip}}<+\infty$.

Given $\Pi\in {\rm M}(\mathbb X) $, we define
\begin{equation}
\label{FM}
\|\Pi\|_{\mathrm{FM}}:=\sup\left\{\left|\int_{\mathbb X}Fd\Pi\right|:\, \|F\|_{{\rm Lip}}\le 1 \right\}.
\end{equation}
It is a norm on  ${\rm M}(\mathbb X)$, see e.g. \cite[Lemma 6,
p. 254]{dudley}. In fact, when restricted to ${\rm M}_1(\mathbb X)$ it defines a
metric, called the Fortet-Mourier metric:
$$
{\rm d}_{\mathrm{FM}}\Big(\Pi_1, \Pi_2\Big):=\|\Pi_1-\Pi_2\|_{\mathrm{FM}},\quad  \Pi_1,\Pi_2\in {\rm M}_1 (\mathbb X).
$$
 The topology of $({\rm M}_1(\mathbb
X), {\rm d}_{\mathrm{FM}})$
topology coincides with the topology of the weak convergence of
probability measures, see e.g. \cite[Theorem 3.2.2, p. 111]{bogachev}.
Suppose that $(\mathbb X,{\rm d})$ is complete
metric space, then the metric space $({\rm M}_1(\mathbb
X), {\rm d}_{\mathrm{FM}})$  is complete. If in addition $(\mathbb
X,{\rm d})$ is separable then $({\rm M}_1(\mathbb
X), {\rm d}_{\mathrm{FM}})$ is also 
separable. 
In this way we can define the norms $\|\cdot\|_{\mathrm{FM},{\cal M}_1}$ and  $\|\cdot\|_{\mathrm{FM},p}$ on the spaces $  {\rm
  M}_1({\cal M}_1(\bT^d)) $ and $  {\rm
  M}_1(D^p(\bT^d)) $ for any $p\in[1,+\infty)$.

\section{Technical lemmas on SHE}
\label{s.SHE}

Recall that $\U(t,x;y)$ is Green's function of the SHE, while
$\U(t,x;\nu)$ is the solution to the SHE with initial data $\nu$
belonging to $
\mathcal{M}_1(\bT^d)$. Let
\begin{equation}
\label{Ztn}
Z_t=\int_{\bT^d} \U(t,x;\nu)dx
\end{equation}
 be the partition function of the
directed polymer.

\subsection{Markov property of the endpoint distribution. Proof
  of Lemma \ref{l.markov}}

\label{MPd}

In the proof we  will write $\U(t;\nu,\omega)$ and $u(t;\nu,\omega)$ to emphasize the dependence of the solution on $\omega\in\Omega$.
Since $\U$ is the solution of a linear
equation, we conclude that
\[
\begin{aligned}
\U(t+s;\nu,\om)&=\U(s;\U(t;\nu,\om),\theta_{t,0}(\om))\\
&=\|\U(t;\nu,\om)\|_{L^1(\bT^d)}\U(s;u(t;\nu,\omega),\theta_{t,0}(\om)).
\end{aligned}
\]
 Therefore
\[
\begin{aligned}
&u(t+s;\nu,\omega)=\frac{{\U}(t+s;\nu,\om)}{\|{\U}(t+s;\nu,\om)\|_{L^1(\bT^d)}}\\
&=\frac{\U(s;u(t;\nu,\omega),\theta_{t,0}(\om))}{\|\U(s;u(t;\nu,\omega),\theta_{t,0}(\om))\|_{L^1(\bT^d)}}=u\Big(s; u(t;\nu,\omega),\theta_{t,0}(\omega)\Big),
\end{aligned}
\]
which completes the proof. 
\qed

\subsection{Some properties of  SHE}
\label{s.shewhite}

In this section, we recall and state several  results on the positive
and negative moments of the solution to SHE. We consider two cases:
either (i) $R(\cdot)\in \tk{C^\infty(\bT^d)}$ and $d\geq1$, or (ii) $R(\cdot)=\delta(\cdot)$ and $d=1$. All results are rather standard here, so we only sketch their proofs.

\begin{lemma}\label{l.estiGreen}
For any $T>0,p\in[1,\infty)$, there exists $C=C(T,p)>0$ such that for  $t\in(0,T],x\in\bT^d$, we have 
\begin{equation}\label{e.estiGreen1}
\EE[\U(t,x;0)^p] \leq \left\{\begin{array}{ll}
Cp_t(x)^p, & \mbox{ in case (i)}, \\
C [1+p_t(0)p_t(x)]^{p/2}, & \mbox{ in case (ii)}.\\
\end{array}
\right.
\end{equation}
Recall that $p_t(x)$ is the heat kernel on the torus, see \eqref{e.defptx}.
Futhermore, for each $\eps>0$ and $\delta\in(0,\tfrac12)$, we have 
\begin{equation}\label{e.modulusCon}
\sup_{t\in (\eps,\eps^{-1}), x\neq y\in \bT^d} \frac{\EE[|\U(t,x;0)-\U(t,y;0)|^p]}{|x-y|^{\delta p}} \leq C(\eps,\delta,p),
\end{equation}
and 
\begin{equation}\label{e.neMMGr}
\begin{aligned}
\sup_{t\in(\eps,\eps^{-1}),x\in\bT^d} \EE[\U(t,x;0)^{-p}] \leq C(\eps,p).
\end{aligned}
\end{equation}
\end{lemma}

\begin{proof}
For \eqref{e.estiGreen1}, case (i) comes from a direct application of
the Feynman-Kac formula. To show it in case (ii), one can invoke
\cite[Theorem 3.3]{davar2}, which was proved for the equation posed on
the whole space, however the proof applies verbatim to our
setting. 
 The estimate
\eqref{e.modulusCon} can be found e.g. in \cite[Chapter 3]{walsh}. For \eqref{e.neMMGr}, one can follow the same proof for
\cite[Corollary 4.9]{khoa}. Since the proofs are almost unchanged, we do
not reproduce them  here.
\end{proof}

\begin{lemma}\label{l.MMZ}
For any $T>0,p\in[1,\infty)$, there exists $C=C(T,p)$ such that 
\begin{equation}
\label{e.negammZ}
\EE[Z_t^p]+\EE[Z_t^{-p}]\leq C, \quad\quad t\in[0,T].
\end{equation}
\end{lemma}

\begin{proof}
Using \eqref{Ztn}, for the positive moments, we conclude by the Jensen
inequality that
\[
\EE[Z_t^p] \leq \int_{\bT^{d}}\EE\left(\int_{\bT^d}\U(t,x;y)dx\right)^p\nu(dy).
\]
For each fixed $y\in \bT^d$ and $t>0$, we know $\int_{\bT^d} \U(t,x;y)dx\stackrel{\rm law}{=}\U(t,0;\1)$, 
which implies 
\[
\EE[Z_t^p] \leq  \EE[\U(t,0;\1)^p].
\]
Similarly, for negative moments we have
\[
\EE[ Z_t^{-p} ]\leq \EE[\U(t,0;\1)^{-p}].
\]
\tk{The positive and negative moments bound on $\U(t,0;\1)$ is a standard
result on SHE, see e.g. \cite[Proposition 4.1]{davar} and
\cite[Theorem 2]{MN08} for the case (ii). We note that although   \cite[Theorem 2]{MN08} deals with the case of Dirichlet boundary
conditions, the proof carries out to the case of periodic boundary
conditions as well, see also \cite{khoa} for a different proof of the negative moment bound.}

\tk{In the case (i), i.e. $R(\cdot)\in C^\infty(\bT^d)$ and $d\geq1$, we can use
formula \eqref{031011-22} that allows us to derive
\begin{align*}
  &\EE[\U(t,0;\1)^{-p}]= \EE\Big\{\int_{\bT^d} p_t(y)\E_{y,0}^t[\exp\left\{\int_0^t
    \xi(s,B_s)ds-\frac12R(0)t\right\}]dy\Big\}^{-p}\\
  &
   \le \EE\Big\{\int_{\bT^d} p_t(y)\E_{y,0}^t[\exp\left\{-p\int_0^t
    \xi(s,B_s)ds+\frac p2R(0)t\right\}]dy\Big\}  .
\end{align*}
The last estimate is a consequence of the Jensen inequality. Denoting
$\tilde{\U}(t,0;\1)$ the solution of the SHE corresponding to the
noise $-p\xi(t,x)$ we conclude that
\begin{align*}
  \EE[\U(t,0;\1)^{-p}]\le e^{p(p+1)R(0)T/2}\EE[\tilde{\U}(t,0;\1)]
    ,\quad \mbox{for }t\in[0,T],
\end{align*}
and the conclusion of the lemma for the negative moments follows from the estimate of positive moments.}
\end{proof}

\begin{lemma}\label{l.mmbdu}
For any $T>0,p\in[1,\infty)$, there exists $C=C(T,p)$ such that for all  $t\in(0,T]$, we have 
\[
\sup_{\nu\in\mathcal{M}_1(\bT^d)}\big(\EE[\|\U(t;\nu)\|_{L^{2p}(\bT^d)}^{2p}]+\EE[\|u(t;\nu)\|_{L^{2p}(\bT^d)}^{2p}]\big)\leq C (1+ t^{-(p-1/4)}). 
\]
\end{lemma}

\begin{proof}
First, we have 
\[
\begin{aligned}
\|u(t;\nu)\|_{L^{2p}(\bT^d)}^{2p} 
=Z_t^{-2p} \int_{\bT^d} \U(t,x;\nu)^{2p}dx.
\end{aligned}
\]
An application of the H\"older inequality leads to
\[
\begin{aligned}
\EE[\|u(t;\nu)\|_{L^{2p}(\bT^d)}^{2p}] \leq &\sqrt{\EE[Z_t^{-4p}]\int_{\bT^d} \EE[\U(t,x;\nu)^{4p}] dx}\\
\leq &\sqrt{\EE[Z_t^{-4p}]\int_{\bT^{2d}} \EE[\U(t,x;y)^{4p}] dx\,\nu(dy)}\\
=&\sqrt{\EE[Z_t^{-4p}]\int_{\bT^{d}} \EE[\U(t,x;0)^{4p}] dx}.
\end{aligned}
\]
\tk{In the last equality we have used the spatial stationarity of the random
  field $(x,y)\mapsto \U(t,x;y)$. 
Thanks to \eqref{e.negammZ} and
\eqref{e.estiGreen1} we conclude that the utmost right hand side of
the expression above can be estimated by
\[
\begin{aligned}
C \left\{\int_{\bT^d}\Big[1+\Big(p_t(0)p_t(x)\Big)^{2p}\Big]dx\right\}^{1/2}\le C \big\{1+[p_t(0)]^{4p-1}\big\}^{1/2}.
\end{aligned}
\]
This leads to the conclusion of the lemma  for $u$ in case (ii). The proof for
$\U$ is similar. Finally the case (i) can be handled in the similar way.}
\end{proof}

\begin{lemma}\label{l.bdQVZ}
For any $p\in[1,\infty)$, we have $\EE\la Z\ra_t^p\leq C(t,p)$.
\end{lemma}

\begin{proof}
Recall that 
\[
\la Z\ra_t=\int_0^t \int_{\bT^{2d}} \U(s,y)\U(s,y')R(y-y')dydy'ds.
\]
In case (i), we have 
\[
\la Z\ra_t \leq \|R\|_{L^\infty(\bT^d)}\int_0^t Z_s^2ds,
\]
then it suffices to apply Lemma~\ref{l.MMZ} to complete the proof. In case (ii), we have $\la Z\ra_t= \int_0^t \|\U(s;\nu)\|_{L^2(\bT^d)}^2 ds$. Applying Lemma~\ref{l.mmbdu}, we have 
\[
\begin{aligned}
\big(\gu{\EE} \la Z\ra_t^p\big)^{1/p} &\leq \int_0^t \big(\gu{\EE} \|\U(s;\nu)\|_{L^2(\bT^d)}^{2p}\big)^{1/p}ds \\
&\gu{\leq\int_0^t \big(\gu{\EE} \|\U(s;\nu)\|_{L^{2p}(\bT^d)}^{2p}\big)^{1/p}ds} \leq  C\int_0^t (1+s^{-(1-\frac{1}{4p})})ds,
\end{aligned}
\]
which concludes the proof. 
\end{proof}

\begin{lemma}\label{l.lshe1}
For any $T>0,p\in[1,\infty)$, there exists a constant $C=C(T,p)$ such that 

(i) for any $0\leq v\in L^2(\bT^d)$, we have 
\begin{equation}\label{e.she11}
\sup_{t\in[0,T]}\EE[\|\U(t;v)\|_{L^2(\bT^d)}^{2p}] \leq C\|v\|_{L^2(\bT^d)}^{2p},
\end{equation}

(ii) for any $0\leq v\in L^\infty(\bT^d)$, we have
\begin{equation}\label{e.she12}
\sup_{t\in[0,T]}\EE[\|\U(t;v)\|_{L^\infty(\bT^d)}^{2p}] \leq C\|v\|_{L^\infty(\bT^d)}^{2p},
\end{equation}

(iii) for any $v\in D^2(\bT^d)$, we have 
\begin{equation}\label{e.she13}
\sup_{t\in[0,T]}\EE[\|u(t;v)\|_{L^2(\bT^d)}^{2p}] \leq C\|v\|_{L^2(\bT^d)}^{2p},
\end{equation}

(iv) for any $v\in D^\infty(\bT^d)$, we have 
\begin{equation}\label{e.she14}
\sup_{t\in[0,T]}\EE[\|u(t;v)\|_{L^\infty(\bT^d)}^{2p}] \leq C\|v\|_{L^\infty(\bT^d)}^{2p},
\end{equation}

(v) for any $h\in L^2(\bT^d), v\in D (\bT^d)$, we have 
\begin{equation}\label{e.she15}
\sup_{t\in[0,T]}\EE[\|U(t;h,v)\|_{L^2(\bT^d)}^{2p}]\leq C\|h\|_{L^2(\bT^d)}^{2p}.
\end{equation}
\end{lemma}

\begin{proof}
By Lemma~\ref{l.MMZ} and the H\"older inequality, \eqref{e.she13}, \eqref{e.she14} and \eqref{e.she15} are consequences of \eqref{e.she11} and  \eqref{e.she12}. So we focus on \eqref{e.she11} and \eqref{e.she12}. 

By comparison principle \cite[Theorem 1.3]{chenhuang}, we have 
\[
\|\U(t;v)\|_{L^\infty(\bT^d)} \leq \|v\|_{L^\infty(\bT^d)} \|\U(t;\1)\|_{L^\infty(\bT^d)}.
\]
By applying Lemma~\ref{l.U1} formulated below, we conclude \eqref{e.she12}.

To prove \eqref{e.she11}, the case of $R(\cdot)\in \tk{C^\infty(\bT^d)}$ can be studied by the Feynman-Kac formula in a straightforward way, so we focus on the  $1+1$ spacetime white noise. By the mild formulation, we have
\[
\U(t,x;v)=p_t\star v(x)+\int_0^t\int_{\bT^d} p_{t-s}(x-y)\U(s,y;v)\xi(s,y)dyds=:p_t\star v(x)+M_t(x),
\]
so 
\[
\int_{\bT^d}\U(t,x;v)^2dx \leq 2\int_{\bT^d} |p_t\star v(x)|^2 dx+2\int_{\bT^d} M_t(x)^2 dx.
\]
The first term on the r.h.s. is bounded from above by $2\|v\|_{L^2(\bT^d)}^2$. For the second term, applying \cite[Proposition 4.4]{davar1}, we have
\[
\|M_t(x)\|_{L^{2p}(\Omega)}^2 \leq C\int_0^t\int_{\bT^d}p_{t-s}(x-y)^2\|\U(s,y;v)\|_{L^{2p}(\Omega)}^2dyds.
\]
By Lemma~\ref{l.estiGreen}, we have 
\[
\begin{aligned}
\|\U(s,y;v)\|_{L^{2p}(\Omega)}=&\left\|\int_{\bT^d} \U(s,y;z)v(z)dz\right\|_{L^{2p}(\Omega)}
\leq \int_{\bT^d} \|\U(s,y;z)\|_{L^{2p}(\Omega)}v(z)dz \\
\leq &C\|v\|_{L^2(\bT^d)}\sqrt{\int_{\bT^d} (1+p_s(0)p_s(y-z))dz}\\
\leq &C\|v\|_{L^2(\bT^d)}(1+s^{-1/4}),
\end{aligned} 
\]
which implies
\[
\|M_t(x)\|_{L^{2p}(\Omega)}^2 \leq C \|v\|_{L^2(\bT^d)}^2\int_0^t \frac{1}{\sqrt{t-s}}(1+\frac{1}{\sqrt{s}})ds \leq C\|v\|_{L^2(\bT^d)}^2,
\]
and 
\[
\EE|\int_{\bT^d} M_t(x)^2 dx|^p \leq \EE\int_{\bT^d}M_t(x)^{2p} dx \leq C\|v\|_{L^2(\bT^d)}^{2p},
\]
which completes the proof.
\end{proof}

\begin{lemma}\label{l.uconti}
Assume $v\in D^2(\bT^d)$, then for any $s\geq0$ and $p\in[1,\infty)$
\[
\lim_{t\to s} \EE[\|u(t;v)-u(s;v)\|_{L^2(\bT^d)}^p]=0.
\]
\end{lemma}

\begin{proof}
We have 
\[
\begin{aligned}
&\|u(t;v)-u(s;v)\|_{L^2(\bT^d)} \\
&\leq \|\U(t;v)-\U(s;v)\|_{L^2(\bT^d)}\big(Z_t^{-1} +\|\U(s;v)\|_{L^2(\bT^d)}Z_t^{-1}Z_s^{-1}\big).
\end{aligned}
\]
To complete the proof it suffices to use the fact  that
\[
\EE[\|\U(t;v)-\U(s;v)\|_{L^2(\bT^d)}^p]\to0,  \quad \mbox{ as } t\to s
\]
for any $p\geq1$, together with Lemma~\ref{l.MMZ} and estimate \eqref{e.she11}.
\end{proof}

\begin{lemma}\label{l.U1}
For any $T\geq0,p\in[1,\infty)$, we have 
\[
\EE\big[\sup_{t\in[0,T],x\in\bT^d} \U(t,x;\1)^p\big]+\EE\big[\sup_{t\in[0,T],x\in\bT^d} \U(t,x;\1)^{-p}\big] \leq C(T,p).
\]
\end{lemma}

\begin{proof}
 The proof is similar to   Lemma~\ref{l.ulbd}, so we do not repeat it here.
\end{proof}

\subsection{Proof of Lemmas~\ref{l.ulbd}, \ref{l.conuN}}
\label{s.sheOther}


\begin{proof}[Proof of Lemma~\ref{l.ulbd}]
From \eqref{e.defZts}, we have $\cZ_{t,0}(x,y)=\U(t,x;y)$, so the goal reduces to prove that there exists $C>0$ only depending on $d,p,R$ such that
\begin{equation}\label{e.ubgreen}
\EE[\sup_{x,y\in\bT^d} \U(1,x;y)^p]\leq C, \quad\quad \EE[(\inf_{x,y\in\bT^d} \U(1,x;y))^{-p}]\leq C.
\end{equation}
For the second inequality, by writing
\[
(\inf_{x,y\in\bT^d} \U(1,x;y))^{-p}=\sup_{x,y\in\bT^d} \U(1,x;y)^{-p},
\]
it suffices to prove 
\begin{equation}\label{e.supw}
\E\big[\sup_{w\in I} X(w)^p\big]\leq C
\end{equation}
with $w=(x,y)$, $I=\bT^{2d}$, and $X(w)=\U(1,x;y)$ or $\U(1,x;y)^{-1}$. By a well-known chaining argument, see e.g. \cite[Proposition 5.8]{davar}, we need to show that there exists $\alpha>0$ such that for any $p\geq 2$,
\begin{equation}\label{e.chaincondition}
\begin{aligned}
&\EE[X(w)^p]<\infty \quad\quad \text{ for some } w\in I,\\
&\sup_{w_1,w_2\in I,w_1\neq w_2} \EE\left[\frac{|X(w_1)-X(w_2)|^p}{|w_1-w_2|^{p\alpha}}\right]<\infty.
\end{aligned}
\end{equation}
This comes from the following estimates:
\begin{equation}\label{e.standardSHE}
\begin{aligned}
&\sup_{x,y\in \bT^d}\EE[\U(1,x;y)^p]+\sup_{x,y\in\bT^d}\EE[\U(1,x;y)^{-p}]\leq C,\\
&\sup_{x\in\bT^d} \EE[ |\U(1,x;y_1)-\U(1,x;y_2)|^p]\leq C|y_1-y_2|^{\delta p},\\
&\sup_{y\in\bT^d} \EE[ |\U(1,x_1;y)-\U(1,x_2;y)|^p]\leq C|x_1-x_2|^{\delta p},
\end{aligned}
\end{equation}
which are direct consequences of Lemma~\ref{l.estiGreen} and the fact that for fixed $x,y\in\bT^d$, we have 
\[
\U(1,x;y)\stackrel{\text{law}}{=} \U(1,x-y;0).
\]
\end{proof}

\begin{proof}[Proof of Lemma~\ref{l.conuN}]
We only sketch the proof for the case of $1+1$ spacetime white noise. First, we have 
\begin{equation}\label{e.bertini}
\sup_{t\in[0,T]}\EE[\|\U_N(t;v)-\U(t;v)\|_{L^2(\bT^d)}^p]\to0,
\end{equation}
which is a consequence of \cite[Theorem 2.2]{bertini1995stochastic}.
While the result of \cite{bertini1995stochastic} is formulated for the equation posed on the whole space, the proof applies verbatim to our setting. Denote $Z_N(t)=\|\U_N(t;v)\|_{L^1(\bT^d)}$, then as in Lemma~\ref{l.MMZ}, we also have
\[
\EE[Z_N(t)^p]+\EE[Z_N(t)^{-p}]\leq C(T,p), \quad\quad t\in[0,T].
\]
The proofs for all cases are similar, so we take $\EE[\|u_N(t;v)-u(t;v)\|_{L^2(\bT^d)}^p]$ as an example. We have
\[
\begin{aligned}
\|u_N(t;v)-u(t;v)\|_{L^2(\bT^d)}\leq &Z_N(t)^{-1}\|\U_N(t;v)-\U(t;v)\|_{L^2(\bT^d)}\\
&+(Z_N(t)Z_t)^{-1}\|\U(t;v)\|_{L^2(\bT^d)}\|\U_N(t;v)-\U(t;v)\|_{L^2(\bT^d)},
\end{aligned}
\]
then it suffices to apply H\"older inequality and \eqref{e.bertini} to conclude the proof.
\end{proof}

\section{\gu{Proof of the invariance of $\pi_\infty$}}

\label{appC}

In this section, we prove the invariance of the measure $\pi_\infty$, which does not immediately come from \eqref{e.452}. 
We first state the following proposition:
\begin{proposition}
  \label{cont}
  Suppose that $f_1,\ldots,f_N\in C_0^\infty(\R)$,
  $x_1,\ldots,x_N\in\bT^d$. Let $F: D^\infty(\bT^d)\to\R$ be given by
  \begin{equation}
    \label{C02}
    F(v):=\prod_{i=1}^N f_i\Big(v(x_i)\Big), \quad v\in
    D^\infty(\bT^d).
  \end{equation}
  Then, under the assumptions on the covariance
  $R(\cdot)$ made in Section \ref{sec1.1}, for each $t>0$, there exists
  a constant $C>0$ such that
  \begin{equation}
\label{C01}
\Big|{\cal P}_tF(v_2) -{\cal P}_tF(v_1)\Big| \leq
C\|v_2-v_1\|_{L^\infty(\bT^d)},\quad v_1,v_2\in D^\infty(\bT^d).
\end{equation}
\end{proposition}

The above result implies, in particular, that ${\cal P}_tF\in
C_b\Big(D^\infty(\bT^d)\Big)$ for any $F$ as in \eqref{C02}. Fix any $t>0$. By
\eqref{e.452} we conclude that
$$
\lim_{s\to\infty}\int_{D^\infty(\bT^d)}{\cal P}_{s}\Big({\cal P}_{t}
F\Big)d\pi_\infty=\int_{D^\infty(\bT^d)} {\cal P}_{t} Fd\pi_\infty=\int_{D^\infty(\bT^d)}Fd(\pi_\infty {\cal P}_{t}) .
$$
On the other hand, the utmost left hand side equals
$$
\lim_{s\to\infty}\int_{D^\infty(\bT^d)}{\cal P}_{s+t}
Fd\pi_\infty=\int_{D^\infty(\bT^d)} Fd\pi_\infty.
$$
We conclude therefore the equality
$$
 \int_{D^\infty(\bT^d)} Fd\pi_\infty=\int_{D^\infty(\bT^d)}Fd(\pi_\infty {\cal P}_{t} )
 $$
 for all $F$ of the form \eqref{C02}. This, in turn implies that
 $\pi_\infty=\pi_\infty {\cal P}_{t}$ and invariance of $\pi_\infty$
 follows. The only thing yet to be shown is Proposition \ref{cont}.

 \subsubsection*{Proof of Proposition \ref{cont}}

We consider only the case $d=1$ and $R=\delta$. To simplify the
notations, we assume also that $N=1$ and $F(v)=f\Big(v(x)\Big)$ for some
$f\in C_0^\infty(\R)$ and $x\in\bT$. Assume that $v_1,v_2\in
D^\infty(\bT)$, then  
\begin{align*}
\Big|{\cal P}_tF(v_2) -{\cal P}_tF(v_1)\Big| \le
               \|f'\|_{L^\infty(\bT)}\EE\Big|u(t,x;v_2) -u(t,x;v_1)
  \Big|\le {\rm I}+{\rm II},
\end{align*}
with
\begin{align*}
  &
    {\rm I}:= \|f'\|_{L^\infty(\bT)}\EE\left[\left(\int_{\bT}\U(t,x';v_2)dx'\right)^{-1}\Big|\U(t,x;v_2) -\U(t,x;v_1)
    \Big|\right],\\
  &
    {\rm II}:= \|f'\|_{L^\infty(\bT)}\EE\left[\U(t,x;v_1)\Big\|\U(t;v_2) -\U(t;v_1)
    \Big\|_{L^1(\bT)}\prod_{i=1}^2\left(\int_{\bT}\U(t,x';v_i)dx'\right)^{-1}\right].
\end{align*}
To estimate ${\rm I}$, we use the H\"older inequality together with
\eqref{e.neMMGr} and conclude that for $p\in(1,2)$ and $1/p+1/p'=1$
there exists $C>0$ such that
\begin{align*}
  &
    {\rm I}\le
    \|f'\|_{L^\infty(\bT)}\left\{\EE\left[\left(\int_{\bT}\U(t,x';v_2)dx'\right)^{-p}
    \right]\right\}^{1/p}\left\{\EE\left[\Big|\U(t,x;v_2) -\U(t,x;v_1)
    \Big|^{p'}\right]\right\}^{1/p'}\\
  &
    \le C\left\{\EE\left[\Big|\U(t,x;v_2) -\U(t,x;v_1)
    \Big|^{p'}\right]\right\}^{1/p'}.
\end{align*}
Invoking \cite[Proposition 4.1]{davar} we conclude that
\begin{align*}
    {\rm I}\le  C\|v_2-v_1\|_{L^\infty(\bT)},\quad v_1,v_2\in D^\infty(\bT).
\end{align*}
The estimate of  ${\rm II}$ is simiilar. We start with H\"older
inequality together with \eqref{e.estiGreen1} and 
\eqref{e.neMMGr} to obtain that
\begin{align*}
    {\rm II} 
    \le C\left\{\EE\left[\Big\|\U(t;v_2) -\U(t;v_1)
    \Big\|_{L^1(\bT)}^{p'}\right]\right\}^{1/p'}
\end{align*}
for some $p'\ge2$. Applying again \cite[Proposition 4.1]{davar} we
conclude
\begin{align*}
    {\rm II}\le  C\|v_2-v_1\|_{L^\infty(\bT)},\quad v_1,v_2\in D^\infty(\bT).
\end{align*}
This ends the proof of the proposition.\qed


\begin{thebibliography}{10}

  

  \bibitem{amir2011probability}
{\sc G.~Amir, I.~Corwin, and J.~Quastel}, {\em Probability distribution of the
  free energy of the continuum directed random polymer in 1+ 1 dimensions},
  Communications on pure and applied mathematics, 64 (2011), pp.~466--537.
  
  \bibitem{zhipeng1}
  {\sc J.~Baik and Z.~Liu}, {\em TASEP on a ring in sub-relaxation time scale}, Journal of Statistical Physics 165.6 (2016), pp.~1051--1085.
  
 \bibitem{zhipeng2}
  {\sc J.~Baik and Z.~Liu}, {\em Fluctuations of TASEP on a ring in relaxation time scale}, Communications on Pure and Applied Mathematics 71.4 (2018), pp.~747--813.
  
  \bibitem{zhipeng}
  {\sc J.~Baik, Z.~Liu and G. L. F. Silva}, {\em
  Limiting one-point distribution of periodic TASEP}, arxiv preprint arXiv:2008.07024 (2020).
  
  
  
  

\bibitem{bakhtin1}
{\sc Y.~Bakhtin and D.~Seo}, {\em Localization of directed polymers in continuous space}, Electronic Journal of Probability,  2020;25.



\bibitem{balazs2011fluctuation}
{\sc M.~Bal{\'a}zs, J.~Quastel, and T.~Sepp{\"a}l{\"a}inen}, {\em Fluctuation
  exponent of the KPZ/stochastic Burgers equation}, Journal of the American
  Mathematical Society, 24 (2011), pp.~683--708.



\bibitem{baco}
{\sc G.~Barraquand, and I.~Corwin}, {\em Stationary measures for the log-gamma polymer and KPZ equation in half-space}, arXiv preprint arXiv:2203.11037 (2022).


\bibitem{barraquand2021steady}
{\sc G.~Barraquand and P.~L. Doussal}, {\em Steady state of the KPZ equation on
  an interval and Liouville quantum mechanics}, arXiv preprint
  arXiv:2105.15178,  (2021).



\bibitem{bates2016endpoint}
{\sc E.~Bates and S.~Chatterjee}, {\em The endpoint distribution of directed
  polymers}, Ann. Probab., 48  (2020), pp.~817--871.

\bibitem{bertini1995stochastic}
{\sc L.~Bertini and N.~Cancrini}, {\em The stochastic heat equation:
  Feynman-Kac formula and intermittence}, Journal of statistical Physics, 78
  (1995), pp.~1377--1401.
  
  
  \bibitem{bertini}
 {\sc L.~Bertini and G.~Giacomin}, {\em Stochastic Burgers and KPZ
   equations from particle systems}, Communications in mathematical
 physics 183.3 (1997), pp.~571--607.


\bibitem{billingsley}   Billingsley, P. {\em Convergence of probability measures.} Second edition. Wiley Series in Probability and Statistics: Probability and Statistics. A Wiley-Interscience Publication. John Wiley \& Sons, Inc., New York, 1999. 


\bibitem{bogachev}   {\sc  V. I. Bogachev}, {\em Weak convergence of measures.} Mathematical Surveys and Monographs, 234. American Mathematical Society, Providence, RI, 2018.


\bibitem{borodin}
{\sc A.~Borodin, I.~Corwin and P.~Ferrari}, {\em Free energy fluctuations for directed polymers in random media in 1+ 1 dimension},  Communications on Pure and Applied Mathematics, 67.7 (2014), pp.~1129--1214.


\bibitem{borodin1}
{\sc A.~Borodin, I.~Corwin, P.~Ferrari and B.~Vet\"o}, {\em Height fluctuations for the stationary KPZ equation}, Mathematical Physics, Analysis and Geometry, 18(1) (2015), pp.~1--95.

\bibitem{mukherjee}
{\sc Y.~Br\"oker and C.~Mukherjee}, {\em Localization of the Gaussian multiplicative chaos in the Wiener space and the stochastic heat equation in strong disorder}, Annals of Applied Probability, 29(6) (2019), pp.~3745--3785.


\bibitem{brunet}
{\sc E.~Brunet}, {\em Fluctuations of the winding number of a directed polymer in a random medium}, Physical Review E 68.4 (2003): 041101.


\bibitem{brunet1}
{\sc E.~Brunet and D.~Bernard}, {\em Probability distribution of the free energy of a directed polymer in a random medium}, Physical Review E 61.6 (2000): 6789.


\bibitem{brunet2}
{\sc E.~Brunet and D.~Bernard}, {\em Ground state energy of a non-integer number of particles with $\delta$ attractive interactions}, Physica A: Statistical Mechanics and its Applications, 279(1-4) (2000), pp.~398--407.


  \bibitem{BKWW21}
{\sc W.~Bryc, A.~Kuznetsov, Y.~Wang, and J.~Wesolowski}, {\em {Markov processes
  related to the stationary measure for the open KPZ equation}}, July 2021,
  arXiv preprint \href{https://arxiv.org/abs/2105.03946v2}{2105.03946v2}.




 
   \bibitem{CRS18}
 {\sc F.~Caravenna, R.~Sun and N.~Zygouras},
 { \em The two-dimensional KPZ equation in the entire subcritical regime},
 Ann. Probab. 48 (2020),  pp.~1086--1127.
  
  \bibitem{chatterjee2018constructing}
{\sc S.~Chatterjee and A.~Dunlap}, {\em Constructing a solution of the $(2+ 1)
  $-dimensional {KPZ} equation}, Ann. Probab. 48 (2020),
pp.~1014–-1055.


\bibitem{chendalang} {\sc L.~Chen; R. C. Dalang,} {\em Hölder-continuity for the nonlinear stochastic heat equation with rough initial conditions.} Stoch. Partial Differ. Equ. Anal. Comput. 2 (2014), no. 3, 316-352.
  
\bibitem{chenhuang}
{\sc L.~Chen and J.~Huang}, {\em Comparison principle for stochastic heat equation on $\mathbb {R}^{d} $}, The Annals of Probability 47.2 (2019), pp.~989--1035.


 %
\bibitem{comets2017directed}
{\sc F.~Comets}, {\em Directed polymers in random environments (Probability in
  Saint-Flour, Lecture Notes in Mathematics 2175)}, Berlin: Springer, 2017.
  
  
  
  \bibitem{cometskpz}
  {\sc F.~Comets, C.~Cosco and C.~Mukherjee}, {\em Space-time fluctuation of the Kardar-Parisi-Zhang equation in $d\geq3$ and the Gaussian free field}, arXiv preprint arXiv:1905.03200 (2019).
  
  
  \bibitem{davar2}
  {\sc
  D.~Conus, M.~Joseph, D.~Khoshnevisan and S. Y.~Shiu}, {\em Initial measures for the stochastic heat equation}, Annales de l'IHP Probabilit\'es et statistiques, Vol. 50, No. 1, (2014) pp.~136--153.

\bibitem{corwin2012kardar}
{\sc I.~Corwin}, {\em The Kardar--Parisi--Zhang equation and universality
  class}, Random matrices: Theory and applications, 1 (2012), p.~1130001.
  
  
  \bibitem{CK21}
{\sc I.~Corwin and A.~Knizel}, {\em {Stationary measure for the open KPZ
  equation}}, Apr. 2021, arXiv preprint
  \href{https://arxiv.org/abs/2103.12253v2}{2103.12253v2}.


  
  \bibitem{haoshen1}
{\sc I.~Corwin,  and H.~Shen}, {\em Open ASEP in the weakly asymmetric regime}, Communications on Pure and Applied Mathematics 71.10 (2018): 2065-2128.
  
\bibitem{haoshen}
{\sc I.~Corwin and H.~Shen},
{\em Some recent progress in singular stochastic PDEs},  Bulletin of the AMS. 57.3 (2020), pp.~409--454.


\bibitem{cosco}
{\sc C.~Cosco, S.~Nakajima and M.~Nakashima}, {\em Law of large numbers and fluctuations in the
sub-critical and $L^2$ regions for SHE and KPZ equation in dimension $d\geq 3$}, arXiv preprint arXiv: 2005.12689v1 (2020).

  
  
\bibitem{daza}  {\sc G. Da Prato, J. Zabczyk}, {\em Stochastic equations in infinite dimensions.} Second edition. Encyclopedia of Mathematics and its Applications, 152. Cambridge University Press, Cambridge, 2014. 



  
\bibitem{dudley}  {\sc R. M. Dudley}, {\em Convergence of Baire measures.} Studia Math. 27 (1966), 251–268. 



\bibitem{duncan}
{\sc D.~Duncan, J.~Ortmann and B.~Vir\'ag}, {\em The directed landscape}, arXiv preprint arXiv:1812.00309 (2018).
  
  





\bibitem{dunlap}
{\sc A.~Dunlap, Y.~Gu, and T.~Komorowski}, {\em Fluctuation exponents of the KPZ equation on a large torus}, arXiv preprint arXiv:2111.03650 (2021).



  \bibitem{kpz1}
  {\sc  A.~Dunlap, Y.~Gu, L.~Ryzhik, and O.~Zeitouni}, {\em Fluctuations of the solutions to the KPZ equation in dimensions three and higher}, Probability Theory and Related Fields, 176 (2020), pp.~1217--1258.
  

%
%
   



\bibitem{funaki}
{\sc
T.~Funaki and J.~Quastel}, {\em KPZ equation, its renormalization and invariant measures}, Stochastic Partial Differential Equations: Analysis and Computations 3.2 (2015), pp.~159--220.

  
  \bibitem{kpz2}
{\sc Y.~Gu},  {\em Gaussian fluctuations from the 2D KPZ equation}, Stochastics and Partial Differential Equations: Analysis and Computations, 8 (2020), pp.~150--185.
  
    \bibitem{gubinelli2015paracontrolled}
{\sc M.~Gubinelli, P.~Imkeller, and N.~Perkowski}, {\em Paracontrolled
  distributions and singular PDEs}, in Forum of Mathematics, Pi, vol.~3,
  Cambridge University Press, 2015.
  
   \bibitem{gubinelli2017kpz}
{\sc M.~Gubinelli and N.~Perkowski}, {\em {KPZ} reloaded}, Communications in
  Mathematical Physics, 349 (2017), pp.~165--269.
  
  
 \bibitem{perkowski}
 {\sc M.~Gubinelli and N.~Perkowski}, {\em The infinitesimal generator of the stochastic Burgers equation}, Probability Theory and Related Fields 178.3 (2020), pp.~1067--1124.
  
  

\bibitem{hairer2013solving}
  {\sc M.~Hairer}, {\em Solving the {KPZ}
equation}, Ann. Math., {\bf 178}
  (2013), pp.~559--664.

\bibitem{hairer2014theory}
\leavevmode\vrule height 2pt depth -1.6pt width 23pt, {\em A theory of
  regularity structures}, Inv. Math., {\bf 198} (2014), pp.~269--504.
  


\bibitem{khoa}
{\sc Y.~Hu and K.~L\^e}, {\em Asymptotics of the density of parabolic Anderson random fields}, arXiv preprint arXiv:1801.03386 (2018).


\bibitem{sasamoto}
{\sc T.~Imamura, M.~Mucciconi, and T.~Sasamoto}, {\em Solvable models in the KPZ class: approach through periodic and free boundary Schur measures}, arXiv preprint arXiv:2204.08420 (2022).

\bibitem{jsbook}
{\sc J.~Jacob and A. N.~Shiryaev}, {\em Limit theorems for stochastic processes}, Second edition. Grundlehren
der Mathematischen Wissenschaften, 288. Springer-Verlag, Berlin, 2003.

\bibitem{kardar1986dynamic}
{\sc M.~Kardar, G.~Parisi, and Y.-C. Zhang}, {\em Dynamic scaling of growing
  interfaces}, Physical Review Letters, 56 (1986), p.~889.




%

  
\bibitem{davar1}
{\sc
D.~Khoshnevisan},  {\em Analysis of stochastic partial differential equations}, Vol. 119. American Mathematical Soc., 2014.
  
  
  \bibitem{davar} {\sc D.~Khoshnevisan, K.~Kim, C.~Mueller and  S.-Y.~Shiu}, {\em
  Dissipation in parabolic SPDEs}. Journal of Statistical Physics volume 179, pages 502-534 (2020).
  

\bibitem{alisa}
{\sc A.~Knizel,  and K.~Matetski}, {\sc The strong Feller property of the open KPZ equation}, arXiv preprint arXiv:2211.04466 (2022).

%
%
%




 


\bibitem{kupiainen2016renormalization}
{\sc A.~Kupiainen}, {\em Renormalization group and stochastic
  {PDEs}}, in Annales
  Henri Poincar{\'e}, vol.~17, Springer, 2016, pp.~497--535.
  
  \bibitem{nikos}
  {\sc D.~Lygkonis and N.~Zygouras}, {\em Edwards-Wilkinson fluctuations for the directed polymer in the
full $L^2-$regime for dimensions $d\geq3$}, arXiv preprint arXiv: 2005.12706 (2020).

\bibitem{magnen2017diffusive}
{\sc J.~Magnen and J.~Unterberger}, {\em The Scaling Limit of the {KPZ} Equation in Space Dimension 3 and Higher}, 
Jour. Stat. Phys., {\bf 171} (2018), pp.~543--598.

\bibitem{kpzfix}
{\sc K.~Matetski, J.~Quastel and D.~Remenik}, {\em The KPZ fixed
  point}, arXiv preprint arXiv:1701.00018, (2017).

\bibitem{MN08}
{\sc C.~Mueller and D.~Nualart}, {\em Regularity of the density for the stochastic heat equation.} Electr. J. Probab. 13(74) (2008), pp.~2248–-2258.


\bibitem{shalin}
{\sc S.~Parekh}, {\em The KPZ limit of ASEP with boundary}, Communications in Mathematical Physics 365.2 (2019): 569-649.

\bibitem{shalin1}
{\sc S.~Parekh}, {\em Ergodicity results for the open KPZ equation
}, arXiv preprint arXiv:2212.08248, (2023).

\bibitem{quastelintroduction}
{\sc J.~Quastel}, {\em Introduction to KPZ}, Current developments in mathematics, 2011(1).


  \bibitem{quastelkpz}
  {\sc J.~Quastel and S.~Sarkar}, {\em  Convergence of exclusion processes and KPZ equation to the KPZ fixed point}, 
 arXiv preprint arXiv:2008.06584, (2020).
 
 
 \bibitem{quastel2015one}
{\sc J.~Quastel and H.~Spohn}, {\em The one-dimensional {KPZ} equation and its
  universality class}, Journal of Statistical Physics, 160 (2015),
  pp.~965--984.
  
  \bibitem{rosati}
  {\sc
  T.~Rosati}, {\em Synchronization for {KPZ}}, arXiv preprint arXiv:1907.06278 (2019).
  
  
  \bibitem{spohn}
  {\sc T.~Sasamoto and H.~Spohn},  {\em Exact height distributions for the KPZ equation with narrow wedge initial condition},  Nuclear Physics B 834.3 (2010), pp.~523--542.
  
  \bibitem{spohn1}
  {\sc T.~Sasamoto and H.~Spohn}, {\em One-dimensional Kardar-Parisi-Zhang equation: an exact solution and its universality}, Physical review letters 104.23 (2010): 230602.
  
  \bibitem{sinai}
  {\sc
  Y. G.~Sinai}, {\em Two results concerning asymptotic behavior of solutions of the Burgers equation with force}, Journal of statistical physics, 64(1), (1991), pp.~1--12.
  
  \bibitem{virag}
  {\sc B.~Vir\'ag}, {\em The heat and the landscape I}, arXiv preprint arXiv: 2008.07241 (2020).
  
  
  \bibitem{walsh}
  {\sc
  J.~Walsh}, {\em An introduction to stochastic partial differential equations}, \'Ecole d'\'Et\'e de Probabilit\'es de Saint Flour XIV-1984. Springer, Berlin, Heidelberg (1986), pp.~265--439.

\end{thebibliography}

\end{document}